# A new analytical formula for the inverse of a square matrix

W. Astar
University of Maryland, Baltimore County(UMBC); Baltimore, Maryland 21250

***Abstract -*** A concise analytical formula is developed for the inverse of an invertible 3 x 3 matrix using a telescoping method, and is generalized to larger square matrices. The formula is confirmed using randomly generated matrices in Matlab.

## 1. Introduction

The inverse of a square matrix can be found using a number of well-known techniques, such as Gaussian elimination, Gauss-Jordan elimination, and LU-decomposition [1], which are among the many modules available in Mathematica®, Matlab®, Maple® and Maxima, among others [2]. The focus of this report, however, is to explore a closed-form analytical formula for the inverse of a square matrix, and based on standard functions.

The product of a complex scalar *a* with its inverse,

$$\frac{1}{a}a = a\frac{1}{a} = 1, \qquad (1.1)$$

can be extended to an invertible $N \times N$ matrix **A** of such scalars,

$$\mathbf{A}^{-1}\mathbf{A} = \mathbf{A}\mathbf{A}^{-1} = \mathbf{I} \qquad (1.2)$$

for which **I** is the identity matrix, an $N \times N$ matrix whose scalar analog is unity. Like that for the complex scalar *a* (1.1), the product (1.2) is also commutative for the matrix **A**,

$$\mathbf{A} = \begin{bmatrix} a_{11} & a_{12} & \cdots & a_{1N} \\ a_{21} & a_{22} & \cdots & a_{2N} \\ \vdots & \vdots & \ddots & \vdots \\ a_{N1} & a_{N2} & \cdots & a_{NN} \end{bmatrix}. \qquad (1.3)$$

In short-form notation, the matrix is often expressed as

$$\mathbf{A} \sim a_{i,j} = a_{ij}, \quad \{i,j\} \in \{\{1,2,3,\cdots,N\} \times \{1,2,3,\cdots,N\}\} \qquad (1.4)$$

with the tilde '~' denoting equivalence. A comma-delimiter is sometimes used to clarify the dependence of the matrix element *a* on the row index (*i*) and the column index (*j*), which is especially helpful when the indices themselves happen to be functions. Both notations will be used interchangeably throughout this report. Although the matrix elements $a_{ij}$ are generally complex in this report, they are not expressed explicitly in terms of their respective real and imaginary parts, and the letters $i$ and $j$ are strictly reserved for the row and the column indices. The inverse of **A** is found from

$$\mathbf{A}^{-1} = \frac{\tilde{\mathbf{A}}}{|\mathbf{A}|}, \tag{1.5}$$

which is the ratio of the adjugate (the adjoint) of **A** to the determinant of **A**. The adjugate of **A** is obtained from the transpose of the matrix of its cofactors $C_{ij}$. A cofactor of a matrix element $a_{ij}$ is obtained from the product of its corresponding minor $M_{ij}$ with the sign-factor $(-1)^{i+j}$. A minor $M_{ij}$ is found by taking the determinant of the minor-matrix **M**,

$$M_{ij} = |\mathbf{M}(i,j)|. \tag{1.6}$$

Since an $N \times N$ matrix **A** possesses $N^2$ elements, it must also possess $N^2$ minors or minor-matrices. More explicitly, the minor-matrix is derived in terms of the original matrix as

$$\mathbf{M}(i,j) = \mathbf{A}(1,\cdots,i-2,i-1,i+1,i+2,\cdots,N;\ 1,\cdots,j-2,j-1,j+1,j+2,\cdots,N). \tag{1.7}$$

The implication is that **M** is a sub-matrix of **A**, and is found from **A** by omitting its $i$-th row and $j$-th column. It should be emphasized here that the elements $m$ of **M** are not indexed to $i$ and $j$ which are given in (1.4), but to $r$ and $s$. However, **M** is a function of the indices $i$ and $j$, whence the nomenclature of $\mathbf{M}(i,j)$. The indices of the matrix elements $m_{rs}$ that make up a minor-matrix **M**, are found from the Cartesian product of the set of row indices with the set of column indices,

$$\{r,s\} \in \{\{1,2,3,\cdots,N-1\} \times \{1,2,3,\cdots,N-1\}\}. \tag{1.8}$$

The convention (1.7, 1.8) is followed in subsequent sections of this report, although the index alphabet may differ.

The determinant of **A** may be found from the following weighted cofactor expansion using the 1st row of **A**, also known as a Laplace Expansion [3],

$$|\mathbf{A}| = \sum_{j=1}^{N} a_{1j} C_{1j} = \sum_{j=1}^{N} (-1)^{1+j} a_{1j} |\mathbf{M}(1,j)|. \tag{1.9}$$

It can also be found using any other row, or any other column of the square matrix **A** (1.3). However, (1.9) is the convention to which this report will adhere. This definition may also be adapted to the minor (1.6), which is the determinant of the minor-matrix **M**.

The inverse matrix (1.5) may also be expressed in terms of its minors (1.6) as

$$\mathbf{A}^{-1} \sim \alpha_{ji} = \frac{(-1)^{i+j} M_{ij}}{|\mathbf{A}|} = \frac{(-1)^{i+j} M_{ij}}{\sum_{j=1}^{N} a_{1j}(-1)^{1+j} M_{1j}}\ ;\ \{i,j\},\{j,i\} \in \{\{1,\cdots,N\} \times \{1,\cdots,N\}\}. \tag{1.10}$$

The denominator, which represents the determinant of **A**, is only computed once, from the sum of the minors corresponding to the elements of the 1st (*i*=1) row of **A**, weighted by those elements. This is perhaps the most general form of a formula for an element of the inverse of a matrix. However, (1.10) is clearly not entirely expressed in terms of the elements of the original matrix **A**. The kernel of the inverse matrix (1.10) is the minor (1.6), which appears in both the numerator and denominator. The minor itself is actually the determinant of the matrix (1.7). Thus, the problem to resolve, is reduced to finding a general analytical formula for the determinant of a (*N*-1) x (*N*-1) matrix **M**.

Apparently, the term "determinant" was first used by C.F. Gauss in 1801 [4], but the modern meaning of the term was not introduced until 1841 by A. Cayley. The Cayley-Hamilton Theorem [5, 6], which was later generalized by Frobenius [7], can be used to find the determinant of a square matrix. The concept of a determinant also appears to have been in use well before the 19th century, beginning with G. Leibniz (1678), and G. Cramer (1750). The latter was apparently the first to publish on the concept, and with whom is associated Cramer's Rule, which still is under active investigation [8, 9]. P.F. Sarrus, developed a graphical mnemonic for the determinant of a 3 x 3 matrix, now known as the Rule of Sarrus. Although originally constrained to 3 x 3 matrices, it has been recently modified, possibly improved [10-12], and also extended to larger 4 x 4 square matrices [13-16].

In 1866, C.L. Dodgson (L. Carroll) [17-19] developed an algorithm termed the Dodgson Condensation, by which the determinant is found by constructing progressively smaller matrices from the original matrix, concluding with a 1 x 1 matrix that carries the determinant of the original matrix itself. There have recently been various mathematical proofs of the algorithm [20-22], and its connection to Chiò's Pivotal Condensation Method (1853) [23], among others, has also been investigated [24, 25]. Among recently introduced recursive methods, has been that due to Rezaifar and Rezaee [26], which has also been generalized [27], but purportedly found to be a consequence of an identity due to the 19th century mathematician A. Capelli [28, 29].

It appears that the most general analytical formula for the determinant of a square matrix was first proposed by Leibniz, one form of which is given by [30-32]

$$|\mathbf{A}| = \sum_{\sigma^{(j)}=\sigma^{(1)}}^{\sigma^{(N!)}} \text{sgn}\left(\sigma^{(j)}\right) \prod_{i=1}^{N} a_{i,\sigma_i^{(j)}} = \sum_{\sigma^{(j)} \in S_N} \text{sgn}\left(\sigma^{(j)}\right) \prod_{i=1}^{N} a_{i,\sigma_i^{(j)}}. \qquad (1.11)$$

Each set $\sigma^{(j)}$ is identical to the *ordered* set *J* of all possible *N* column index values, but is a different (*j*-th) permutation of *J*. The sign function sgn(•) yields positive (negative) unity when the set $\sigma^{(j)}$ is obtained from an even (odd) number of permutations of the set *J*. Unlike the row index *i*, the column index *j* is not explicitly used as such in the equation. Moreover, the elements of the permutated column index set $\sigma^{(j)}$ make no appearance in the formula, and this set gets larger for larger square matrices. Indeed, this set has to be identified prior to its usage in (1.11), since the identification process is not explicitly accounted for in (1.11). Thus, although (1.11) is an acceptable formula, it is not formulaic to the extent of an explicit dependence on the independent variables[1], represented by the row and column indices of the elements $a_{ij}$ (1.4).

---

[1] examples of which, are Newton's second law $F = m \cdot a$, and Einstein's mass-energy equivalence $E = m \cdot c^2$

The aim of this report is to derive an analytical formula for any element of the inverse of an invertible $N \times N$ matrix in terms of the elements of that matrix, using standard functions, demonstrating explicit dependence on the row and column indices, but with nothing known *a priori* about the values of the elements $a_{ij}$. No algebraic relation connecting any of the elements is assumed at the outset, unlike in [33, 34]. The report intends to explore the limitations of such an approach, and to also identify the mathematical hurdles to surmount in order to generalize it to arbitrary square matrices. It will be found that the approach is a telescoping method, like the Dodgson Condensation, but has more in common with the Leibniz formula for determinants (1.11). The method will be demonstrated for $N$ up to 5, after which it will be generalized for any $N$.

## 2. The 2 x 2 matrix D

The 2 x 2 matrix is generally given by

$$\mathbf{D} \sim d_{ij} = \begin{bmatrix} d_{11} & d_{12} \\ d_{21} & d_{22} \end{bmatrix}, \quad \{i,j\} \in \{\{1,2\} \times \{1,2\}\} \tag{2.1}$$

and is well known to have the inverse matrix

$$\mathbf{D}^{-1} = \frac{\begin{bmatrix} d_{22} & -d_{12} \\ -d_{21} & d_{11} \end{bmatrix}}{d_{11}d_{22} - d_{12}d_{21}}. \tag{2.2}$$

Relative to the original matrix **D**, the diagonal elements are interchanged, whereas the skew-diagonal elements are only negated. It has the formulaic form of

$$\mathbf{D}^{-1} \sim \delta_{ji} = \frac{(-1)^{i+j} d_{3-i,3-j}}{\sum_{j=1}^{2} \delta_{i-1} d_{ij} (-1)^{i+j} d_{3-i,3-j}} = \frac{(-1)^{i+j} d_{3-i,3-j}}{\sum_{j=1}^{2} (-1)^{1+j} d_{1,j} d_{2,3-j}}, \quad \{j,i\} \in \{\{1,2\} \times \{1,2\}\}, \tag{2.3}$$

with an expression for the determinant of **D** found in the denominator,

$$|\mathbf{D}| = \sum_{j=1}^{2} \delta_{i-1} d_{ij} (-1)^{i+j} d_{3-i,3-j} = \sum_{j=1}^{2} (-1)^{1+j} d_{1,j} d_{2,3-j} = d_{11}d_{22} - d_{12}d_{21}. \tag{2.4}$$

The formula (2.3)[2] reproduces the indicial change required for the diagonal elements of **D**, while maintaining the indices of the skew-diagonal elements intact. The formula also negates the skew-diagonal elements. This formula is not unique, and other forms are possible. The point of deriving the formulas in this section, serves only to illustrate an approach, which is to be adapted to the larger $N \times N$ matrices, and which are more complicated due to the fact that the inverse of an $N \times N$ matrix involves the evaluations of $(N-1) \times (N-1)$ minors, $N$ times.

---

[2] **The matrix element $\delta_{ij}$, should not be confused with the Kronecker delta-function δ which is not italicized, and which is used in the denominator of (2.3), in (2.4), and in subsequent sections**

## 3. The 3x3 matrix C

The 3 x 3 matrix is generally given by

$$\mathbf{C} \sim c_{kl} \sim \begin{bmatrix} c_{11} & c_{12} & c_{13} \\ c_{21} & c_{22} & c_{23} \\ c_{31} & c_{32} & c_{33} \end{bmatrix}, \quad \{k,l\} \in \{\{1,2,3\} \times \{1,2,3\}\}. \tag{3.1}$$

Unlike that of the 2 x 2 matrix D, the inverse of this matrix is not so easy to commit to memory, but is still analytical and compact, and is given by applying (1.6, 7) 9 times in (1.10) with $\mathbf{A} = \mathbf{C}$,

$$\mathbf{C}^{-1} \sim \gamma_{lk} \sim \frac{\begin{bmatrix} c_{22}c_{33} - c_{23}c_{32} & -(c_{21}c_{33} - c_{23}c_{31}) & c_{21}c_{32} - c_{22}c_{31} \\ -(c_{12}c_{33} - c_{13}c_{32}) & c_{11}c_{33} - c_{13}c_{31} & -(c_{11}c_{32} - c_{12}c_{31}) \\ c_{12}c_{23} - c_{13}c_{22} & -(c_{11}c_{23} - c_{13}c_{21}) & c_{11}c_{22} - c_{12}c_{21} \end{bmatrix}}{c_{11}(c_{22}c_{33} - c_{23}c_{32}) - c_{12}(c_{21}c_{33} - c_{23}c_{31}) + c_{13}(c_{21}c_{32} - c_{22}c_{31})}. \tag{3.2}$$

The adjugate of $\mathbf{C}$ is actually the transpose of the numerator of (3.2), which is represented by $\gamma_{kl}|\mathbf{C}|$. However, the transpose operation is omitted due to the use of $\gamma_{lk}$ instead of $\gamma_{kl}$ on the LHS, since an element of a transposed matrix is identical to the element of the original matrix, if its indices are interchanged. A general, compact formula for any element $\gamma_{lk}$ of the inverse of the matrix $\mathbf{C}$ will now be found.

Up to this point, the derivation of the formula has been based on patterns observed in the indices of the matrix elements and/or its cofactors. A more systematic approach will be developed for a 3 x 3 matrix, that can be adapted to larger square matrices. It is based on the identification of the general form of the 2 x 2 minor-matrix that leads to the evaluation of all the minors required for the inverse of the matrix (3.1).

For the 1st row of $\mathbf{C}$ (3.1), the minor-matrices are

$$\begin{bmatrix} \not{c}_{11} & \not{c}_{12} & \not{c}_{13} \\ \not{c}_{21} & c_{22} & c_{23} \\ \not{c}_{31} & c_{32} & c_{33} \end{bmatrix} \to \begin{bmatrix} c_{22} & c_{23} \\ c_{32} & c_{33} \end{bmatrix}; \begin{bmatrix} \not{c}_{11} & \not{c}_{12} & \not{c}_{13} \\ c_{21} & \not{c}_{22} & c_{23} \\ c_{31} & \not{c}_{32} & c_{33} \end{bmatrix} \to \begin{bmatrix} c_{21} & c_{23} \\ c_{31} & c_{33} \end{bmatrix}; \begin{bmatrix} \not{c}_{11} & \not{c}_{12} & \not{c}_{13} \\ c_{21} & c_{22} & \not{c}_{23} \\ c_{31} & c_{32} & \not{c}_{33} \end{bmatrix} \to \begin{bmatrix} c_{21} & c_{22} \\ c_{31} & c_{32} \end{bmatrix}$$

(3.3)

and using the Kronecker delta-function, defined as

$$\delta_{s-s_0} = \delta_{s,s_0} = \delta[s-s_0] = \begin{cases} 1, & s = s_0 \\ 0, & s \neq s_0 \end{cases}; \quad s, s_0 \in \mathbb{Z} \tag{3.4}$$

yields the *general* minor-matrix, a 2 x 2 matrix corresponding to any element on the 1st row of (3.1), in accordance with the scheme shown by (3.3),

$$\mathbf{D_C}(1,l) = \mathbf{C}(2,3; 1+\delta_{l-1}, 2+\delta_{l-1}+\delta_{l-2}) = \begin{bmatrix} c_{2,1+\delta_{l-1}} & c_{2,2+\delta_{l-1}+\delta_{l-2}} \\ c_{3,1+\delta_{l-1}} & c_{3,2+\delta_{l-1}+\delta_{l-2}} \end{bmatrix}, \; l \in \{1,2,3\}. \quad (3.5)$$

The Kronecker delta-function δ is not to be confused with the (italicized) element *δ* of the inverse (2.3) of the 2 x 2 matrix. On the RHS of this expression, a 2 x 2 sub-matrix is being expressed in terms of the 2nd and 3rd rows, and of the $(1+\delta_{l-1})$-th and the $(2+\delta_{l-1}+\delta_{l-2})$-th columns, of the original 3 x 3 matrix **C** (3.1). Then for the 1st row of **C**, the 3 minors are obtained from (3.5) by setting in turn the column index $l = 1, 2$, then 3. Adapting this procedure for the remaining elements of **C**, yields the general 2 x 2 minor-matrix for the entire matrix **C** (3.1),

$$\mathbf{D_C}(k,l) = \mathbf{C}(1+\delta_{k-1}, 2+\delta_{k-1}+\delta_{k-2}; 1+\delta_{l-1}, 2+\delta_{l-1}+\delta_{l-2}) = \begin{bmatrix} c_{1+\delta_{k-1},1+\delta_{l-1}} & c_{1+\delta_{k-1},2+\delta_{l-1}+\delta_{l-2}} \\ c_{2+\delta_{k-1}+\delta_{k-2},1+\delta_{l-1}} & c_{2+\delta_{k-1}+\delta_{k-2},2+\delta_{l-1}+\delta_{l-2}} \end{bmatrix},$$

$$\{k,l\} \in \{\{1,2,3\} \times \{1,2,3\}\}.$$

(3.6)

The minor-matrix also has the more compact, alternative forms

$$\mathbf{D_C}(k,l) = \mathbf{C}(2-\mathrm{H}_{k-2}, 3-\mathrm{H}_{k-3}; 2-\mathrm{H}_{l-2}, 3-\mathrm{H}_{l-3}), \quad (3.7)$$

$$\mathbf{D_C}(k,l) = \begin{bmatrix} d_{11} & d_{12} \\ d_{21} & d_{22} \end{bmatrix} \sim d_{ij}(k,l) = c_{i'j'} = c_{i+1-\mathrm{H}_{k-i-1}, j+1-\mathrm{H}_{l-j-1}} \sim \begin{bmatrix} c_{2-\mathrm{H}_{k-2},2-\mathrm{H}_{l-2}} & c_{2-\mathrm{H}_{k-2},3-\mathrm{H}_{l-3}} \\ c_{3-\mathrm{H}_{k-3},2-\mathrm{H}_{l-2}} & c_{3-\mathrm{H}_{k-3},3-\mathrm{H}_{l-3}} \end{bmatrix},$$

$$\{i,j\} \in \{\{1,2\} \times \{1,2\}\}, \& \{k,l\} \in \{\{1,2,3\} \times \{1,2,3\}\},$$

(3.8)

using the discrete Heaviside step-function, which, like the Kronecker delta-function whence it is derived, is another generalized function, defined as

$$\mathrm{H}_{s-s_0} = \mathrm{H}[s-s_0] = \sum_{\sigma=s_0}^{\infty} \delta[s-\sigma] = \begin{cases} 1, & s \geq s_0 \\ 0, & s < s_0 \end{cases}, \; \sigma, s, s_0 \in \mathbb{Z}. \quad (3.9)$$

The 2 x 2 minor-matrix $\mathbf{D_C}$ is actually indexed to *i* and *j*, as in the previous section, but its elements $d_{ij}$ are functions of the couple $(k, l)$, specified for the minor of the 3 x 3 matrix (3.1). However, since (3.7, 3.8) is also a 2 x 2 matrix, it should make it amenable to some of the results found in the previous section, as will be seen later. The use of Kronecker delta-functions or Heaviside step-functions in the indices of the matrix elements is a method of incorporating conditional statements in the indices themselves. For instance in (3.6), a row index of $1+\delta_{k-1}$ is interpreted as 1, unless $k = 1$, in which case that row index becomes a 2. Similarly, a column index of $3-\mathrm{H}_{l-3}$ is meant to be construed as a 3, unless $l \geq 3$, in which case the column index is reduced to 2. Lastly, the Kronecker delta-function and the discrete Heaviside step-function are not standard functions, but may be considered to be the discrete analogues of the Dirac delta-function, and the continuous Heaviside step-function, which are generalized functions, or distributions [35-38]. However, it will be shown in §**7** that these functions can be expressed in terms of standard functions over the required values of the indices.

There are now 2 approaches to attaining the inverse of **C**, and they both take advantage of the conclusions reached in the previous section for the 2 x 2 matrix **D**. For the 1st approach, termed the surrogate matrix approach, (1.10) is implemented using the LHS of (3.8) in the form of the matrix $\mathbf{D_C}$ as the general minor-matrix for **C**, yielding

$$\mathbf{C}^{-1} \sim \gamma_{lk} = \frac{(-1)^{k+l}|\mathbf{D_C}(k,l)|}{|\mathbf{C}|} = \frac{(-1)^{k+l}|\mathbf{D_C}(k,l)|}{\sum_{l=1}^{3}\delta_{k-1}c_{kl}(-1)^{k+l}|\mathbf{D_C}(k,l)|}, \quad \{l,k\} \in \{\{1,2,3\} \times \{1,2,3\}\}. \quad (3.10)$$

The equation requires the determinant of the 2 x 2 minor-matrix $\mathbf{D_C}$, whose general form is given by the denominator of the inverse of the 2 x 2 matrix **D** considered in §**2**, specifically (2.4):

$$|\mathbf{D}| = \sum_{j=1}^{2}\delta_{i-1}(-1)^{i+j}d_{ij}d_{3-i,3-j} = d_{11}d_{22} - d_{12}d_{21}. \quad (3.11)$$

Thus, using the surrogate matrix approach, a more explicit version of (3.10) is

$$\mathbf{C}^{-1} \sim \gamma_{lk} = \frac{(-1)^{k+l}\sum_{j=1}^{2}\delta_{i-1}(-1)^{i+j}\left(d_{ij}d_{3-i,3-j}\right)(k,l)}{\sum_{l=1}^{3}\sum_{j=1}^{2}\delta_{i-1}\delta_{k-1}(-1)^{i+j+k+l}c_{kl}\left(d_{ij}d_{3-i,3-j}\right)(k,l)}, \quad \{l,k\} \in \{\{1,2,3\} \times \{1,2,3\}\} \quad (3.12)$$

whose denominator now carries 2 nested Laplace Expansions, one for $|\mathbf{D_C}|$, nested within another for $|\mathbf{C}|$. After carrying out the Kronecker delta-functions where relevant, the equation simplifies to

$$\mathbf{C}^{-1} \sim \gamma_{lk} = \frac{\sum_{j=1}^{2}(-1)^{1+j+k+l}\left(d_{1,j}d_{2,3-j}\right)(k,l)}{\sum_{l=1}^{3}\sum_{j=1}^{2}(-1)^{2+j+l}c_{1,l}\left(d_{1,j}d_{2,3-j}\right)(1,l)}, \quad \{l,k\} \in \{\{1,2,3\} \times \{1,2,3\}\}. \quad (3.13)$$

As stated in §**1**, and based on the convention used in this report, the determinant of the matrix **C** is found from the sum of the minors due to the ($k = 1$) first row of the matrix (3.1), weighted by their corresponding, signed elements $(-1)^{k+l}c_{1l}$ on that row, and is more explicitly given by the denominator of (3.12) or (3.13),

$$|\mathbf{C}| = \sum_{l=1}^{3}\sum_{j=1}^{2}\delta_{i-1}\delta_{k-1}(-1)^{i+j+k+l}c_{kl}\left(d_{ij}d_{3-i,3-j}\right)(k,l) = \sum_{l=1}^{3}\sum_{j=1}^{2}(-1)^{j+l}c_{1,l}\left(d_{1,j}d_{2,3-j}\right)(1,l). \quad (3.14)$$

To obtain a typical element $\gamma_{lk}$ of the inverse matrix $\mathbf{C}^{-1}$, the general minor-matrix $\mathbf{D_C}$ can be found from (3.8) for a specified $\{l, k\}$, after which its elements $d_{ij}$ are substituted into either (3.12) or (3.13). It is tacitly understood that all the elements $d_{ij}$ in the numerator are functions of the indices $\{k, l\}$ due to (3.8), and are expressed as $(k, l)$ in (3.12-14).

In the second, direct approach for formulating the inverse of $\mathbf{C}$, and which encapsulates the first approach, the general minor-matrix found on the RHS of (3.8) is used instead of the LHS, the surrogate matrix $\mathbf{D_C}$. This step is more challenging because the elements now have variable indices, specified as $i'$ and $j'$ in (3.8). For a given $\{k, l\}$, there is a bijective relation between the indices $\{i, j\}$ of the surrogate matrix, and $\{i', j'\}$,

$$i' = i + \sum_{u=1}^{i}\delta_{k-u} = i + 1 - H_{k-i-1} = \kappa(i,k), \tag{3.15}$$

$$j' = j + \sum_{u=1}^{j}\delta_{l-u} = j + 1 - H_{l-j-1} = \kappa(j,l), \tag{3.16}$$

for which $\kappa$ is not to be confused with the index $k$. In the new, primed indicial system $(i', j')$, the standard mathematical operations are not satisfied, examples of which are

$$\begin{aligned}1' + 1' &= 2 + 2\delta_{s-1} \neq 2' = 2 + \delta_{s-1} + \delta_{s-2}, \\ 2' - 1' &= 2 + \delta_{s-1} + \delta_{s-2} - (1 + \delta_{s-1}) = 1 + \delta_{s-2} \neq 1' = 1 + \delta_{s-1}.\end{aligned} \tag{3.17}$$

Eqs. (3.12) and (3.13) are still applicable, but with the elements $d$ replaced by $c$, and with the indices now primed, as follows,

$$\mathbf{C}^{-1} \sim \gamma_{lk} = \frac{\sum_{j=1}^{2}\delta_{i-1}(-1)^{i+j+k+l}c_{i'j'}c_{(3-i)',(3-j)'}}{\sum_{l=1}^{3}\sum_{j=1}^{2}\delta_{i-1}\delta_{k-1}(-1)^{i+j+k+l}c_{kl}c_{i'j'}c_{(3-i)',(3-j)'}}, \quad \{l,k\} \in \{\{1,2,3\} \times \{1,2,3\}\}. \tag{3.18}$$

Based on (3.15),

$$(3-i)' = \begin{cases} 2', & i=1 \\ 1', & i=2 \end{cases} = 4 - i - H_{k-4+i} = \lambda(i,k), \quad i \in \{1,2\}. \tag{3.19}$$

However, the result is also *formally* obtainable from (3.15) by the *discrete convolution*

$$(3-i)' = i' \otimes \delta_{3-i}, \tag{3.20}$$

which basically entails the replacement of $i$ with $(3 - i)$ in (3.15). A similar equation can be found for the corresponding column index $(3-j)'$ using (3.16), by the same logic used for (3.19). However, this equation is also easily found from (3.20) as follows

$$(3-j)' = (3-i)'\delta_{i-j}\delta_{k-l} = 4 - j - H_{l-4+j} = \lambda(j,l), \quad j \in \{1,2\}. \tag{3.21}$$

Substituting (3.15), (3.16), (3.19) and (3.21) into (3.18) yields

$$\mathbf{C}^{-1} \sim \gamma_{lk} = \frac{\sum_{j=1}^{2} \delta_{i-1}(-1)^{i+j+k+l} c_{\kappa(i,k),\kappa(j,l)} c_{\lambda(i,k),\lambda(j,l)}}{\sum_{l=1}^{3}\sum_{j=1}^{2} \delta_{i-1}\delta_{k-1}(-1)^{i+j+k+l} c_{k,l} c_{\kappa(i,k),\kappa(j,l)} c_{\lambda(i,k),\lambda(j,l)}}, \quad \{l,k\} \in \{\{1,2,3\} \times \{1,2,3\}\} \tag{3.22}$$

but more explicitly as,

$$\mathbf{C}^{-1} \sim \gamma_{lk} = \frac{\sum_{j=1}^{2} \delta_{i-1}(-1)^{i+j+k+l} c_{i+1-H_{k-i-1},\, j+1-H_{l-j-1}} \cdot c_{4-i-H_{k-4+i},\, 4-j-H_{l-4+j}}}{\sum_{l=1}^{3}\sum_{j=1}^{2} \delta_{i-1}\delta_{k-1}(-1)^{i+j+k+l} c_{k,l} \cdot c_{i+1-H_{k-i-1},\, j+1-H_{l-j-1}} \cdot c_{4-i-H_{k-4+i},\, 4-j-H_{l-4+j}}}. \tag{3.23}$$

The denominator, which represents the determinant, is only computed once, since it is independent of the specified indices $\{l, k\}$ due to the extant $\delta_{k-1}$ and the summation over $l$. After applying all Kronecker delta-functions, a simpler expression for any element $\gamma_{lk}$ of $\mathbf{C}^{-1}$ is found explicitly in terms of the elements $c_{kl}$ of $\mathbf{C}$ (3.1)[3],

$$\mathbf{C}^{-1} \sim \gamma_{lk} = \frac{\sum_{j=2}^{3}(-1)^{j+k+l} c_{2-H_{k-2},\, j-H_{l-j}} \cdot c_{3-H_{k-3},\, 5-j-H_{l-5+j}}}{\sum_{l=1}^{3}\sum_{j=2}^{3}(-1)^{j+1+l} c_{1,l} \cdot c_{2,\, j-H_{l-j}} \cdot c_{3,\, 5-j-H_{l-5+j}}}, \tag{3.24}$$

after a change of bounds on the $j$-summation. Thus, an expression alternative to (3.14) for the determinant of $\mathbf{C}$, and entirely in terms of its elements, is the denominator of (3.24),

$$|\mathbf{C}| = -\sum_{l=1}^{3}\sum_{j=2}^{3}(-1)^{j+l} c_{1,l} \cdot c_{2,\, j-H_{l-j}} \cdot c_{3,\, 5-j-H_{l-5+j}}. \tag{3.25}$$

Expanding the $j$-summations in the numerator and denominator using (3.15) and (3.16),

$$\mathbf{C}^{-1} \sim \gamma_{lk} = \frac{(-1)^{k+l}\left(c_{1+\delta_{k-1},1+\delta_{l-1}} c_{3-H_{k-3},3-H_{l-3}} - c_{1+\delta_{k-1},3-H_{l-3}} c_{3-H_{k-3},1+\delta_{l-1}}\right)}{\sum_{l=1}^{3}(-1)^{l} c_{1,l}\left(c_{2,3-H_{l-3}} c_{3,1+\delta_{l-1}} - c_{2,1+\delta_{l-1}} c_{3,3-H_{l-3}}\right)}, \quad \{l,k\} \in \{\{1,2,3\} \times \{1,2,3\}\}. \tag{3.26}$$

In §**7**, it is shown that (3.24-26) can be re-cast in terms of standard functions.

---

[3] In Matlab, the summand $c_{2-H_{k-2},\, j-H_{l-j}} \cdot c_{3-H_{k-3},\, 5-j-H_{l-5+j}}$ can be simply expressed as: C(2-(k>=2), j-(l>=j)) * C(3-(k>=3), 5-j-(l>=(5-j)))

## 4. The 4 x 4 matrix B

The 4 x 4 matrix is generally given by[4]

$$\mathbf{B} \sim b_{mn} \sim \begin{bmatrix} b_{11} & b_{12} & b_{13} & b_{14} \\ b_{21} & b_{22} & b_{23} & b_{24} \\ b_{31} & b_{32} & b_{33} & b_{34} \\ b_{41} & b_{42} & b_{43} & b_{44} \end{bmatrix}, \quad \{m,n\} \in \{\{1,2,3,4\} \times \{1,2,3,4\}\}. \tag{4.1}$$

The inverse of this matrix is once again given by applying (1.6, 7) 16 times in (1.10) with $\mathbf{A} = \mathbf{B}$, but this time with unresolved minors,

$$\mathbf{B}^{-1} \sim \beta_{nm} \sim \frac{\begin{bmatrix} \begin{vmatrix} b_{22} & b_{23} & b_{24} \\ b_{32} & b_{33} & b_{34} \\ b_{42} & b_{43} & b_{44} \end{vmatrix} & -\begin{vmatrix} b_{21} & b_{23} & b_{24} \\ b_{31} & b_{33} & b_{34} \\ b_{41} & b_{43} & b_{44} \end{vmatrix} & \begin{vmatrix} b_{21} & b_{22} & b_{24} \\ b_{31} & b_{32} & b_{34} \\ b_{41} & b_{42} & b_{44} \end{vmatrix} & -\begin{vmatrix} b_{21} & b_{22} & b_{23} \\ b_{31} & b_{32} & b_{33} \\ b_{41} & b_{42} & b_{43} \end{vmatrix} \\ -\begin{vmatrix} b_{12} & b_{13} & b_{14} \\ b_{32} & b_{33} & b_{34} \\ b_{42} & b_{43} & b_{44} \end{vmatrix} & \begin{vmatrix} b_{11} & b_{13} & b_{14} \\ b_{31} & b_{33} & b_{34} \\ b_{41} & b_{43} & b_{44} \end{vmatrix} & -\begin{vmatrix} b_{11} & b_{12} & b_{14} \\ b_{31} & b_{32} & b_{34} \\ b_{41} & b_{42} & b_{44} \end{vmatrix} & \begin{vmatrix} b_{11} & b_{12} & b_{13} \\ b_{31} & b_{32} & b_{33} \\ b_{41} & b_{42} & b_{43} \end{vmatrix} \\ \begin{vmatrix} b_{12} & b_{13} & b_{14} \\ b_{22} & b_{23} & b_{24} \\ b_{42} & b_{43} & b_{44} \end{vmatrix} & -\begin{vmatrix} b_{11} & b_{13} & b_{14} \\ b_{21} & b_{23} & b_{24} \\ b_{41} & b_{43} & b_{44} \end{vmatrix} & \begin{vmatrix} b_{11} & b_{12} & b_{14} \\ b_{21} & b_{22} & b_{24} \\ b_{41} & b_{42} & b_{44} \end{vmatrix} & -\begin{vmatrix} b_{11} & b_{12} & b_{13} \\ b_{21} & b_{22} & b_{23} \\ b_{41} & b_{42} & b_{43} \end{vmatrix} \\ -\begin{vmatrix} b_{12} & b_{13} & b_{14} \\ b_{22} & b_{23} & b_{24} \\ b_{32} & b_{33} & b_{34} \end{vmatrix} & \begin{vmatrix} b_{11} & b_{13} & b_{14} \\ b_{21} & b_{23} & b_{24} \\ b_{31} & b_{33} & b_{34} \end{vmatrix} & -\begin{vmatrix} b_{11} & b_{12} & b_{14} \\ b_{21} & b_{22} & b_{24} \\ b_{31} & b_{32} & b_{34} \end{vmatrix} & \begin{vmatrix} b_{11} & b_{12} & b_{13} \\ b_{21} & b_{22} & b_{23} \\ b_{31} & b_{32} & b_{33} \end{vmatrix} \end{bmatrix}}{b_{11}\begin{vmatrix} b_{22} & b_{23} & b_{24} \\ b_{32} & b_{33} & b_{34} \\ b_{42} & b_{43} & b_{44} \end{vmatrix} - b_{12}\begin{vmatrix} b_{21} & b_{23} & b_{24} \\ b_{31} & b_{33} & b_{34} \\ b_{41} & b_{43} & b_{44} \end{vmatrix} + b_{13}\begin{vmatrix} b_{21} & b_{22} & b_{24} \\ b_{31} & b_{32} & b_{34} \\ b_{41} & b_{42} & b_{44} \end{vmatrix} - b_{14}\begin{vmatrix} b_{21} & b_{22} & b_{23} \\ b_{31} & b_{32} & b_{33} \\ b_{41} & b_{42} & b_{43} \end{vmatrix}}. \tag{4.2}$$

The adjugate of $\mathbf{B}$ is actually the transpose of the numerator of (4.2), which is represented by $\beta_{nm}|\mathbf{B}|$, as was the case in the previous section for $\mathbf{C}$. The goal is now to find a general, compact expression for any element $\beta_{nm}$ of the inverse of the matrix $\mathbf{B}$, beginning with its general 3 x 3 minor-matrix, which is now derived.

For the 1st row, the corresponding 3 x 3 minor-matrices (4.3), found on the next page, yield the general 3 x 3 minor-matrix $\mathbf{C_B}$, valid for the 1st row and any $n$-th column

$$\mathbf{C_B}(1,n) = \mathbf{B}\left(2,3,4; 1+\delta_{n-1}, 2+\delta_{n-1}+\delta_{n-2}, 3+\delta_{n-1}+\delta_{n-2}+\delta_{n-3}\right) = \begin{bmatrix} b_{2,1+\delta_{n-1}} & b_{2,2+\delta_{n-1}+\delta_{n-2}} & b_{2,3+\delta_{n-1}+\delta_{n-2}+\delta_{n-3}} \\ b_{3,1+\delta_{n-1}} & b_{3,2+\delta_{n-1}+\delta_{n-2}} & b_{3,3+\delta_{n-1}+\delta_{n-2}+\delta_{n-3}} \\ b_{4,1+\delta_{n-1}} & b_{4,2+\delta_{n-1}+\delta_{n-2}} & b_{4,3+\delta_{n-1}+\delta_{n-2}+\delta_{n-3}} \end{bmatrix}$$
(4.4)

---

[4] The index $m$ is not to be confused with an element of the minor-matrix $\mathbf{M}$, which is not used in this section

$$\begin{bmatrix} b_{11} & b_{12} & b_{13} & b_{14} \\ b_{21} & b_{22} & b_{23} & b_{24} \\ b_{31} & b_{32} & b_{33} & b_{34} \\ b_{41} & b_{42} & b_{43} & b_{44} \end{bmatrix} = \begin{bmatrix} b_{22} & b_{23} & b_{24} \\ b_{32} & b_{33} & b_{34} \\ b_{42} & b_{43} & b_{44} \end{bmatrix}; \quad \begin{bmatrix} b_{11} & b_{12} & b_{13} & b_{14} \\ b_{21} & b_{22} & b_{23} & b_{24} \\ b_{31} & b_{32} & b_{33} & b_{34} \\ b_{41} & b_{42} & b_{43} & b_{44} \end{bmatrix} = \begin{bmatrix} b_{21} & b_{23} & b_{24} \\ b_{31} & b_{33} & b_{34} \\ b_{41} & b_{43} & b_{44} \end{bmatrix};$$

$$\begin{bmatrix} b_{11} & b_{12} & b_{13} & b_{14} \\ b_{21} & b_{22} & b_{23} & b_{24} \\ b_{31} & b_{32} & b_{33} & b_{34} \\ b_{41} & b_{42} & b_{43} & b_{44} \end{bmatrix} = \begin{bmatrix} b_{21} & b_{22} & b_{24} \\ b_{31} & b_{32} & b_{34} \\ b_{41} & b_{42} & b_{44} \end{bmatrix}; \quad \begin{bmatrix} b_{11} & b_{12} & b_{13} & b_{14} \\ b_{21} & b_{22} & b_{23} & b_{24} \\ b_{31} & b_{32} & b_{33} & b_{34} \\ b_{41} & b_{42} & b_{43} & b_{44} \end{bmatrix} = \begin{bmatrix} b_{21} & b_{22} & b_{23} \\ b_{31} & b_{32} & b_{33} \\ b_{41} & b_{42} & b_{43} \end{bmatrix}.$$

(4.3)

which always involves the 2nd, 3rd, and 4th rows of (4.1). When $n = 1$, the minor-matrix only involves the 2nd, 3rd, and 4th columns, for $n = 2$, the 1st, 3rd, and 4th columns, and lastly, for $n = 3$, the 1st, 2nd, and 4th columns. Carrying out this process for the remaining 12 elements of (4.1) yields the general minor-matrix $\mathbf{C_B}$ for the matrix $\mathbf{B}$:

$$\mathbf{C_B}(m,n) = \mathbf{B}\left(1+\delta_{m-1}, 2+\delta_{m-1}+\delta_{m-2}, 3+\delta_{m-1}+\delta_{m-2}+\delta_{m-3}; 1+\delta_{n-1}, 2+\delta_{n-1}+\delta_{n-2}, 3+\delta_{n-1}+\delta_{n-2}+\delta_{n-3}\right)$$

(4.5)

or in long-form as

$$\mathbf{C_B}(m,n) \sim c_{kl}(m,n) \sim \begin{bmatrix} b_{1+\delta_{m-1},1+\delta_{n-1}} & b_{1+\delta_{m-1},2+\delta_{n-1}+\delta_{n-2}} & b_{1+\delta_{m-1},3+\delta_{n-1}+\delta_{n-2}+\delta_{n-3}} \\ b_{2+\delta_{m-1}+\delta_{m-2},1+\delta_{n-1}} & b_{2+\delta_{m-1}+\delta_{m-2},2+\delta_{n-1}+\delta_{n-2}} & b_{2+\delta_{m-1}+\delta_{m-2},3+\delta_{n-1}+\delta_{n-2}+\delta_{n-3}} \\ b_{3+\delta_{m-1}+\delta_{m-2}+\delta_{m-3},1+\delta_{n-1}} & b_{3+\delta_{m-1}+\delta_{m-2}+\delta_{m-3},2+\delta_{n-1}+\delta_{n-2}} & b_{3+\delta_{m-1}+\delta_{m-2}+\delta_{m-3},3+\delta_{n-1}+\delta_{n-2}+\delta_{n-3}} \end{bmatrix},$$

$$\{k,l\} \in \{\{1,2,3\} \times \{1,2,3\}\}, \quad \{m,n\} \in \{\{1,2,3,4\} \times \{1,2,3,4\}\}.$$

(4.6)

Since $\mathbf{C_B}$ is a 3 x 3 matrix like that analyzed in the previous section, its elements $c$ are indexed using $k$ and $l$, although these elements are still functions of $(m, n)$, as evidenced by the RHS of (4.6). Consequently, a representative element of $\mathbf{C_B}$ is expressed as $c_{kl}(m,n)$. Unlike the indices of the elements of the RHS of (4.5), or (4.6), the elements $c_{kl}$ are always indexed beginning with $\{k, l\} = \{1,1\}$, and are limited to a maximum of $\{3, 3\}$. The 2 expressions can be simplified to the more compact forms

$$\mathbf{C_B}(m,n) = \mathbf{B}(2-\mathrm{H}_{m-2}, 3-\mathrm{H}_{m-3}, 4-\mathrm{H}_{m-4}; 2-\mathrm{H}_{n-2}, 3-\mathrm{H}_{n-3}, 4-\mathrm{H}_{n-4}),$$  (4.7)

$$\mathbf{C_B}(m,n) \sim c_{kl}(m,n) = b_{k+1-\mathrm{H}_{m-k-1}, l+1-\mathrm{H}_{n-l-1}} \sim \begin{bmatrix} b_{2-\mathrm{H}_{m-2}, 2-\mathrm{H}_{n-2}} & b_{2-\mathrm{H}_{m-2}, 3-\mathrm{H}_{n-3}} & b_{2-\mathrm{H}_{m-2}, 4-\mathrm{H}_{n-4}} \\ b_{3-\mathrm{H}_{m-3}, 2-\mathrm{H}_{n-2}} & b_{3-\mathrm{H}_{m-3}, 3-\mathrm{H}_{n-3}} & b_{3-\mathrm{H}_{m-3}, 4-\mathrm{H}_{n-4}} \\ b_{4-\mathrm{H}_{m-4}, 2-\mathrm{H}_{n-2}} & b_{4-\mathrm{H}_{m-4}, 3-\mathrm{H}_{n-3}} & b_{4-\mathrm{H}_{m-4}, 4-\mathrm{H}_{n-4}} \end{bmatrix}.$$

(4.8)

For any given couple $\{m, n\}$, the corresponding minor-matrix $\mathbf{C_B}$ can be found from (4.8). It is possible to telescope to a 2 x 2 matrix $\mathbf{D_C}$ that serves as the minor-matrix, to the minor-matrix $\mathbf{C_B}$ itself, in short-form as

$$\mathbf{D_C}(k,l) = \mathbf{C_B}\left(1+\delta_{k-1}, 2+\delta_{k-1}+\delta_{k-2}; 1+\delta_{l-1}, 2+\delta_{l-1}+\delta_{l-2}\right), \tag{4.9}$$

and in long-form,

$$\mathbf{D_C}(k,l) \sim d_{ij}(k,l) = c_{i+1-H_{k-i-1},\, j+1-H_{l-j-1}} = \begin{bmatrix} c_{2-H_{k-2},\, 2-H_{l-2}} & c_{2-H_{k-2},\, 3-H_{l-3}} \\ c_{3-H_{k-3},\, 2-H_{l-2}} & c_{3-H_{k-3},\, 3-H_{l-3}} \end{bmatrix}, \tag{4.10}$$

$$\{i, j\} \in \{\{1,2\} \times \{1,2\}\}, \& \{k,l\} \in \{\{1,2,3\} \times \{1,2,3\}\}.$$

As for the previous section, there are 2 approaches to attaining the inverse of $\mathbf{B}$, and they both take advantage of the conclusions reached in the previous section for the matrix $\mathbf{C}$.

For the first, surrogate matrix approach, the inverse matrix is found by substituting (4.7) into (1.10) with $\mathbf{A} = \mathbf{B}$, and with $\mathbf{M} = \mathbf{C_B}$,

$$\mathbf{B}^{-1} \sim \beta_{nm} = \frac{(-1)^{m+n}|\mathbf{C_B}(m,n)|}{|\mathbf{B}|} = \frac{(-1)^{m+n}|\mathbf{C_B}(m,n)|}{\displaystyle\sum_{n=1}^{4}\delta_{m-1}b_{mn}(-1)^{m+n}|\mathbf{C_B}(m,n)|}, \{n,m\} \in \{\{1,2,3,4\} \times \{1,2,3,4\}\}.$$

(4.11)

The determinant of a 3 x 3 matrix $\mathbf{C_B}$ was already found in §3, and is given by (3.14),

$$|\mathbf{C_B}| = \sum_{l=1}^{3}\delta_{k-1}c_{kl}(-1)^{k+l}|\mathbf{D_C}(k,l)| = \sum_{l=1}^{3}\sum_{j=1}^{2}\delta_{i-1}\delta_{k-1}(-1)^{i+j+k+l}c_{kl}(m,n)\left(d_{ij}d_{3-i,3-j}\right)(k,l). \tag{4.12}$$

In one form, the inverse (4.11) may then be re-cast using the LHS of (4.12),

$$\mathbf{B}^{-1} \sim \beta_{nm} = \frac{(-1)^{m+n}|\mathbf{C_B}(m,n)|}{|\mathbf{B}|} = \frac{\displaystyle\sum_{l=1}^{3}\delta_{k-1}(-1)^{k+l+m+n}c_{kl}|\mathbf{D_C}(k,l)|}{\displaystyle\sum_{n=1}^{4}\sum_{l=1}^{3}\delta_{k-1}\delta_{m-1}(-1)^{k+l+m+n}b_{mn}c_{kl}|\mathbf{D_C}(k,l)|} \tag{4.13}$$

which is an expression of the inverse of $\mathbf{B}$ in terms of the determinant of the 2 x 2 matrix $\mathbf{D_C}$. Alternatively, substituting the RHS of (4.12) into (4.11), $\mathbf{B}^{-1}$ is more explicitly

$$\mathbf{B}^{-1} \sim \beta_{nm} = \frac{(-1)^{m+n}\displaystyle\sum_{l=1}^{3}\sum_{j=1}^{2}\delta_{i-1}\delta_{k-1}(-1)^{i+j+k+l}c_{kl}(m,n)\cdot\left(d_{ij}d_{3-i,3-j}\right)(k,l)}{\displaystyle\sum_{n=1}^{4}\sum_{l=1}^{3}\sum_{j=1}^{2}\delta_{i-1}\delta_{k-1}\delta_{m-1}(-1)^{i+j+k+l+m+n}b_{mn}\cdot c_{kl}(m,n)\cdot\left(d_{ij}d_{3-i,3-j}\right)(k,l)}, \tag{4.14}$$

$$\{n,m\} \in \{\{1,2,3,4\} \times \{1,2,3,4\}\}.$$

It should be clarified that the elements $c$ of the surrogate minor-matrix $\mathbf{C_B}$ are indexed to $\{k, l\}$, but are functions of the indices $m$ and $n$ of the original 4 × 4 matrix $\mathbf{B}$ (4.1), and are thus represented as $c_{kl}(m, n)$, in accordance with (4.8). Furthermore, the elements $d$ of the surrogate 2 × 2 matrix $\mathbf{D_C}$ are indexed to $\{i, j\}$ but are functions of the indices $k$ and $l$ of the minor-matrix $\mathbf{C_B}$ and so are represented as $d_{ij}(k, l)$, in accordance with (4.10). The indices of a matrix element are always expressed as subscripts in this report. The denominator, which is the determinant of $\mathbf{B}$, is effectively 3 nested Laplace Expansions. On the other hand, the numerator, which is the determinant of the minor-matrix $\mathbf{C_B}$, is comprised of 2 nested Laplace Expansions. Applying the various Kronecker delta-functions yields the simplified expression,

$$\mathbf{B}^{-1} \sim \beta_{nm} = \frac{(-1)^{m+n} \sum_{l=1}^{3} \sum_{j=1}^{2} (-1)^{2+j+l} c_{1,l}(m,n) \cdot \left(d_{1,j} d_{2,3-j}\right)(1,l)}{\sum_{n=1}^{4} \sum_{l=1}^{3} \sum_{j=1}^{2} (-1)^{3+j+l+n} b_{1,n} \cdot c_{1,l}(1,n) \cdot \left(d_{1,j} d_{2,3-j}\right)(1,l)}, \quad (4.15)$$

$$\{n, m\} \in \{\{1, 2, 3, 4\} \times \{1, 2, 3, 4\}\},$$

which, from the denominator, yields the determinant

$$|\mathbf{B}| = -\sum_{n=1}^{4} \sum_{l=1}^{3} \sum_{j=1}^{2} (-1)^{j+l+n} b_{1,n} \cdot c_{1,l}(1,n) \cdot \left(d_{1,j} d_{2,3-j}\right)(1,l). \quad (4.16)$$

In the second, direct approach for formulating the inverse of $\mathbf{B}$, the elements $d$ of the surrogate matrix $\mathbf{D_C}$ are replaced by the elements $b$ of the original matrix $\mathbf{B}$. This is possible since the matrix $\mathbf{D_C}$ is obtained from the matrix $\mathbf{C_B}$, which in turn, is obtained from the original matrix $\mathbf{B}$. Similarly, the elements $c$ are also replaced by $b$'s. As in the previous section, this approach is significantly more challenging, since both $\mathbf{C_B}$ and $\mathbf{D_C}$ are variable in indices, as is evident from (4.8) and (4.10). After taking the 1st step, in which both elements $c$ and $d$ are replaced by $b$'s in (4.14),

$$\mathbf{B}^{-1} \sim \beta_{nm} = \frac{(-1)^{m+n} \sum_{l=1}^{3} \sum_{j=1}^{2} \delta_{i-1} \delta_{k-1} (-1)^{i+j+k+l} b_{k'l'}(m,n) \cdot \left(b_{i''j''} b_{(3-i)'',(3-j)''}\right)(k,l)}{\sum_{n=1}^{4} \sum_{l=1}^{3} \sum_{j=1}^{2} \delta_{i-1} \delta_{k-1} \delta_{m-1} (-1)^{i+j+k+l+m+n} b_{mn} \cdot b_{k'l'}(m,n) \cdot \left(b_{i''j''} b_{(3-i)'',(3-j)''}\right)(k,l)},$$

$$\{n, m\} \in \{\{1, 2, 3, 4\} \times \{1, 2, 3, 4\}\}.$$
(4.17)

A new convention is now adopted for the indices of the elements $b$ that replace both $c$ and $d$ in (4.14): The elements $b$ retain the indices of the elements being replaced, but they are primed a number of times dependent on the number of telescoping steps associated with the matrix whose element is being replaced by $b$. After replacing the elements $c$ with $b$ in (4.14), its indices now carry single primes since the replaced elements $c$ represent $\mathbf{C_B}$, which is found from the original matrix $\mathbf{B}$ after a single telescoping step. Following the approach used for (3.15, 16), the indices now have the following expressions,

$$k' = k+1-H_{m-k-1} = \kappa(k,m), \tag{4.18}$$

$$l' = l+1-H_{n-l-1} = \kappa(l,n), \tag{4.19}$$

which are obtained from the LHS of (4.8). The Latin index $k$ should not be confused with the function $\kappa$, which is Greek. Upon replacing the elements $d$ in (4.14) with $b$, its indices acquire double-primes, since the replaced elements $d$ are obtained from $\mathbf{D_C}$, which is found from the minor-matrix $\mathbf{C_B}$, which in turn is found from $\mathbf{B}$ - a total of two telescoping steps. There results a total of 4 double-primed indices that have to be resolved. However, all the double-primed indices may be obtainable through simple operations on just one of the indices. The row index $i''$ is arbitrarily considered first:

$$i'' = (i+1-H_{k-i-1})' = \begin{cases} (2-H_{k-2})', i=1 \\ (3-H_{k-3})', i=2 \end{cases} = \begin{cases} 2'H_{1-k} + 1'H_{k-2}, i=1 \\ 3'H_{2-k} + 2'H_{k-3}, i=2 \end{cases} = \begin{cases} (2-H_{m-2})H_{k-2} + (3-H_{m-3})H_{1-k}, i=1 \\ (3-H_{m-3})H_{k-3} + (4-H_{m-4})H_{2-k}, i=2 \end{cases} \tag{4.20}$$

which can be summarized as

$$i'' = (i+1-H_{m-i-1})H_{k-i-1} + (i+2-H_{m-i-2})H_{i-k}; \ i \in \{1,2\}, \ k \in \{1,2,3\}, \ m \in \{1,2,3,4\}. \tag{4.21}$$

It can be simplified to the compact expression

$$i'' = \sum_{u=1}^{2}(i+u-H_{m-i-u})H_{(-1)^u(i-u-k+2)} = \lambda(i,k,m) \tag{4.22}$$

with the same indicial conditions as the ones in (4.21). For the double-primed column index $j''$ in (4.17), it may be obtained after simple operations on (4.21) and/or (4.22),

$$j'' = i''\delta_{i-j}\delta_{k-l}\delta_{m-n} \tag{4.23}$$

with the result

$$j'' = (j+1-H_{n-j-1})H_{l-j-1} + (j+2-H_{n-j-2})H_{j-l} = \sum_{u=1}^{2}(j+u-H_{n-j-u})H_{(-1)^u(j-u-l+2)} = \lambda(j,l,n);$$
$$j \in \{1,2\}, \ l \in \{1,2,3\}, \ n \in \{1,2,3,4\}. \tag{4.24}$$

As for $(3-i)''$, it is obtainable from (4.22) by the discrete convolution

$$(3-i)'' = i'' \otimes \delta_{3-i}, \tag{4.25}$$

which performs a reflection in *i*, followed by a shift of 3 units. It yields the compact expression

$$(3-i)'' = (4-i-H_{m-4+i})H_{k-4+i} + (5-i-H_{m-5+i})H_{3-i-k} = \sum_{u=1}^{2}(-i+u+3-H_{m+i-u-3})H_{(-1)^u(-i-k-u+5)} = \mu(i,k,m);$$

$$i \in \{1,2\}, k \in \{1,2,3\}, m \in \{1,2,3,4\}.$$

(4.26)

The corresponding column index may be obtained from the following simple operations on (4.26):

$$(3-j)'' = (3-i)'' \delta_{i-j} \delta_{k-l} \delta_{m-n} \tag{4.27}$$

with the implication that *i* is to be replaced *j*, *k* by *l*, and *m* by *n*. The result of these operations yields the expression

$$(3-j)'' = (4-j-H_{n-4+j})H_{l-4+j} + (5-j-H_{n-5+j})H_{3-j-l} = \sum_{u=1}^{2}(-j+u+3-H_{n+j-u-3})H_{(-1)^u(-j-l-u+5)} = \mu(j,l,n);$$

$$j \in \{1,2\}, l \in \{1,2,3\}, n \in \{1,2,3,4\}.$$

(4.28)

After substituting (4.18), (4.19), (4.22), (4.24), (4.26) and (4.28) into (4.17),

$$\mathbf{B}^{-1} \sim \beta_{nm} = \frac{(-1)^{m+n} \sum_{l=1}^{3}\sum_{j=1}^{2} \delta_{i-1}\delta_{k-1}(-1)^{i+j+k+l} b_{\kappa(k,m),\kappa(l,n)} \cdot b_{\lambda(i,k,m),\lambda(j,l,n)} \cdot b_{\mu(i,k,m),\mu(j,l,n)}}{\sum_{n=1}^{4}\sum_{l=1}^{3}\sum_{j=1}^{2} \delta_{i-1}\delta_{k-1}\delta_{m-1}(-1)^{i+j+k+l+m+n} b_{mn} \cdot b_{\kappa(k,m),\kappa(l,n)} \cdot b_{\lambda(i,k,m),\lambda(j,l,n)} \cdot b_{\mu(i,k,m),\mu(j,l,n)}},$$

$$\{n,m\} \in \{\{1,2,3,4\} \times \{1,2,3,4\}\}.$$

(4.29)

which is now explicitly in terms of the elements of the original matrix **B**. After enforcing the various Kronecker delta-functions,

$$\mathbf{B}^{-1} \sim \beta_{nm} = \frac{(-1)^{m+n} \sum_{l=1}^{3}\sum_{j=1}^{2}(-1)^{2+j+l} b_{\kappa(1,m),\kappa(l,n)} \cdot b_{\lambda(1,1,m),\lambda(j,l,n)} \cdot b_{\mu(1,1,m),\mu(j,l,n)}}{\sum_{n=1}^{4}\sum_{l=1}^{3}\sum_{j=1}^{2}(-1)^{3+j+l+n} b_{1,n} \cdot b_{\kappa(1,1),\kappa(l,n)} \cdot b_{\lambda(1,1,1),\lambda(j,l,n)} \cdot b_{\mu(1,1,1),\mu(j,l,n)}},$$

$$\{n,m\} \in \{\{1,2,3,4\} \times \{1,2,3,4\}\}$$

(4.30)

which yields the simplified indicial functions

$$\kappa(1,m) = 2 - \mathrm{H}_{m-2} \tag{4.31}$$

$$\kappa(1,1) = 2 - \mathrm{H}_{-1} = 2 \tag{4.32}$$

$$\lambda(1,1,m) = (i+1-\mathrm{H}_{m-i-1})\mathrm{H}_{-1} + (3-\mathrm{H}_{m-3})\mathrm{H}_0 = 3-\mathrm{H}_{m-3}, \tag{4.33}$$

$$\lambda(1,1,1) = 3 \tag{4.34}$$

$$\mu(1,1,m) = (4-i-\mathrm{H}_{m-4+i})\mathrm{H}_{-2} + (4-\mathrm{H}_{m-4})\mathrm{H}_1 = 4-\mathrm{H}_{m-4}, \tag{4.35}$$

$$\mu(1,1,1) = 4 \tag{4.36}$$

since the Heaviside step-function vanishes for negative arguments, according to its definition (3.9). Eqs. (4.33) and (4.35) are actually reduced beginning with the outermost Heaviside step-functions, and once their arguments are found to be negative, they are neglected along with any parenthesized expression with which they are multiplied, since they are always finite. Finally, an expression for the inverse is found explicitly in terms of the elements of the original matrix $\mathbf{B}$:

$$\mathbf{B}^{-1} \sim \beta_{nm} = -\frac{(-1)^{m+n}\sum_{l=1}^{3}\sum_{j=1}^{2}(-1)^{j+l} b_{2-\mathrm{H}_{m-2},\kappa(l,n)} \cdot b_{3-\mathrm{H}_{m-3},\lambda(j,l,n)} \cdot b_{4-\mathrm{H}_{m-4},\mu(j,l,n)}}{\sum_{n=1}^{4}\sum_{l=1}^{3}\sum_{j=1}^{2}(-1)^{j+l+n} b_{1,n} \cdot b_{2,\kappa(l,n)} \cdot b_{3,\lambda(j,l,n)} \cdot b_{4,\mu(j,l,n)}}, \tag{4.37}$$

$$\{n,m\} \in \{\{1,2,3,4\} \times \{1,2,3,4\}\},$$

with the determinant given by the denominator

$$|\mathbf{B}| = \sum_{n=1}^{4}\sum_{l=1}^{3}\sum_{j=1}^{2}(-1)^{j+l+n} b_{1,n} \cdot b_{2,\kappa(l,n)} \cdot b_{3,\lambda(j,l,n)} \cdot b_{4,\mu(j,l,n)} \tag{4.38}$$

and with the indicial functions given by

$$\kappa(l,n) = l+1-\mathrm{H}_{n-l-1}, \tag{4.39}$$

$$\lambda(j,l,n) = \sum_{u=1}^{2}(j+u-\mathrm{H}_{n-j-u})\mathrm{H}_{(-1)^u(j-u-l+2)} = \sum_{u=1}^{2}\kappa(j+u-1,n)\mathrm{H}_{(-1)^u(j-u-l+2)}, \tag{4.40}$$

$$\mu(j,l,n) = \sum_{u=1}^{2}(-j+u+3-\mathrm{H}_{n+j-u-3})\mathrm{H}_{(-1)^u(-j-l-u+5)} = \sum_{u=1}^{2}\kappa(-j+u+2,n)\mathrm{H}_{(-1)^u(-j-u-l+5)}. \tag{4.41}$$

## 5. The 5x5 matrix A

The 5 × 5 matrix is generally given by

$$\mathbf{A} \sim a_{pq} \sim \begin{bmatrix} a_{11} & a_{12} & a_{13} & a_{14} & a_{15} \\ a_{21} & a_{22} & a_{23} & a_{24} & a_{25} \\ a_{31} & a_{32} & a_{33} & a_{34} & a_{35} \\ a_{41} & a_{42} & a_{43} & a_{44} & a_{45} \\ a_{51} & a_{52} & a_{53} & a_{54} & a_{55} \end{bmatrix}, \ \{p,q\} \in \{\{1,2,3,4,5\} \times \{1,2,3,4,5\}\}. \qquad (5.1)$$

The expression for the inverse of this matrix is even larger than that (4.2) of the 4 × 4 **B** matrix (4.1), and is provided here only for completeness:

$$\mathbf{A}^{-1} \sim \alpha_{qp} \sim \frac{\begin{bmatrix} \text{[25 cofactor minors}\ 4\times 4] \end{bmatrix}}{a_{11}\begin{vmatrix}\cdot\end{vmatrix} - a_{12}\begin{vmatrix}\cdot\end{vmatrix} + a_{13}\begin{vmatrix}\cdot\end{vmatrix} - a_{14}\begin{vmatrix}\cdot\end{vmatrix} + a_{15}\begin{vmatrix}\cdot\end{vmatrix}}. \qquad (5.2)$$

The general minor-matrices of the 3 × 3 and 4 × 4 matrices (3.1) and (4.1) were respectively found to be (3.6) and (4.5), reproduced here as

$$\mathbf{C}(1+\delta_{k-1},2+\delta_{k-1}+\delta_{k-2};\ 1+\delta_{l-1},2+\delta_{l-1}+\delta_{l-2}),\ \{k,l\} \in \{\{1,2,3\} \times \{1,2,3\}\}, \qquad (5.3)$$

$$\mathbf{B}(1+\delta_{m-1},2+\delta_{m-1}+\delta_{m-2},3+\delta_{m-1}+\delta_{m-2}+\delta_{m-3};1+\delta_{n-1},2+\delta_{n-1}+\delta_{n-2},3+\delta_{n-1}+\delta_{n-2}+\delta_{n-3}),$$

$$\{m,n\} \in \{\{1,2,3,4\} \times \{1,2,3,4\}\}.$$

$$(5.4)$$

The matrix version for each, can be found from the Cartesian product of the set of row indices, with the set of column indices, defined by the semi-colon delimiter '**;**'. For the 5 x 5 matrix (5.1) then, the following 4 x 4 minor-matrix expression may be deduced by a process of induction,

$$\mathbf{A}\begin{pmatrix} 1+\delta_{p-1}, 2+\delta_{p-1}+\delta_{p-2}, 3+\delta_{p-1}+\delta_{p-2}+\delta_{p-3}, 4+\delta_{p-1}+\delta_{p-2}+\delta_{p-3}+\delta_{p-4}; \\ 1+\delta_{q-1}, 2+\delta_{q-1}+\delta_{q-2}, 3+\delta_{q-1}+\delta_{q-2}+\delta_{q-3}, 4+\delta_{q-1}+\delta_{q-2}+\delta_{q-3}+\delta_{q-4} \end{pmatrix},$$

$$\{p,q\} \in \{\{1,2,3,4,5\} \times \{1,2,3,4,5\}\}.$$

(5.5)

It was also found that the general minor-matrices of the 3 x 3 and 4 x 4 matrices (3.1) and (4.1) had the more compact forms of (3.7) and (4.7), reproduced here as

$$\mathbf{C}(2-H_{k-2}, 3-H_{k-3}; 2-H_{l-2}, 3-H_{l-3}), \quad \{k,l\} \in \{\{1,2,3\} \times \{1,2,3\}\},$$

(5.6)

$$\mathbf{B}(2-H_{m-2}, 3-H_{m-3}, 4-H_{m-4}; 2-H_{n-2}, 3-H_{n-3}, 4-H_{n-4}), \quad \{m,n\} \in \{\{1,2,3,4\} \times \{1,2,3,4\}\},$$

(5.7)

which have been found from (5.3) and (5.4) using the equivalence of

$$\sum_{u=1}^{P}(1+\delta_{s-u}) = P+1-H_{s-P-1} \; ; \; u, P \in \mathbb{Z}.$$

(5.8)

This alternative form, which exclusively uses the Heaviside step-function, is more compact and more amenable to the standard matrix format. Based on (5.6) and (5.7), it is not difficult to see that for the 5 x 5 matrix (5.1), its *general* minor-matrix must be

$$\mathbf{A}(2-H_{p-2}, 3-H_{p-3}, 4-H_{p-4}, 5-H_{p-5}; 2-H_{q-2}, 3-H_{q-3}, 4-H_{q-4}, 5-H_{q-5}),$$

$$\{p,q\} \in \{\{1,2,3,4,5\} \times \{1,2,3,4,5\}\},$$

(5.9)

which is valid for any $\{p, q\}$, and which can also be expressed in matrix form as

$$a_{m',n'} \sim a_{m+1-H_{p-m-1}, n+1-H_{q-n-1}} \sim \begin{bmatrix} a_{2-H_{p-2}, 2-H_{q-2}} & a_{2-H_{p-2}, 3-H_{q-3}} & a_{2-H_{p-2}, 4-H_{q-4}} & a_{2-H_{p-2}, 5-H_{q-5}} \\ a_{3-H_{p-3}, 2-H_{q-2}} & a_{3-H_{p-3}, 3-H_{q-3}} & a_{3-H_{p-3}, 4-H_{q-4}} & a_{3-H_{p-3}, 5-H_{q-5}} \\ a_{4-H_{p-4}, 2-H_{q-2}} & a_{4-H_{p-4}, 3-H_{q-3}} & a_{4-H_{p-4}, 4-H_{q-4}} & a_{4-H_{p-4}, 5-H_{q-5}} \\ a_{5-H_{p-5}, 2-H_{q-2}} & a_{5-H_{p-5}, 3-H_{q-3}} & a_{5-H_{p-5}, 4-H_{q-4}} & a_{5-H_{p-5}, 5-H_{q-5}} \end{bmatrix},$$

$$\{m,n\} \in \{\{1,2,3,4\} \times \{1,2,3,4\}\}, \{p,q\} \in \{\{1,2,3,4,5\} \times \{1,2,3,4,5\}\}.$$

(5.10)

Using the surrogate matrix approach, the following closed-form expressions were found for the inverses of the 2 x 2, the 3 x 3, and the 4 x 4 matrices, respectively given by (2.3), (3.12), and (4.14),

$$\mathbf{D}^{-1} \sim \delta_{ji} = \frac{(-1)^{i+j} d_{3-i,3-j}}{\sum_{j=1}^{2} \delta_{i-1}(-1)^{i+j} d_{ij} d_{3-i,3-j}}, \quad \{j,i\} \in \{\{1,2\} \times \{1,2\}\}, \tag{5.11}$$

$$\mathbf{C}^{-1} \sim \gamma_{lk} = \frac{(-1)^{k+l} \sum_{j=1}^{2} \delta_{i-1}(-1)^{i+j} \left(d_{ij} d_{3-i,3-j}\right)(k,l)}{\sum_{l=1}^{3} \sum_{j=1}^{2} \delta_{i-1} \delta_{k-1}(-1)^{i+j+k+l} c_{kl} \cdot \left(d_{ij} d_{3-i,3-j}\right)(k,l)}, \quad \{l,k\} \in \{\{1,2,3\} \times \{1,2,3\}\}, \tag{5.12}$$

$$\mathbf{B}^{-1} \sim \beta_{nm} = \frac{(-1)^{m+n} \sum_{l=1}^{3} \sum_{j=1}^{2} \delta_{i-1} \delta_{k-1}(-1)^{i+j+k+l} c_{kl}(m,n) \cdot \left(d_{ij} d_{3-i,3-j}\right)(k,l)}{\sum_{n=1}^{4} \sum_{l=1}^{3} \sum_{j=1}^{2} \delta_{i-1} \delta_{k-1} \delta_{m-1}(-1)^{i+j+k+l+m+n} b_{mn} \cdot c_{kl}(m,n) \cdot \left(d_{ij} d_{3-i,3-j}\right)(k,l)}, \tag{5.13}$$

$$\{n,m\} \in \{\{1,2,3,4\} \times \{1,2,3,4\}\}.$$

The expressions can be significantly simplified by enforcing the Kronecker delta-functions within each summand. However, they are all retained to establish the trend that should lead to the desired expression for the inverse of **A**. When expressed as a quotient then, an element of the inverse of a $N$ x $N$ matrix is comprised of ($N$-2) summations in the numerator, and ($N$-1) summations in the denominator, the outermost of which depends on the column index of the $N$ x $N$ matrix. The numerator of the inverse element of a $N$ x $N$ matrix is of the form of the denominator of the inverse element of the ($N$-1) x ($N$-1) matrix. The denominator of the inverse element, which is the matrix determinant, is the sum of its numerator evaluated over the $N$ columns, and weighted by the signed elements of the first row of the matrix. By induction for the 5 x 5 matrix **A** (5.1), any element of its inverse $\mathbf{A}^{-1}$ should be a rational expression with 3 summations in its numerator, and 4 summations in its denominator. Furthermore, its numerator should be of the same form as the denominator of the inverse element (5.13) of the 4 x 4 matrix **B**, signed according to $(-1)^{p+q}$, whereas its denominator should be its numerator summed over all ($q = 5$) five columns of **A**, and weighted by that row's signed elements, $(-1)^{p+q} a_{1q}$. Consequently, the inverse of **A** (5.1) is deduced to have the form of

$$\mathbf{A}^{-1} \sim \alpha_{qp} = \frac{(-1)^{p+q} \sum_{n=1}^{4} \sum_{l=1}^{3} \sum_{j=1}^{2} \delta_{i-1} \delta_{k-1} \delta_{m-1}(-1)^{i+j+k+l+m+n} b_{mn}(p,q) \cdot c_{kl}(m,n) \cdot \left(d_{ij} d_{3-i,3-j}\right)(k,l)}{\sum_{q=1}^{5} \sum_{n=1}^{4} \sum_{l=1}^{3} \sum_{j=1}^{2} \delta_{i-1} \delta_{k-1} \delta_{m-1} \delta_{p-1}(-1)^{i+j+k+l+m+n+p+q} a_{pq} \cdot b_{mn}(p,q) \cdot c_{kl}(m,n) \cdot \left(d_{ij} d_{3-i,3-j}\right)(k,l)},$$

$$\{n,m\} \in \{\{1,2,3,4,5\} \times \{1,2,3,4,5\}\}.$$

(5.14)

In this expression, the elements *b*, *c*, and *d* are respectively obtained from the surrogate matrices **B**$_A$, **C**$_B$, and **D**$_C$. The matrix **B**$_A$ ~ $b_{mn}(p, q)$ is a 4 × 4 matrix that serves as a surrogate for the minor-matrix of **A**, and is a direct function of the indices (*p*, *q*) of the matrix **A**. The matrix **C**$_B$ ~ $c_{kl}(m, n)$ is a 3 × 3 matrix that is used as the surrogate for the minor-matrix of the minor-matrix **B**$_A$, is a direct function of its indices (*m*, *n*), but is also an indirect function of (*p*, *q*). Lastly, **D**$_C$ ~ $d_{ij}(k, l)$ is a 2 × 2 matrix that is used as the surrogate for the minor-matrix of the matrix **C**$_B$. They are all given by

$$\mathbf{B}_A = \mathbf{A}\left(2-H_{p-2}, 3-H_{p-3}, 4-H_{p-4}, 5-H_{p-5}; 2-H_{q-2}, 3-H_{q-3}, 4-H_{q-4}, 5-H_{q-5}\right),$$

$$\{p,q\} \in \{\{1,2,3,4,5\} \times \{1,2,3,4,5\}\},$$

(5.15)

$$\mathbf{C}_B = \mathbf{B}_A\left(2-H_{m-2}, 3-H_{m-3}, 4-H_{m-4}; 2-H_{n-2}, 3-H_{n-3}, 4-H_{n-4}\right),$$

$$\{m,n\} \in \{\{1,2,3,4\} \times \{1,2,3,4\}\},$$

(5.16)

$$\mathbf{D}_C = \mathbf{C}_B\left(2-H_{k-2}, 3-H_{k-3}; 2-H_{l-2}, 3-H_{l-3}\right),$$

$$\{k,l\} \in \{\{1,2,3\} \times \{1,2,3\}\}.$$

(5.17)

More explicitly in terms of the elements of the original matrices themselves, the following expressions for the inverse matrix elements were also found:

$$\mathbf{D}^{-1} \sim \delta_{ji} = \frac{(-1)^{i+j} d_{3-i, 3-j}}{\sum_{j=1}^{2}(-1)^{1+j} d_{1,j} d_{2, 3-j}}, \quad \{j,i\} \in \{\{1,2\} \times \{1,2\}\},$$

(5.18)

$$\mathbf{C}^{-1} \sim \gamma_{lk} = \frac{\sum_{j=1}^{2} \delta_{i-1}(-1)^{i+j+k+l} c_{\kappa(i,k),\kappa(j,l)} c_{\lambda(i,k),\lambda(j,l)}}{\sum_{l=1}^{3}\sum_{j=1}^{2} \delta_{i-1}\delta_{k-1}(-1)^{i+j+k+l} c_{k,l} c_{\kappa(i,k),\kappa(j,l)} c_{\lambda(i,k),\lambda(j,l)}}, \quad \{l,k\} \in \{\{1,2,3\} \times \{1,2,3\}\}$$

(5.19)

$$\mathbf{B}^{-1} \sim \beta_{nm} = \frac{(-1)^{m+n} \sum_{l=1}^{3}\sum_{j=1}^{2} \delta_{i-1}\delta_{k-1}(-1)^{i+j+k+l} b_{\kappa(k,m),\kappa(l,n)} \cdot b_{\lambda(i,k,m),\lambda(j,l,n)} \cdot b_{\mu(i,k,m),\mu(j,l,n)}}{\sum_{n=1}^{4}\sum_{l=1}^{3}\sum_{j=1}^{2} \delta_{i-1}\delta_{k-1}\delta_{m-1}(-1)^{i+j+k+l+m+n} b_{mn} \cdot b_{\kappa(k,m),\kappa(l,n)} \cdot b_{\lambda(i,k,m),\lambda(j,l,n)} \cdot b_{\mu(i,k,m),\mu(j,l,n)}},$$

$$\{n,m\} \in \{\{1,2,3,4\} \times \{1,2,3,4\}\}.$$

(5.20)

A more explicit expression for **A** is now found, in which all the elements in (5.14) are replaced by *a*'s, but with all indices primed where required,

$$\mathbf{A}^{-1} \sim \alpha_{qp} = \frac{(-1)^{p+q} \sum_{n=1}^{4} \sum_{l=1}^{3} \sum_{j=1}^{2} \delta_{i-1} \delta_{k-1} \delta_{m-1} (-1)^{i+j+k+l+m+n} a_{m'n'}(p,q) \cdot a_{k''l''}(m,n) \cdot \left(a_{i'''j'''} a_{(3-i)''',(3-j)'''}\right)(k,l)}{\sum_{q=1}^{5} \sum_{n=1}^{4} \sum_{l=1}^{3} \sum_{j=1}^{2} \delta_{i-1} \delta_{k-1} \delta_{m-1} \delta_{p-1} (-1)^{i+j+k+l+m+n+p+q} a_{pq} \cdot a_{m'n'}(p,q) \cdot a_{k''l''}(m,n) \cdot \left(a_{i'''j'''} a_{(3-i)''',(3-j)'''}\right)(k,l)},$$

$$\{q, p\} \in \{\{1,2,3,4,5\} \times \{1,2,3,4,5\}\}. \tag{5.21}$$

As previously explained, an index is primed depending on the number of telescoping steps required to attain its corresponding matrix element, which is being replaced by $a$. For instance, the indices $i$ and $j$ are primed three times since 3 telescoping steps are required to derive the 2 × 2 matrix $\mathbf{D_C}$ (5.17) from the original 5 × 5 matrix $\mathbf{A}$ (1.5). The indices can be simplified as follows, beginning with

$$m' = m + 1 - \mathrm{H}_{p-m-1} = \kappa(m, p), \tag{5.22}$$

$$n' = m'\delta_{m-n}\delta_{p-q} = n + 1 - \mathrm{H}_{q-n-1} = \kappa(n, q), \tag{5.23}$$

from the LHS of (5.10). For the double-primed row index $k''$, it simplifies as

$$k'' = (k + 1 - \mathrm{H}_{m-k-1})' = \begin{cases} (2 - \mathrm{H}_{m-2})', k = 1 \\ (3 - \mathrm{H}_{m-3})', k = 2 \end{cases} = \begin{cases} 2'\mathrm{H}_{1-m} + 1'\mathrm{H}_{m-2}, k = 1 \\ 3'\mathrm{H}_{2-m} + 2'\mathrm{H}_{m-3}, k = 2 \end{cases} = \begin{cases} (2 - \mathrm{H}_{p-2})\mathrm{H}_{m-2} + (3 - \mathrm{H}_{p-3})\mathrm{H}_{1-m}, k = 1 \\ (3 - \mathrm{H}_{p-3})\mathrm{H}_{m-3} + (4 - \mathrm{H}_{p-4})\mathrm{H}_{2-m}, k = 2 \end{cases}$$
$$\tag{5.24}$$

which can be summarized as

$$k'' = (k + 1 - \mathrm{H}_{p-k-1})\mathrm{H}_{m-k-1} + (k + 2 - \mathrm{H}_{p-k-2})\mathrm{H}_{k-m}; \ k \in \{1,2,3\}, \ m \in \{1,\cdots,4\}, \ p \in \{1,\cdots,5\}. \tag{5.25}$$

It can be simplified to the compact expression

$$k'' = \sum_{u=1}^{2}(k + u - \mathrm{H}_{p-k-u})\mathrm{H}_{(-1)^u(k-u-m+2)} = \sum_{u=1}^{2}\kappa(k+u-1, p)\mathrm{H}_{(-1)^u(k-u-m+2)} = \lambda(k, m, p); \tag{5.26}$$

$$k \in \{1,2,3\}, \ m \in \{1,2,3,4\}, \ p \in \{1,2,3,4,5\}.$$

Applying the following operations on the double-primed row index $k''$ yields the double-primed column index $l''$:

$$l'' = k''\delta_{k-l}\delta_{m-n}\delta_{p-q} = \sum_{u=1}^{2}(l + u - \mathrm{H}_{q-l-u})\mathrm{H}_{(-1)^u(l-u-n+2)} = \sum_{u=1}^{2}\kappa(l+u-1, q)\mathrm{H}_{(-1)^u(l-u-n+2)} = \lambda(l, n, q);$$

$$l \in \{1,2,3\}, \ n \in \{1,2,3,4\}, \ q \in \{1,2,3,4,5\}. \tag{5.27}$$

As for the triple-primed row index

$$i''' = \left(i+1-\mathrm{H}_{k-i-1}\right)'' = \begin{cases} \left(2-\mathrm{H}_{k-2}\right)'', i=1 \\ \left(3-\mathrm{H}_{k-3}\right)'', i=2 \end{cases} = \begin{cases} 2''\mathrm{H}_{1-k}+1''\mathrm{H}_{k-2}, i=1 \\ 3''\mathrm{H}_{2-k}+2''\mathrm{H}_{k-3}, i=2 \end{cases} = \begin{cases} \left(2-\mathrm{H}_{m-2}\right)'\mathrm{H}_{k-2}+\left(3-\mathrm{H}_{m-3}\right)'\mathrm{H}_{1-k}, i=1 \\ \left(3-\mathrm{H}_{m-3}\right)'\mathrm{H}_{k-3}+\left(4-\mathrm{H}_{m-4}\right)'\mathrm{H}_{2-k}, i=2 \end{cases}$$

$$= \begin{cases} \left(2'\mathrm{H}_{1-m}+1'\mathrm{H}_{m-2}\right)\mathrm{H}_{k-2}+\left(3'\mathrm{H}_{2-m}+2'\mathrm{H}_{m-3}\right)\mathrm{H}_{1-k}, i=1 \\ \left(3'\mathrm{H}_{2-m}+2'\mathrm{H}_{m-3}\right)\mathrm{H}_{k-3}+\left(4'\mathrm{H}_{3-m}+3'\mathrm{H}_{m-4}\right)\mathrm{H}_{2-k}, i=2 \end{cases}$$

$$= \begin{cases} \left(\left(2-\mathrm{H}_{p-2}\right)\mathrm{H}_{m-2}+\left(3-\mathrm{H}_{p-3}\right)\mathrm{H}_{1-m}\right)\mathrm{H}_{k-2}+\left(\left(3-\mathrm{H}_{p-3}\right)\mathrm{H}_{m-3}+\left(4-\mathrm{H}_{p-4}\right)\mathrm{H}_{2-m}\right)\mathrm{H}_{1-k}, i=1 \\ \left(\left(3-\mathrm{H}_{p-3}\right)\mathrm{H}_{m-3}+\left(4-\mathrm{H}_{p-4}\right)\mathrm{H}_{2-m}\right)\mathrm{H}_{k-3}+\left(\left(4-\mathrm{H}_{p-4}\right)\mathrm{H}_{m-4}+\left(5-\mathrm{H}_{p-5}\right)\mathrm{H}_{3-m}\right)\mathrm{H}_{2-k}, i=2 \end{cases}$$

(5.28)

which can be generalized for any value of *i* being 1 or 2 as

$$i''' = \left(\left(i+1-\mathrm{H}_{p-i-1}\right)\mathrm{H}_{m-i-1}+\left(i+2-\mathrm{H}_{p-i-2}\right)\mathrm{H}_{i-m}\right)\mathrm{H}_{k-i-1}+\left(\left(i+2-\mathrm{H}_{p-i-2}\right)\mathrm{H}_{m-i-2}+\left(i+3-\mathrm{H}_{p-i-3}\right)\mathrm{H}_{i+1-m}\right)\mathrm{H}_{i-k}$$

(5.29)

and which can also be contracted to:

$$i''' = \sum_{u=1}^{2}\sum_{v=0}^{1}\left(i+u+v-\mathrm{H}_{p-i-u-v}\right)\mathrm{H}_{(-1)^u(i-u-m+v+2)}\mathrm{H}_{(-1)^v(k-i+v-1)} = \mu(i,k,m,p). \qquad (5.30)$$

The corresponding column index is obtainable from the above results as follows

$$j''' = i'''\delta_{i-j}\delta_{k-l}\delta_{m-n}\delta_{p-q} \qquad (5.31)$$

which involves the replacement in $i'''$ (5.30), of *i* with *j*, *k* with *l*, *m* with *n*, and *p* with *q*, with the result,

$$j''' = \sum_{u=1}^{2}\sum_{v=0}^{1}\left(j+u+v-\mathrm{H}_{q-j-u-v}\right)\mathrm{H}_{(-1)^u(j-u-n+v+2)}\mathrm{H}_{(-1)^v(l-j+v-1)} = \mu(j,l,n,q). \qquad (5.32)$$

As for

$$(3-i)''' = i''' \otimes \delta_{3-i} = \sum_{u=1}^{2}\sum_{v=0}^{1}\left(-i+u+v+3-\mathrm{H}_{p+i-u-v-3}\right)\mathrm{H}_{(-1)^u(-i-u+v-m+5)}\mathrm{H}_{(-1)^v(i+k+v-4)} = \nu(i,k,m,p)$$

(5.33)

and is obtained by discrete convolution on (5.30), as explained in the previous section. The function $\nu(\bullet)$, which is Greek, should not be confused with the summation index *v* on the LHS. Lastly, the corresponding column index is obtained from (5.32) after a few simple operations,

$$(3-j)''' = (3-i)'''\delta_{i-j}\delta_{k-l}\delta_{m-n}\delta_{p-q} = \nu(j,l,n,q). \qquad (5.34)$$

Consequently,

$$\mathbf{A}^{-1} \sim \alpha_{qp} = \frac{(-1)^{p+q} \sum_{n=1}^{4} \sum_{l=1}^{3} \sum_{j=1}^{2} \delta_{i-1} \delta_{k-1} \delta_{m-1} (-1)^{i+j+k+l+m+n} a_{\kappa(m,p),\kappa(n,q)} \cdot a_{\lambda(k,m,p),\lambda(l,n,q)} \cdot a_{\mu(i,k,m,p),\mu(j,l,n,q)} \cdot a_{\nu(i,k,m,p),\nu(j,l,n,q)}}{\sum_{q=1}^{5} \sum_{n=1}^{4} \sum_{l=1}^{3} \sum_{j=1}^{2} \delta_{i-1} \delta_{k-1} \delta_{m-1} \delta_{p-1} (-1)^{i+j+k+l+m+n+p+q} a_{pq} \cdot a_{\kappa(m,p),\kappa(n,q)} \cdot a_{\lambda(k,m,p),\lambda(l,n,q)} \cdot a_{\mu(i,k,m,p),\mu(j,l,n,q)} \cdot a_{\nu(i,k,m,p),\nu(j,l,n,q)}}.$$

(5.35)

After enforcing the various Kronecker delta-functions, the expression simplifies to

$$\mathbf{A}^{-1} \sim \alpha_{qp} = \frac{(-1)^{p+q} \sum_{n=1}^{4} \sum_{l=1}^{3} \sum_{j=1}^{2} (-1)^{1+j+l+n} a_{\kappa(1,p),\kappa(n,q)} \cdot a_{\lambda(1,1,p),\lambda(l,n,q)} \cdot a_{\mu(1,1,1,p),\mu(j,l,n,q)} \cdot a_{\nu(1,1,1,p),\nu(j,l,n,q)}}{\sum_{q=1}^{5} \sum_{n=1}^{4} \sum_{l=1}^{3} \sum_{j=1}^{2} (-1)^{j+l+n+q} a_{1q} \cdot a_{\kappa(1,1),\kappa(n,q)} \cdot a_{\lambda(1,1,1),\lambda(l,n,q)} \cdot a_{\mu(1,1,1,1),\mu(j,l,n,q)} \cdot a_{\nu(1,1,1,1),\nu(j,l,n,q)}},$$

(5.36)

where

$$\kappa(1,p) = 2 - \mathrm{H}_{p-2}, \tag{5.37}$$

$$\kappa(1,1) = 2 - \mathrm{H}_{-1} = 2, \tag{5.38}$$

$$\lambda(1,1,p) = \left(2 - \mathrm{H}_{p-2}\right)\mathrm{H}_{-1} + \left(3 - \mathrm{H}_{p-3}\right)\mathrm{H}_0 = 3 - \mathrm{H}_{p-3} \tag{5.39}$$

$$\lambda(1,1,1) = 3, \tag{5.40}$$

$$\mu(1,1,1,p) = \left(\left(2 - \mathrm{H}_{p-2}\right)\mathrm{H}_{-1} + \left(3 - \mathrm{H}_{p-3}\right)\mathrm{H}_0\right)\mathrm{H}_{-1} + \left(\left(3 - \mathrm{H}_{p-3}\right)\mathrm{H}_{-2} + \left(4 - \mathrm{H}_{p-4}\right)\mathrm{H}_1\right)\mathrm{H}_0 = 4 - \mathrm{H}_{p-4} \tag{5.41}$$

$$\mu(1,1,1,1) = 4, \tag{5.42}$$

$$\nu(1,1,1,p) = \left(\left(4 - i - \mathrm{H}_{p-4+i}\right)\mathrm{H}_{m-4+i} + \left(5 - i - \mathrm{H}_{p-5+i}\right)\mathrm{H}_{3-i-m}\right)\mathrm{H}_{-2} + \left(\left(5 - i - \mathrm{H}_{p-5+i}\right)\mathrm{H}_{-3} + \left(5 - \mathrm{H}_{p-5}\right)\mathrm{H}_2\right)\mathrm{H}_1 = 5 - \mathrm{H}_{p-5} \tag{5.43}$$

$$\nu(1,1,1,1) = 5, \tag{5.44}$$

and finally,

$$\mathbf{A}^{-1} \sim \alpha_{qp} = \frac{(-1)^{p+q} \sum_{n=1}^{4} \sum_{l=1}^{3} \sum_{j=1}^{2} (-1)^{1+j+l+n} a_{2-\mathrm{H}_{p-2},\kappa(n,q)} \cdot a_{3-\mathrm{H}_{p-3},\lambda(l,n,q)} \cdot a_{4-\mathrm{H}_{p-4},\mu(j,l,n,q)} \cdot a_{5-\mathrm{H}_{p-5},\nu(j,l,n,q)}}{\sum_{q=1}^{5} \sum_{n=1}^{4} \sum_{l=1}^{3} \sum_{j=1}^{2} (-1)^{j+l+n+q} a_{1q} \cdot a_{2,\kappa(n,q)} \cdot a_{3,\lambda(l,n,q)} \cdot a_{4,\mu(j,l,n,q)} \cdot a_{5,\nu(j,l,n,q)}},$$

$$\{q, p\} \in \{\{1,2,3,4,5\} \times \{1,2,3,4,5\}\}.$$

(5.45)

where $\kappa(n,q), \lambda(l,n,q), \mu(j,l,n,q),$ and $\nu(j,l,n,q)$ are respectively given by (5.23), (5.27), (5.32) and (5.34). The determinant is found in the denominator of (5.45),

$$|\mathbf{A}| = \sum_{q=1}^{5}\sum_{n=1}^{4}\sum_{l=1}^{3}\sum_{j=1}^{2}(-1)^{j+l+n+q} a_{1,q} \cdot a_{2,\kappa(n,q)} \cdot a_{3,\lambda(l,n,q)} \cdot a_{4,\mu(j,l,n,q)} \cdot a_{5,\nu(j,l,n,q)}. \qquad (5.46)$$

No mathematical proof is proffered for the conjecture (5.45) at this juncture[5]. However, in order to validate (5.45), one million 5 x 5 matrices were generated in Matlab, using the randn(5)-function, which produced 5 x 5 matrices comprised of pseudo-random elements selected from the standard normal distribution. The inverse $(\alpha_{qp})$ of an $r$-th matrix was found using (5.45), as well as by Matlab's inv-function, represented by $\xi_{pq}$. To expedite the execution of the program, the matrix elements were constrained to be real. The computer that hosted the Matlab software, was equipped with 8 GB of RAM, an Intel® Core™2 Duo Processor T9900 at a base frequency of 3.06 GHz and with a 6-MB L2 Cache, and the Intel 64® architecture.

The mean-squared error (MSE) $\varepsilon_r$ for each $r$-th matrix was found using

$$\varepsilon_r = \frac{1}{5^2}\sum_{p=1}^{5}\sum_{q=1}^{5}\left|\alpha_{qp} - \xi_{pq}\right|_r^2, \; r \in \{1,2,3,\cdots,999998,999999,10^6\}. \qquad (5.47)$$

A histogram, of 10 times the base-10-logarithm of the MSE, was then generated using the Matlab hist- or histogram-function. In order to avoid zero MSE values that would yield infinite outcomes under the logarithm, the MSE vector $\boldsymbol{\varepsilon_r}$ was initialized to constants of the order of $10^{-100}$. The result is found in fig. 1. The inset of fig. 1 demonstrates the Gaussian behavior of the Matlab randn-function for a sample size of $10^4$, whereas the sample size of each matrix element was $10^6$.

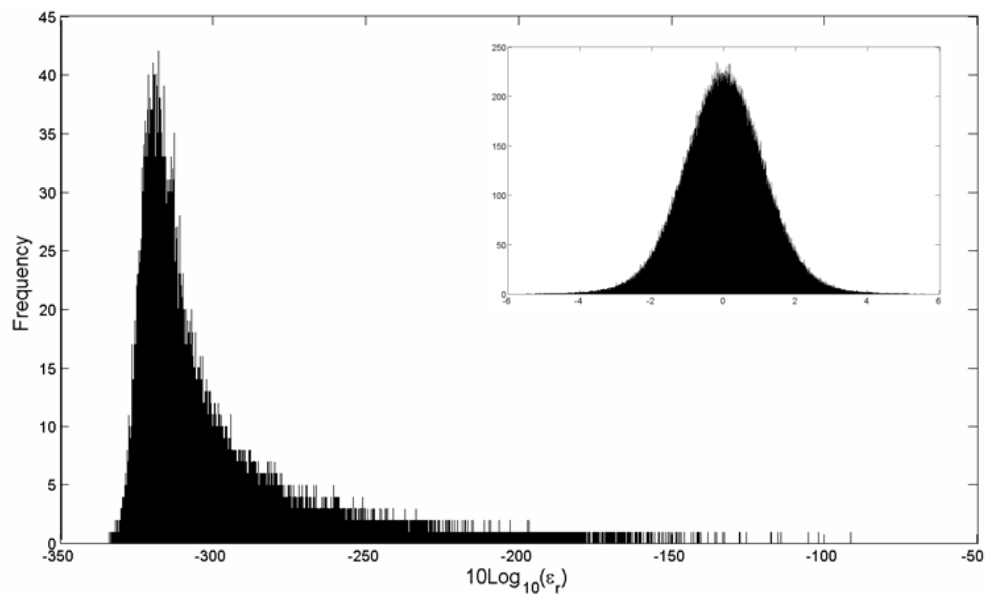

Fig. 1. A histogram of the distribution of $10\log_{10}\varepsilon_r$ (5.47) for one million, 5 x 5 random matrices. The inset shows a typical histogram generated by the Matlab randn-function.

---
[5] A Matlab implementation of (5.45) is provided in the **Appendix**

The minimum MSE, which was non-zero, was found to be $\approx 4 \times 10^{-34}$ (-334 in fig. 1), while the maximum MSE was no more than $\approx 7 \times 10^{-10}$ (-91.6 in fig. 1). The mode or most-frequent MSE, is seen to be of the order of $10^{-32}$ (-320 in fig. 1). The non-zero MSE is attributed to the different computational methods used in the 2 approaches, but is quite small regardless, validating (5.45) for the inverse of a randomly generated 5 x 5 matrix. The majority of the MSE is observed to be $< 10^{-20}$ (-200 in fig. 1), and the distribution exhibits a non-negative skew-ness. The MSE distribution would not be expected to be Gaussian, since the determinant of a square matrix comprised of i.i.d. Gaussian random variables, would not be Gaussian in general. Therefore, the MSE (5.47) would be even less likely to be Gaussian, as it involves non-linear operations such as the logarithm and the modulus-squared.

It is known that if at least one row (or column) of a square matrix is comprised entirely of zeros, the determinant resolves to zero, and the associated matrix would not be invertible. However, the determinant also resolves to zero when this is not the case for some inversion techniques on sparse, but invertible square matrices, such as the following matrices,

$$\begin{bmatrix} 0 & 0 & 0 & x_1 & 0 \\ 0 & x_2 & 0 & 0 & 0 \\ x_3 & 0 & 0 & 0 & 0 \\ 0 & 0 & x_4 & 0 & 0 \\ 0 & 0 & 0 & 0 & x_5 \end{bmatrix}, \begin{bmatrix} x_1 & 0 & 0 & 0 & 0 \\ 0 & x_2 & 0 & 0 & 0 \\ 0 & 0 & 0 & 0 & x_3 \\ 0 & 0 & x_4 & 0 & 0 \\ 0 & 0 & 0 & x_5 & 0 \end{bmatrix}, \begin{bmatrix} 0 & x_1 & 0 & 0 & 0 \\ x_2 & 0 & 0 & 0 & 0 \\ 0 & 0 & x_3 & 0 & 0 \\ 0 & 0 & 0 & 0 & x_4 \\ 0 & 0 & 0 & x_5 & 0 \end{bmatrix}.$$

(5.48)

Preliminary work indicates that the technique reported herein is immune to such cases. As exemplified by (5.48), an invertible $N \times N$ matrix need not support more than $N$ non-zero elements in order to remain invertible, or it can equivalently support as many as ($N^2$-$N$) appropriately located zeros and still be invertible. There are no restrictions on the values of $x$, and they may all be identical or dissimilar, as long as the matrix remains invertible. Matrices such as those shown in (5.48) carry many inner zero-elements, which can be problematic for certain determinant techniques such as Dodgson Condensation [8, 18, 19], and the recently reported recursive technique due to Razaifar and Rezaee [26]. The problem is usually alleviated by some additional (pre-)processing steps on the matrix in question [8, 27]. In order to establish the viability of the technique (5.45) proposed, but for invertible sparse 5 x 5 matrices, (5.48) were selected at random, and used as test cases. There was no disagreement between (5.46) and Matlab's det-function, for any of the 3 matrices. Furthermore, (5.45) produced identical results to Matlab's inv-function, again for all 3 matrices. Although this evaluation by no means exhausts the set of all possible invertible, maximally sparse 5 x 5 matrices, the technique is at least valid for the test cases (5.48). A more rigorous evaluation is perhaps required, but the results reported here represent a satisfactory preliminary effort.

## 6. Generalization to a N x N matrix A

Based on the previous results, it can be surmised that, given a general $N \times N$ square matrix $\mathbf{A}$,

$$\mathbf{A} = \begin{bmatrix} a_{11} & a_{12} & \cdots & a_{1N} \\ a_{21} & a_{22} & \cdots & a_{2N} \\ \vdots & \vdots & \ddots & \vdots \\ a_{N1} & a_{N2} & \cdots & a_{NN} \end{bmatrix} \tag{6.1}$$

which is indexed to $\{r_0, s_0\}$, and with an inverse of the form of

$$\mathbf{A}^{-1} \sim \alpha_{s_0, r_0} \sim \text{(matrix of cofactor minors as shown)} \tag{6.2}$$

its general $(N-1) \times (N-1)$ minor-matrix $\mathbf{M}$ would be expected to have the form of

$$\mathbf{M}_{r_1, s_1}(r_0, s_0) = \mathbf{A}\left(1+\delta_{r_0-1}, 2+\delta_{r_0-1}+\delta_{r_0-2}, \cdots, N-1+\sum_{x=1}^{N-1}\delta_{r_0-x}; 1+\delta_{s_0-1}, 2+\delta_{s_0-1}+\delta_{s_0-2}, \cdots, N-1+\sum_{y=1}^{N-1}\delta_{s_0-y}\right),$$

$$\{r_1, s_1\} \in \{\{1, \cdots, N-1\} \times \{1, \cdots, N-1\}\}, \{r_0, s_0\} \in \{\{1, \cdots, N\} \times \{1, \cdots, N\}\}, N \geq 2 \ \& \ N \in \mathbb{Z}. \tag{6.3}$$

The matrix $\mathbf{M}$ is a function of the indices $\{r_0, s_0\}$ of the original matrix $\mathbf{A}$, but $\mathbf{M}$ itself is actually indexed to $\{r_1, s_1\}$. Eq. (6.3) may be re-stated concisely as

$$\mathbf{M}(r_0, s_0) = \mathbf{A}\left(\bigcup_{\sigma=1}^{N-1}\left(\sigma + \sum_{x=1}^{\sigma}\delta_{r_0-x}\right); \bigcup_{\tau=1}^{N-1}\left(\tau + \sum_{y=1}^{\tau}\delta_{s_0-y}\right)\right), \tag{6.4}$$

where in general,

$$\bigcup_{\xi=1}^{N-1}\left(\xi+\sum_{w=1}^{\xi}\delta_{t-w}\right)=\left\{1+\delta_{t-1},2+\delta_{t-1}+\delta_{t-2},\cdots,N-1+\sum_{w=1}^{N-1}\delta_{t-w}\right\},\ \{t,w,\xi\}\in\{\{r_0,x,\sigma\},\{s_0,y,\tau\}\}.$$
(6.5)

However, the braces would be discarded due to their redundancy, if (6.5) were utilized as an argument in the minor-matrix expression (6.3). For instance, applying this conclusion to a 2 x 2 matrix **D** for which $N = 2$, its general minor-matrix (6.4) is reduced to

$$\mathbf{D}\left(\bigcup_{\sigma=1}^{1}\left(\sigma+\sum_{x=1}^{\sigma}\delta_{r_0-x}\right);\bigcup_{\tau=1}^{1}\left(\tau+\sum_{y=1}^{\tau}\delta_{s_0-y}\right)\right)=\mathbf{D}\left(1+\delta_{r_0-1},1+\delta_{s_0-1}\right)=d_{1+\delta_{r_0-1},1+\delta_{s_0-1}},\ \{r_0,s_0\}\in\{\{1,2\}\times\{1,2\}\}$$
(6.6)

and the minor-matrix trivially reduces to the elements of the matrix **D** itself, for a given $(r_0, s_0) = (i, j)$. The expression (6.4) yields the correct result for $N = 2$, and also reproduces (3.6) and (4.5), respectively for $N = 3$, and 4.

Using the equivalence,

$$\xi+\sum_{w=1}^{\xi}\delta_{t-w}=\sum_{w=1}^{\xi}(w+\delta_{t-w})=\xi+1-\mathrm{H}_{t-\xi-1}$$
(6.7)

an alternative version of the general minor-matrix (6.4) is found to be

$$\mathbf{M}(r_0,s_0)=\mathbf{A}\left(\bigcup_{\sigma=2}^{N}(\sigma-\mathrm{H}_{r_0-\sigma});\bigcup_{\tau=2}^{N}(\tau-\mathrm{H}_{s_0-\tau})\right),\ \{r_0,s_0\}\in\{\{1,\cdots,N\}\times\{1,\cdots,N\}\}$$
(6.8)

or in long form, as

$$\mathbf{M}(r_0,s_0)=\mathbf{A}\left(2-\mathrm{H}_{r_0-2},3-\mathrm{H}_{r_0-3},\cdots,N-\mathrm{H}_{r_0-N};2-\mathrm{H}_{s_0-2},3-\mathrm{H}_{s_0-3},\cdots,N-\mathrm{H}_{s_0-N}\right),$$
$$\{r_0,s_0\}\in\{\{1,\cdots,N\}\times\{1,\cdots,N\}\},\ N\geq 2.$$
(6.9)

The matrix version is given by

$$\mathbf{M}(r_0,s_0)\sim a_{r_1+1-\mathrm{H}_{r_0-r_1-1},s_1+1-\mathrm{H}_{s_0-s_1-1}}\sim\begin{bmatrix}a_{2-\mathrm{H}_{r_0-2},2-\mathrm{H}_{s_0-2}} & a_{2-\mathrm{H}_{r_0-2},3-\mathrm{H}_{s_0-3}} & \cdots & a_{2-\mathrm{H}_{r_0-2},N-\mathrm{H}_{s_0-N}}\\ a_{3-\mathrm{H}_{r_0-3},2-\mathrm{H}_{s_0-2}} & a_{3-\mathrm{H}_{r_0-3},3-\mathrm{H}_{s_0-3}} & \cdots & a_{3-\mathrm{H}_{r_0-3},N-\mathrm{H}_{s_0-N}}\\ \vdots & \vdots & \ddots & \vdots\\ a_{N-\mathrm{H}_{r_0-N},2-\mathrm{H}_{s_0-2}} & a_{N-\mathrm{H}_{r_0-N},3-\mathrm{H}_{s_0-3}} & \cdots & a_{N-\mathrm{H}_{r_0-N},N-\mathrm{H}_{s_0-N}}\end{bmatrix},$$
$$\{r_1,s_1\}\in\{\{1,\cdots,N-1\}\times\{1,\cdots,N-1\}\},\ \{r_0,s_0\}\in\{\{1,\cdots,N\}\times\{1,\cdots,N\}\},\ N\geq 2.$$
(6.10)

At this juncture, it seems that the matrix indices are unnecessarily over-subscripted.

However, it will be seen later that using these indices simplify the ultimate general expressions for the inverse of **A** (6.1). In order to use the telescoping method, a new indexing system is adopted, since the alphabet set is easily exhausted. Beginning with the primary minor-matrix $\mathbf{M}^{(N-1)}_{r_1,s_1}(r_0, s_0)$, whose determinants are the entries of the numerator of the inverse of **A** (6.2),

$$\mathbf{M}^{(N-1)}_{r_1,s_1} = \mathbf{A}\left(\bigcup_{\sigma_1=2}^{N}(\sigma_1 - \mathrm{H}_{r_0-\sigma_1}); \bigcup_{\tau_1=2}^{N}(\tau_1 - \mathrm{H}_{s_0-\tau_1})\right), \{r_1,s_1\} \in \{\{1,\cdots,N-1\}\times\{1,\cdots,N-1\}\},$$

$$\mathbf{M}^{(N-2)}_{r_2,s_2} = \mathbf{M}^{(N-1)}_{r_1,s_1}\left(\bigcup_{\sigma_2=2}^{N-1}(\sigma_2 - \mathrm{H}_{r_1-\sigma_2}); \bigcup_{\tau_2=2}^{N-1}(\tau_2 - \mathrm{H}_{s_1-\tau_2})\right), \{r_2,s_2\} \in \{\{1,\cdots,N-2\}\times\{1,\cdots,N-2\}\},$$

$$\vdots$$

$$\mathbf{M}^{(3)}_{r_{N-3},s_{N-3}} = \mathbf{M}^{(4)}_{r_{N-4},s_{N-4}}\left(\bigcup_{\sigma_{N-3}=2}^{4}(\sigma_{N-3} - \mathrm{H}_{r_{N-4}-\sigma_{N-3}}); \bigcup_{\tau_{N-3}=2}^{4}(\tau_{N-3} - \mathrm{H}_{s_{N-4}-\tau_{N-3}})\right), \{r_{N-3},s_{N-3}\} \in \{\{1,2,3\}\times\{1,2,3\}\},$$

$$\mathbf{M}^{(2)}_{r_{N-2},s_{N-2}} = \mathbf{M}^{(3)}_{r_{N-3},s_{N-3}}\left(\bigcup_{\sigma_{N-2}=2}^{3}(\sigma_{N-2} - \mathrm{H}_{r_{N-3}-\sigma_{N-2}}); \bigcup_{\tau_{N-2}=2}^{3}(\tau_{N-2} - \mathrm{H}_{s_{N-3}-\tau_{N-2}})\right), \{r_{N-2},s_{N-2}\} \in \{\{1,2\}\times\{1,2\}\}.$$

(6.11)

which are valid for $N \geq 3$. By convention, the dimension of the matrix **A** is stated as $N \times N = N^2$, and in the above progression, the *superscript* of a minor-matrix **M** indicates the square-root of its dimension. The *index subscript* indicates the number of telescoping steps used to attain its corresponding minor-matrix **M**. For instance, if **A** (6.1) were a 3 × 3 matrix, the couple $(r_{N-2}, s_{N-2}) = (r_1, s_1)$ for which the indices are both subscripted with a 1 implies that its corresponding 2 × 2 minor-matrix $\mathbf{M}^{(2)}_{r_{N-2},s_{N-2}} = \mathbf{M}^{(2)}_{r_1,s_1}$ is obtained after 1 telescoping step from **A**. In general, the index subscript is obtained by negating the minor-matrix superscript, then adding the square-root of the dimension of **A**.

Once all the minor-matrices are found, they are processed in the formula for the inverse of matrix **A**[6]

$$\mathbf{A}^{-1} \sim \alpha_{s_0,r_0} = \frac{(-1)^{r_0+s_0}\sum_{s_1=1}^{N-1}\cdots\sum_{s_{N-4}=1}^{4}\sum_{s_{N-3}=1}^{3}\sum_{s_{N-2}=1}^{2}\delta_{r_1-1}\cdots\delta_{r_{N-4}-1}\delta_{r_{N-3}-1}\delta_{r_{N-2}-1}\cdot(-1)^{\sum_{k=2}^{N-1}r_{N-k}+s_{N-k}}m^{(N-1)}_{r_1,s_1}(r_0,s_0)\cdots m^{(4)}_{r_{N-4},s_{N-4}}(r_{N-5},s_{N-5})m^{(3)}_{r_{N-3},s_{N-3}}(r_{N-4},s_{N-4})\left(m^{(2)}_{r_{N-2},s_{N-2}}m^{(2)}_{3-r_{N-2},3-s_{N-2}}\right)(r_{N-3},s_{N-3})}{\sum_{s_0=1}^{N}\sum_{s_1=1}^{N-1}\cdots\sum_{s_{N-4}=1}^{4}\sum_{s_{N-3}=1}^{3}\sum_{s_{N-2}=1}^{2}\delta_{r_1-1}\delta_{r_1-1}\cdots\delta_{r_{N-4}-1}\delta_{r_{N-3}-1}\delta_{r_{N-2}-1}\cdot(-1)^{\sum_{k=2}^{N}r_{N-k}+s_{N-k}}a_{r_0,s_0}m^{(N-1)}_{r_1,s_1}(r_0,s_0)\cdots m^{(4)}_{r_{N-4},s_{N-4}}(r_{N-5},s_{N-5})m^{(3)}_{r_{N-3},s_{N-3}}(r_{N-4},s_{N-4})\left(m^{(2)}_{r_{N-2},s_{N-2}}m^{(2)}_{3-r_{N-2},3-s_{N-2}}\right)(r_{N-3},s_{N-3})},$$

$$\{s_0,r_0\} \in \{\{1,\cdots,N\}\times\{1,\cdots,N\}\}, \quad N \geq 3.$$

(6.12)

with the implication that there are (*N*-2) summations in the numerator, and (*N*-1) summations in the denominator, as evidenced by the respective upper bounds of the summations. Moreover, a Kronecker delta-function is associated with each summation, and an extant (-1) is raised to the sum of all the row and column indices of all the minors involved in telescoping toward the 2 × 2 minor-matrix $\mathbf{M}^{(2)}_{r_{N-2},s_{N-2}}$. Eq. (6.12) also shows the functional dependence of each matrix element $m^{(\rho)}_{r_{N-\rho},s_{N-\rho}} \in \mathbf{M}^{(\rho)}_{r_{N-\rho},s_{N-\rho}}$. Eq. (6.12) can be re-expressed more concisely as

$$\mathbf{A}^{-1} \sim \alpha_{s_0,r_0} = \frac{(-1)^{r_0+s_0} \prod_{\rho=2}^{N-1} \sum_{s_{N-\rho}=1}^{\rho} \delta_{r_{N-\rho}-1} \cdot (-1)^{r_{N-\rho}+s_{N-\rho}} m^{(\rho)}_{r_{N-\rho},s_{N-\rho}} \cdot \left(\mathrm{H}_{\rho-3} + \delta_{\rho-2} m^{(2)}_{3-r_{N-2},3-s_{N-2}}\right)}{\prod_{\rho=2}^{N} \sum_{s_{N-\rho}=1}^{\rho} \delta_{r_{N-\rho}-1} \cdot (-1)^{r_{N-\rho}+s_{N-\rho}} \left(\delta_{\rho-N} a_{r_0,s_0} + \mathrm{H}_{N-1-\rho} m^{(\rho)}_{r_{N-\rho},s_{N-\rho}}\right)\left(\mathrm{H}_{\rho-3} + \delta_{\rho-2} m^{(2)}_{3-r_{N-2},3-s_{N-2}}\right)},$$

$$\{s_0, r_0\} \in \{\{1,\cdots,N\} \times \{1,\cdots,N\}\}, \ N \geq 3.$$

(6.13)

where the direct functional dependence of the various matrix elements is as follows:

$$m^{(\rho)}_{r_{N-\rho},s_{N-\rho}} \sim m^{(\rho)}_{r_{N-\rho},s_{N-\rho}}\left(r_{N-\rho-1}, s_{N-\rho-1}\right), \tag{6.14}$$

$$m^{(2)}_{3-r_{N-2},3-s_{N-2}} \sim m^{(2)}_{3-r_{N-2},3-s_{N-2}}\left(r_{N-3}, s_{N-3}\right). \tag{6.15}$$

---

[6] May have to magnify (6.12) for a better view

An application of (6.13) is now illustrated using the 4 x 4 matrix **B** examined in §**4**, with the aim of reproducing (4.15). In this case, $N = 4$, and

$$\mathbf{B}^{-1} \sim \beta_{s_0,r_0} = \frac{(-1)^{r_0+s_0} \prod_{\rho=2}^{3} \sum_{s_{4-\rho}=1}^{\rho} \delta_{r_{4-\rho}-1} \cdot (-1)^{r_{4-\rho}+s_{4-\rho}} m^{(\rho)}_{r_{4-\rho},s_{4-\rho}}\left(r_{3-\rho}, s_{3-\rho}\right) \cdot \left(\mathrm{H}_{\rho-3} + \delta_{\rho-2} m^{(2)}_{3-r_2,3-s_2}\left(r_1,s_1\right)\right)}{\prod_{\rho=2}^{4} \sum_{s_{4-\rho}=1}^{\rho} \delta_{r_{4-\rho}-1} \cdot (-1)^{r_{4-\rho}+s_{4-\rho}} \left(\delta_{\rho-4} b_{r_0,s_0} + \mathrm{H}_{3-\rho} m^{(\rho)}_{r_{4-\rho},s_{4-\rho}}\left(r_{3-\rho}, s_{3-\rho}\right)\right)\left(\mathrm{H}_{\rho-3} + \delta_{\rho-2} m^{(2)}_{3-r_2,3-s_2}\left(r_1,s_1\right)\right)},$$

$$\{s_0, r_0\} \in \{\{1,2,3,4\} \times \{1,2,3,4\}\}.$$

(6.16)

After expanding the products,

$$\mathbf{B}^{-1} \sim \beta_{s_0,r_0} = \frac{(-1)^{r_0+s_0} \sum_{s_1=1}^{3} \delta_{r_1-1} \cdot (-1)^{r_1+s_1} m^{(3)}_{r_1,s_1}(r_0,s_0) \sum_{s_2=1}^{2} \delta_{r_2-1} \cdot (-1)^{r_2+s_2} \left(m^{(2)}_{r_2,s_2} m^{(2)}_{3-r_2,3-s_2}\right)(r_1,s_1)}{\sum_{s_0=1}^{4} \delta_{r_0-1} \cdot (-1)^{r_0+s_0} b_{r_0,s_0} \sum_{s_1=1}^{3} \delta_{r_1-1} \cdot (-1)^{r_1+s_1} m^{(3)}_{r_1,s_1}(r_0,s_0) \sum_{s_2=1}^{2} \delta_{r_2-1} \cdot (-1)^{r_2+s_2} \left(m^{(2)}_{r_2,s_2} m^{(2)}_{3-r_2,3-s_2}\right)(r_1,s_1)}$$

(6.17)

and simplifying

$$\mathbf{B}^{-1} \sim \beta_{s_0,r_0} = \frac{(-1)^{r_0+s_0} \sum_{s_1=1}^{3} \sum_{s_2=1}^{2} \delta_{r_1-1} \delta_{r_2-1} \cdot (-1)^{r_1+s_1+r_2+s_2} m^{(3)}_{r_1,s_1}(r_0,s_0) \cdot \left(m^{(2)}_{r_2,s_2} m^{(2)}_{3-r_2,3-s_2}\right)(r_1,s_1)}{\sum_{s_0=1}^{4} \sum_{s_1=1}^{3} \sum_{s_2=1}^{2} \delta_{r_0-1} \delta_{r_1-1} \delta_{r_2-1} \cdot (-1)^{r_0+s_0+r_1+s_1+r_2+s_2} b_{r_0,s_0} \cdot m^{(3)}_{r_1,s_1}(r_0,s_0) \cdot \left(m^{(2)}_{r_2,s_2} m^{(2)}_{3-r_2,3-s_2}\right)(r_1,s_1)}$$

(6.18)

Eq. (4.15) is reproduced after enforcing the various Kronecker delta-functions in (6.18),

$$\mathbf{B}^{-1} \sim \beta_{s_0, r_0} = \frac{(-1)^{r_0+s_0} \sum_{s_1=1}^{3} \sum_{s_2=1}^{2} (-1)^{2+s_1+s_2} m_{1,s_1}^{(3)}(r_0, s_0) \cdot \left(m_{1,s_2}^{(2)} m_{2,3-s_2}^{(2)}\right)(1, s_1)}{\sum_{s_0=1}^{4} \sum_{s_1=1}^{3} \sum_{s_2=1}^{2} (-1)^{3+s_0+s_1+s_2} b_{1,s_0} \cdot m_{1,s_1}^{(3)}(1, s_0) \cdot \left(m_{1,s_2}^{(2)} m_{2,3-s_2}^{(2)}\right)(1, s_1)} \tag{6.19}$$

and making the element re-assignment of $\{m^{(3)}, m^{(2)}\}$ to $\{c, d\}$, and of the indices $\{r_0, s_0\}$ to $\{m, n\}$, $\{r_1, s_1\}$ to $\{k, l\}$, and $\{r_2, s_2\}$ to $\{i, j\}$,

$$\mathbf{B}^{-1} \sim \beta_{nm} = \frac{(-1)^{m+n} \sum_{l=1}^{3} \sum_{j=1}^{2} (-1)^{2+j+l} c_{1,l}(m, n) \cdot \left(d_{1,j} d_{2,3-j}\right)(1, l)}{\sum_{n=1}^{4} \sum_{l=1}^{3} \sum_{j=1}^{2} (-1)^{3+j+l+n} b_{1,n} \cdot c_{1,l}(1, n) \cdot \left(d_{1,j} d_{2,3-j}\right)(1, l)}, \tag{6.20}$$

$$\{n, m\} \in \{\{1, 2, 3, 4\} \times \{1, 2, 3, 4\}\}.$$

It is not difficult to demonstrate the reverse, i.e. obtaining (6.16) from (6.20). The matrix elements $\{m^{(3)}, m^{(2)}\}$ should not be confused with the row index $m$ used in (6.20).

In the second approach, which is an extension of the above telescoping approach, all the matrix elements $m$ are replaced with the elements of the original matrix $\mathbf{A}$ (6.1), with their corresponding indices primed a number of times identical with the number of telescoping steps required to find them,

$$\mathbf{A}^{-1} \sim \alpha_{s_0, r_0} = \frac{(-1)^{r_0+s_0} \sum_{s_1=1}^{N-1} \cdots \sum_{s_{N-4}=1}^{4} \sum_{s_{N-3}=1}^{3} \sum_{s_{N-2}=1}^{2} \delta_{r_1-1} \cdots \delta_{r_{N-4}-1} \delta_{r_{N-3}-1} \delta_{r_{N-2}-1} \cdot (-1)^{\sum_{\Lambda=2}^{N-1} r_{N-\Lambda}+s_{N-\Lambda}} a_{r_1^{(1)}, s_1^{(1)}} \cdots a_{r_{N-3}^{(N-3)}, s_{N-3}^{(N-3)}} \cdot \left(a_{r_{N-2}^{(N-2)}, s_{N-2}^{(N-2)}} a_{(3-r_{N-2})^{(N-2)}, (3-s_{N-2})^{(N-2)}}\right)}{\sum_{s_0=1}^{N} \sum_{s_1=1}^{N-1} \cdots \sum_{s_{N-4}=1}^{4} \sum_{s_{N-3}=1}^{3} \sum_{s_{N-2}=1}^{2} \delta_{r_0-1} \delta_{r_1-1} \cdots \delta_{r_{N-4}-1} \delta_{r_{N-3}-1} \delta_{r_{N-2}-1} \cdot (-1)^{\sum_{\Lambda=2}^{N} r_{N-\Lambda}+s_{N-\Lambda}} a_{r_0, s_0} a_{r_1^{(1)}, s_1^{(1)}} \cdots a_{r_{N-3}^{(N-3)}, s_{N-3}^{(N-3)}} \cdot \left(a_{r_{N-2}^{(N-2)}, s_{N-2}^{(N-2)}} a_{(3-r_{N-2})^{(N-2)}, (3-s_{N-2})^{(N-2)}}\right)},$$

$$\{s_0, r_0\} \in \{\{1, \cdots, N\} \times \{1, \cdots, N\}\}, \quad N \geq 3,$$
(6.21)

or as

$$\mathbf{A}^{-1} \sim \alpha_{s_0, r_0} = \frac{(-1)^{r_0+s_0} \prod_{\rho=2}^{N-1} \sum_{s_{N-\rho}=1}^{\rho} \delta_{r_{N-\rho}-1} \cdot (-1)^{r_{N-\rho}+s_{N-\rho}} a_{r_{N-\rho}^{(N-\rho)}, s_{N-\rho}^{(N-\rho)}} \cdot \left(\mathrm{H}_{\rho-3} + \delta_{\rho-2} a_{(3-r_{N-2})^{(N-2)}, (3-s_{N-2})^{(N-2)}}\right)}{\prod_{\rho=2}^{N} \sum_{s_{N-\rho}=1}^{\rho} \delta_{r_{N-\rho}-1} \cdot (-1)^{r_{N-\rho}+s_{N-\rho}} \left(\delta_{\rho-N} a_{r_0, s_0} + \mathrm{H}_{N-1-\rho} a_{r_{N-\rho}^{(N-\rho)}, s_{N-\rho}^{(N-\rho)}}\right) \left(\mathrm{H}_{\rho-3} + \delta_{\rho-2} a_{(3-r_{N-2})^{(N-2)}, (3-s_{N-2})^{(N-2)}}\right)}$$
(6.22)

The functional dependence of the elements $a$ have been suppressed in (6.21) and (6.22) for compactness. In both cases, the number of times an index is primed is expressed as a parenthesized *superscript* on that index. For $\rho = N$-2 for instance, $a_{r_{N-\rho}^{(N-\rho)}, s_{N-\rho}^{(N-\rho)}} = a_{r_2'', s_2''}$. There are actually 4 indices to generalize, $r_{N-\rho}^{(N-\rho)}, s_{N-\rho}^{(N-\rho)}, (3-r_{N-2})^{(N-2)}$, and $(3-s_{N-2})^{(N-2)}$. However, and as in the previous sections, only one of them need be generalized, and the rest can be obtained from it by simple operations. In this section, $r_{N-\rho}^{(N-\rho)}$ is chosen as the one to first generalize. The progression of primed row indices follows that covered in the

previous sections, with

$$r'_1 = r_1 + 1 - H_{r_0-r_1-1}, \qquad (6.23)$$

$$r''_2 = (r_2 + 1 - H_{r_0-r_2-1})H_{r_1-r_2-1} + (r_2 + 2 - H_{r_0-r_2-2})H_{r_2-r_1} = \sum_{u=1}^{2}(r_2 + u - H_{r_0-r_2-u})H_{(-1)^u(r_2-r_1-u+2)}, \qquad (6.24)$$

$$r'''_3 = \sum_{u=1}^{2}(r_3 + u - H_{r_0-r_3-u})H_{(-1)^u(r_3-r_1-u+2)}H_{r_2-r_3-1} + \sum_{u=1}^{2}(r_3 + u + 1 - H_{r_0-r_3-u-1})H_{(-1)^u(r_3-r_1-u+3)}H_{r_3-r_2}$$
$$= \sum_{v=0}^{1}\sum_{u=1}^{2}(r_3 + v + u - H_{r_0-r_3-v-u})H_{(-1)^u(r_3-r_1+v-u+2)}H_{(-1)^v(r_2-r_3+v-1)}, \qquad (6.25)$$

$$r''''_4 = \sum_{v=0}^{1}\sum_{u=1}^{2}(r_4 + u + v - H_{r_0-r_4-u-v})H_{(-1)^u(r_4-r_1-u+v+2)}H_{(-1)^v(r_2-r_4+v-1)}H_{r_3-r_4-1} + \sum_{v=0}^{1}\sum_{u=1}^{2}(r_4 + u + v + 1 - H_{r_0-r_4-u-v-1})H_{(-1)^u(r_4-r_1-u+v+3)}H_{(-1)^v(r_2-r_4+v-2)}H_{r_4-r_3}$$
$$= \sum_{w=0}^{1}\sum_{v=0}^{1}\sum_{u=1}^{2}(r_4 + u + v + w - H_{r_0-r_4-u-v-w})H_{(-1)^u(r_4-r_1-u+v+w+2)}H_{(-1)^v(r_2-r_4+v-w-1)}H_{(-1)^w(r_3-r_4+w-1)}. \qquad (6.26)$$

In general then, it may be deduced that for any row index $r_K$ primed $K$ times,

$$r_K^{(K)} = \prod_{k=1}^{K-1}\sum_{u_k=\delta_{k-1}}^{1+\delta_{k-1}}\left[H_{k-2} + \delta_{k-1}\left(r_K + \sum_{l=1}^{K-1}u_l - H\left(r_0 - r_K - \sum_{l=1}^{K-1}u_l\right)\right)\right]H\left((-1)^{u_k+\delta_{k-1}}\left(r_k - r_K + 2u_k - \sum_{l=k}^{K-1}u_l - \delta_{k-1} - 1\right)\right),$$
$$K \in \{2, 3, \cdots, N-3, N-2\}, \; N \geq 3. \qquad (6.27)$$

Throughout this report, and for compactness reasons, the argument of the Heaviside step-function has been displayed as a subscript, instead of being parenthesized as would be the case for a standard function. However, the arguments have grown relatively large, and for the sake of clarity, they are now displayed parenthesized. In (6.27), and for $K = 4$, the set $\{u_1, u_2, u_3\}$ corresponds bijectively to the set $\{u, v, w\}$ used in (6.26), for instance. The above expression is true for $K \geq 2$, but fails for $K = 1$ (6.23). Therefore, (6.27) must be combined with (6.23), in the event that $K = 1$,

$$r_K^{(K)} = \begin{cases} r_1^{(1)} = r'_1 = r_1 + 1 - H_{r_0-r_1-1}, \; K = 1; \\ \prod_{k=1}^{K-1}\sum_{u_k=\delta_{k-1}}^{1+\delta_{k-1}}\left[H_{k-2} + \delta_{k-1}\left(r_K + \sum_{l=1}^{K-1}u_l - H\left(r_0 - r_K - \sum_{l=1}^{K-1}u_l\right)\right)\right]H\left((-1)^{u_k+\delta_{k-1}}\left(r_k - r_K + 2u_k - \sum_{l=k}^{K-1}u_l - \delta_{k-1} - 1\right)\right) \\ K \in \{2, 3, \cdots, N-3, N-2\}, \; N \geq 3. \end{cases} \qquad (6.28)$$

However, the expression can also be reduced to the more compact form of

$$r_K^{(K)} = \prod_{k=1}^{K-H_{K-2}}\sum_{u_k=\delta_{k-1}}^{\delta_{k-1}+H_{K-2}}\left[H_{k-2} + \delta_{k-1}\left(r_K + \sum_{l=1}^{K-H_{K-2}}u_l - H\left(r_0 - r_K - \sum_{l=1}^{K-H_{K-2}}u_l\right)\right)\right]H\left((-1)^{u_k+\delta_{k-1}}\left(r_k - r_K + 2u_k - \sum_{l=k}^{K-H_{K-2}}u_l - \delta_{k-1} - H_{K-2}\right)\right)$$
$$K \in \{1, 2, 3, \cdots, N-3, N-2\}, \; N \geq 3. \qquad (6.29)$$

The expression is now valid for $K \geq 1$: It reduces to (6.23) for $K = 1$, and to (6.27) otherwise, since the Heaviside step-function $H_{K-2} = H(K-2)$ vanishes for $K = 1$, but is unity otherwise. Ultimately, whether (6.28) or (6.29) is used, is a matter of preference. Eqs. (6.28) and (6.29) also reveal the functional dependence of a row index $r_K$ primed $K$-times, which involves all row indices up to and including the $K$-th row index itself,

$$r_K^{(K)} \sim r_K^{(K)}\left(r_0, r_1, r_2, \cdots, r_{K-1}, r_K\right). \tag{6.30}$$

In order to obtain a similar relation for the column index, the following operation is carried out on (6.29), based on (6.30),

$$s_K^{(K)} = r_K^{(K)} \cdot \prod_{k=0}^{K} \delta_{r_k - s_k}. \tag{6.31}$$

The row index $(3 - r_K)^{(K)}$ can be obtained by convolution from (6.29), as follows,

$$(3 - r_K)^{(K)} = r_K^{(K)} \otimes \delta_{3 - r_K} \tag{6.32}$$

with a dependence of

$$(3 - r_K)^{(K)} \sim r_K^{(K)}\left(r_0, r_1, r_2, \cdots, r_{K-1}, 3 - r_K\right). \tag{6.33}$$

Lastly,

$$(3 - s_K)^{(K)} = (3 - r_K)^{(K)} \cdot \prod_{k=0}^{K} \delta_{r_k - s_k}. \tag{6.34}$$

Using (6.29, 31-34), along with the inverse matrix expression

$$\mathbf{A}^{-1} \sim \alpha_{s_0, r_0} = \frac{(-1)^{r_0 + s_0} \prod_{\rho=2}^{N-1} \sum_{s_{N-\rho}=1}^{\rho} \delta_{r_{N-\rho}-1} \cdot (-1)^{r_{N-\rho} + s_{N-\rho}} a_{r_{N-\rho}^{(N-\rho)}, s_{N-\rho}^{(N-\rho)}}\left(r_{N-\rho-1}, s_{N-\rho-1}\right) \cdot \left(H_{\rho-3} + \delta_{\rho-2} a_{(3-r_{N-\rho})^{(N-\rho)}, (3-s_{N-\rho})^{(N-\rho)}}\left(r_{N-\rho-1}, s_{N-\rho-1}\right)\right)}{\prod_{\rho=2}^{N} \sum_{s_{N-\rho}=1}^{\rho} \delta_{r_{N-\rho}-1} \cdot (-1)^{r_{N-\rho} + s_{N-\rho}} \left(\delta_{\rho-N} a_{r_0, s_0} + H_{N-1-\rho} a_{r_{N-\rho}^{(N-\rho)}, s_{N-\rho}^{(N-\rho)}}\left(r_{N-\rho-1}, s_{N-\rho-1}\right)\right)\left(H_{\rho-3} + \delta_{\rho-2} a_{(3-r_{N-\rho})^{(N-\rho)}, (3-s_{N-\rho})^{(N-\rho)}}\left(r_{N-3}, s_{N-3}\right)\right)} \tag{6.35}$$

with the arguments of its elements fully restored, it should be possible to find an explicit expression for the inverse of a matrix in terms of the elements of the original matrix. An application is now demonstrated for $N = 5$. The desired targets are (5.21), along with (5.22, 23, 26, 27, 30, 32 - 34). Then

$$\mathbf{A}^{-1} \sim \alpha_{s_0, r_0} = \frac{(-1)^{r_0 + s_0} \prod_{\rho=2}^{4} \sum_{s_{5-\rho}=1}^{\rho} \delta_{r_{5-\rho}-1} \cdot (-1)^{r_{5-\rho} + s_{5-\rho}} a_{r_{5-\rho}^{(5-\rho)}, s_{5-\rho}^{(5-\rho)}}\left(r_{4-\rho}, s_{4-\rho}\right) \cdot \left(H_{\rho-3} + \delta_{\rho-2} a_{(3-r_{5-\rho})^{(5-\rho)}, (3-s_{5-\rho})^{(5-\rho)}}\left(r_{4-\rho}, s_{4-\rho}\right)\right)}{\prod_{\rho=2}^{5} \sum_{s_{5-\rho}=1}^{\rho} \delta_{r_{5-\rho}-1} \cdot (-1)^{r_{5-\rho} + s_{5-\rho}} \left(\delta_{\rho-5} a_{r_0, s_0} + H_{4-\rho} a_{r_{5-\rho}^{(5-\rho)}, s_{5-\rho}^{(5-\rho)}}\left(r_{4-\rho}, s_{4-\rho}\right)\right)\left(H_{\rho-3} + \delta_{\rho-2} a_{(3-r_{5-\rho})^{(5-\rho)}, (3-s_{5-\rho})^{(5-\rho)}}\left(r_2, s_2\right)\right)} \tag{6.36}$$

Expanding the product symbol,

$$\mathbf{A}^{-1} \sim \alpha_{s_0, r_0} = \frac{(-1)^{r_0+s_0} \sum_{s_1=1}^{4} \delta_{r_1-1} \cdot (-1)^{r_1+s_1} a_{r_1', s_1'}(r_0, s_0) \cdot \sum_{s_2=1}^{3} \delta_{r_2-1} \cdot (-1)^{r_2+s_2} a_{r_2'', s_2''}(r_1, s_1) \cdot \sum_{s_3=1}^{2} \delta_{r_3-1} \cdot (-1)^{r_3+s_3} \left( a_{r_3''', s_3'''} \cdot a_{(3-r_3)''', (3-s_3)'''} \right)(r_2, s_2)}{\sum_{s_0=1}^{5} \delta_{r_0-1} \cdot (-1)^{r_0+s_0} a_{r_0, s_0} \cdot \sum_{s_1=1}^{4} \delta_{r_1-1} \cdot (-1)^{r_1+s_1} a_{r_1', s_1'}(r_0, s_0) \cdot \sum_{s_2=1}^{3} \delta_{r_2-1} \cdot (-1)^{r_2+s_2} a_{r_2'', s_2''}(r_1, s_1) \cdot \sum_{s_3=1}^{2} \delta_{r_3-1} \cdot (-1)^{r_3+s_3} \left( a_{r_3''', s_3'''} \cdot a_{(3-r_3)''', (3-s_3)'''} \right)(r_2, s_2)}$$

(6.37)

and simplifying further,

$$\mathbf{A}^{-1} \sim \alpha_{s_0, r_0} = \frac{(-1)^{r_0+s_0} \sum_{s_1=1}^{4} \sum_{s_2=1}^{3} \sum_{s_3=1}^{2} \delta_{r_1-1} \delta_{r_2-1} \delta_{r_3-1} \cdot (-1)^{r_1+s_1+r_2+s_2+r_3+s_3} \cdot a_{r_1', s_1'}(r_0, s_0) \cdot a_{r_2'', s_2''}(r_1, s_1) \cdot \left( a_{r_3''', s_3'''} \cdot a_{(3-r_3)''', (3-s_3)'''} \right)(r_2, s_2)}{\sum_{s_0=1}^{5} \sum_{s_1=1}^{4} \sum_{s_2=1}^{3} \sum_{s_3=1}^{2} \delta_{r_0-1} \delta_{r_1-1} \delta_{r_2-1} \delta_{r_3-1} \cdot (-1)^{r_0+s_0+r_1+s_1+r_2+s_2+r_3+s_3} a_{r_0, s_0} \cdot a_{r_1', s_1'}(r_0, s_0) \cdot a_{r_2'', s_2''}(r_1, s_1) \cdot \left( a_{r_3''', s_3'''} \cdot a_{(3-r_3)''', (3-s_3)'''} \right)(r_2, s_2)}$$

(6.38)

which indeed reproduces (5.21), with the bijective assignment of variables to be clarified later. Then for the primed indices, using (6.29)

$$r_1^{(1)} = r_1' = \prod_{k=1}^{1} \sum_{u_k=\delta_{k-1}}^{\delta_{k-1}} \left[ \mathrm{H}_{k-2} + \delta_{k-1}\left( r_1 + \sum_{l=1}^{1} u_l - \mathrm{H}\left( r_0 - r_1 - \sum_{l=1}^{1} u_l \right) \right) \right] \mathrm{H}\left( (-1)^{u_k+\delta_{k-1}} \left( r_k - r_1 + 2u_k - \sum_{l=k}^{1} u_l - \delta_{k-1} \right) \right) = r_1 + 1 - \mathrm{H}_{r_0-r_1-1}$$

(6.39)

$$s_1' = s_1^{(1)} = r_1^{(1)} \cdot \prod_{k=0}^{1} \delta_{r_k - s_k} = s_1 + 1 - \mathrm{H}_{s_0 - s_1 - 1}$$

(6.40)

$$r_2'' = r_2^{(2)} = \prod_{k=1}^{1} \sum_{u_k=\delta_{k-1}}^{1+\delta_{k-1}} \left[ \mathrm{H}_{k-2} + \delta_{k-1}\left( r_2 + \sum_{l=1}^{1} u_l - \mathrm{H}\left( r_0 - r_2 - \sum_{l=1}^{1} u_l \right) \right) \right] \mathrm{H}\left( (-1)^{u_k+\delta_{k-1}} \left( r_k - r_2 + 2u_k - \sum_{l=k}^{1} u_l - \delta_{k-1} - 1 \right) \right) = \sum_{u_1=1}^{2} \left( r_2 + u_1 - \mathrm{H}_{r_0 - r_2 - u_1} \right) \mathrm{H}_{(-1)^{u_1}(r_2 - r_1 - u_1 + 2)},$$

(6.41)

$$s_2'' = s_2^{(2)} = r_2^{(2)} \cdot \prod_{k=0}^{2} \delta_{r_k - s_k} = \sum_{u_1=1}^{2} \left( s_2 + u_1 - \mathrm{H}_{s_0 - s_2 - u_1} \right) \mathrm{H}_{(-1)^{u_1}(s_2 - s_1 - u_1 + 2)},$$

(6.42)

$$r_3^{(3)} = r_3''' = \prod_{k=1}^{2} \sum_{u_k=\delta_{k-1}}^{1+\delta_{k-1}} \left[ \mathrm{H}_{k-2} + \delta_{k-1}\left( r_3 + \sum_{l=1}^{2} u_l - \mathrm{H}\left( r_0 - r_3 - \sum_{l=1}^{2} u_l \right) \right) \right] \mathrm{H}\left( (-1)^{u_k+\delta_{k-1}} \left( r_k - r_3 + 2u_k - \sum_{l=k}^{2} u_l - \delta_{k-1} - 1 \right) \right)$$
$$= \sum_{u_1=1}^{2} \left[ r_3 + u_1 + u_2 - \mathrm{H}_{r_0 - r_3 - u_1 - u_2} \right] \mathrm{H}_{(-1)^{u_1}(r_3 - r_1 - u_1 + u_2 + 2)} \sum_{u_2=0}^{1} \mathrm{H}_{(-1)^{u_2}(r_2 - r_3 + u_2 - 1)} = \sum_{u_2=0}^{1} \sum_{u_1=1}^{2} \left( r_3 + u_1 + u_2 - \mathrm{H}_{r_0 - r_3 - u_1 - u_2} \right) \mathrm{H}_{(-1)^{u_1}(r_3 - r_1 - u_1 + u_2 + 2)} \mathrm{H}_{(-1)^{u_2}(r_2 - r_3 + u_2 - 1)}$$

(6.43)

$$s_3''' = s_3^{(3)} = r_3^{(3)} \cdot \prod_{k=0}^{3} \delta_{r_k - s_k} = \sum_{u_2=0}^{1} \sum_{u_1=1}^{2} \left( s_3 + u_1 + u_2 - \mathrm{H}_{s_0 - s_3 - u_1 - u_2} \right) \mathrm{H}_{(-1)^{u_1}(s_3 - s_1 - u_1 + u_2 + 2)} \mathrm{H}_{(-1)^{u_2}(s_2 - s_3 + u_2 - 1)}$$

(6.44)

Using (6.32),

$$(3-r_3)^{(3)} = (3-r_3)''' = \prod_{k=1}^{2} \sum_{u_k=\delta_{k-1}}^{1+\delta_{k-1}} \left[ H_{k-2} + \delta_{k-1}\left(3 - r_3 + \sum_{l=1}^{2} u_l - H\left(r_0 + r_3 - \sum_{l=1}^{2} u_l + 3\right)\right) \right] H\left((-1)^{u_k + \delta_{k-1}}\left(r_k + r_3 + 2u_k - \sum_{l=k}^{2} u_l - \delta_{k-1} - 4\right)\right)$$

$$= \sum_{u_1=1}^{2} \left[-r_3 + u_1 + u_2 + 3 - H_{r_0+r_3-u_1-u_2-3}\right] H_{(-1)^{u_1}(-r_3-r_1-u_1+u_2+5)} \sum_{u_2=0}^{1} H_{(-1)^{u_2}(r_2+r_3+u_2-4)} = \sum_{u_2=0}^{1} \sum_{u_1=1}^{2} \left(-r_3 + u_1 + u_2 + 3 - H_{r_0+r_3-u_1-u_2-3}\right) H_{(-1)^{u_1}(-r_3-r_1-u_1+u_2+5)} H_{(-1)^{u_2}(r_2+r_3+u_2-4)}$$
(6.45)

but is also obtainable by applying (6.32) to (6.43), and

$$(3-s_3)''' = (3-s_3)^{(3)} = (3-r_3)^{(3)} \cdot \prod_{k=0}^{3} \delta_{r_k-s_k} = \sum_{u_2=0}^{1}\sum_{u_1=1}^{2} \left(-s_3 + u_1 + u_2 + 3 - H_{s_0+s_3-u_1-u_2-3}\right) H_{(-1)^{u_1}(-s_3-s_1-u_1+u_2+5)} H_{(-1)^{u_2}(s_2+s_3+u_2-4)}$$
(6.46)

The arguments of the Heaviside step-functions have been relocated to the subscript of these indicial functions, for compactness. After the bijective assignment of the indices $\{\{r_0,s_0\},\{r_1,s_1\},\{r_2,s_2\},\{r_3,s_3\}\}$ to $\{\{p,q\},\{m,n\},\{k,l\},\{i,j\}\}$, and of the summation indices $\{u_1, u_2\}$ to $\{u, v\}$, (6.38) indeed reproduces (5.21). Furthermore, (6.39-46) reproduce (5.22, 23, 26, 27, 30, 32-34). For (6.43-46) in particular, the 2 summands must be combined into a single term under 2 sequential summations as in (5.30, 32-34), since each summand is dependent on *both* summation indices. This step has been carried out last in each of (6.43, 45).

An alternative expression to (6.29) for the generally primed row-index $r_K^{(K)}$ is now derived. Beginning once again with $K = 1$,

$$r_1' = r_1 + 1 - H_{r_0-r_1-1} = \kappa(r_0, r_1),$$
(6.47)

and continuing with

$$r_2'' = \left(r_2 + 1 - H_{r_0-r_2-1}\right) H_{r_1-r_2-1} + \left(r_2 + 2 - H_{r_0-r_2-2}\right) H_{r_2-r_1} = \sum_{u=1}^{2} \left(r_2 + u - H_{r_0-r_2-u}\right) H_{(-1)^u(r_2-r_1-u+2)},$$
(6.48)

it can be synthesized from (6.47) using a multiplication followed by a discrete convolution, as follows,

$$r_2'' = \sum_{u=0}^{1} \kappa(r_0, r_2 + u) H_{(-1)^u(r_1-r_2+u-1)} = \kappa(r_0, r_1) \cdot \sum_{u=0}^{1} \delta_{r_1-r_2-u} \otimes H_{(-1)^u(r_1+2u-1)}$$
(6.49)

which is a non-associative triple product involving a multiplication and a discrete convolution. The argument $r_1$ is replaced by $r_2 + u$ in the process of multiplication (•), rendering $\kappa$ independent of $r_1$. Subsequently, the Kronecker delta-function is convolved ($\otimes$) with the Heaviside step-function, which is also a function of $r_1$, and results in the replacement of $r_1$ by ($r_1 - r_2 - u$). The convolution is strictly over $r_1$ alone. It is known that by themselves, convolution and multiplication are commutative. However, it would not be valid to exchange multiplication with discrete convolution in this product. It would also be incorrect to circularly shift the operands within the product. It is however possible to carry out the following permutations, without affecting the outcome,

$$\left(\kappa(r_0,r_1) \bullet \sum_{u=0}^{1} \delta_{r_1-r_2-u}\right) \otimes H_{(-1)^u(r_1+2u-1)} = \left(\sum_{u=0}^{1} \delta_{r_1-r_2-u} \bullet \kappa(r_0,r_1)\right) \otimes H_{(-1)^u(r_1+2u-1)} = H_{(-1)^u(r_1+2u-1)} \otimes \left(\sum_{u=0}^{1} \delta_{r_1-r_2-u} \bullet \kappa(r_0,r_1)\right)$$
(6.50)

If the convolution were executed first, it would not be possible to subsequently replace $r_1$ by $r_2 + u$ in $\kappa$. Continuing with the progression from (6.48),

$$r_3''' = \sum_{u=0}^{1} \left(r_3 + u + H_{r_3+u-r_0}\right) H_{(-1)^u(r_1-r_3-1+u)} H_{r_2-(r_3+1)} + \sum_{u=0}^{1} \left(r_3 + 1 + u + H_{r_3+1+u-r_0}\right) H_{(-1)^u(r_1-r_3-2+u)} H_{r_3-r_2}$$
(6.51)

which after combining the 2 terms using a second, $v$-summation, simplifies to

$$r_3''' = \sum_{v=0}^{1} \sum_{u=0}^{1} \left(r_3 + u + v + H_{r_3+u+v-r_0}\right) H_{(-1)^u(r_1-r_3-1+u-v)} H_{(-1)^v(r_2-r_3-1+v)} .$$
(6.52)

It can be decomposed in terms of multiplications and convolutions as

$$r_3''' = \left\{\left[\left(\kappa(r_0,r_1) \bullet \sum_{u=0}^{1} \delta_{r_1-r_2-u}\right) \otimes H_{(-1)^u(r_1+2u-1)}\right] \bullet \sum_{v=0}^{1} \delta_{r_2-r_3-v}\right\} \otimes H_{(-1)^v(r_2+2v-1)}$$
(6.53)

which demonstrates the correct association. The quantity enclosed by the square bracket is clearly $r_2''$ (6.49). The 1st convolution is over $r_1$, whereas the 2nd convolution, over $r_2$. The convolution is always over the 1st variable in the arguments of the Kronecker delta-function and the discrete Heaviside step-function. For an arbitrary positive integer $K$,

$$r_K^{(K)} = \kappa(r_0,r_1) \bullet \prod_{i=1}^{K-1} \sum_{u_i=0}^{1} \delta_{r_i-r_{i+1}-u_i} \otimes H_{(-1)^{u_i}(r_i+2u_i-1)} , \quad K \in \{2,3,\cdots,N-3,N-2\}, N \geq 3$$
(6.54)

with the understanding that a multiplication is always carried out *before* a convolution from left to right, after expanding the product over $i$,

$$r_K^{(K)} = \kappa(r_0,r_1) \bullet \sum_{u_1=0}^{1} \delta_{r_1-r_2-u_1} \otimes H_{(-1)^{u_1}(r_1+2u_1-1)} \bullet \sum_{u_2=0}^{1} \delta_{r_2-r_3-u_2} \otimes H_{(-1)^{u_2}(r_2+2u_2-1)} \cdots \sum_{u_{K-1}=0}^{1} \delta_{r_{K-1}-r_K-u_{K-1}} \otimes H_{(-1)^{u_{K-1}}(r_{K-1}+2u_{K-1}-1)}$$
(6.55)

and making the associations in the following order,

$$r_K^{(K)} = \left\{ \left[ \left[ \left( \kappa(r_0,r_1) \cdot \sum_{u_1=0}^{1} \delta_{r_1-r_2-u_1} \right) \otimes H_{(-1)^{u_1}(r_1+2u_1-1)} \right] \cdot \sum_{u_2=0}^{1} \delta_{r_2-r_3-u_2} \right\} \otimes H_{(-1)^{u_2}(r_2+2u_2-1)} \right] \cdots \sum_{u_{K-1}=0}^{1} \delta_{r_{K-1}-r_K-u_{K-1}} \right\} \otimes H_{(-1)^{u_{K-1}}(r_{K-1}+2u_{K-1}-1)}$$

(6.56)

The integer $K$, which is Latin, should not be confused with the function $\kappa$ (6.47), which is Greek. In order to account for $r_1'$, the following modification is made to (6.54):

$$r_K^{(K)} = \begin{cases} \kappa(r_0, r_1), \ K=1 \\ \kappa(r_0, r_1) \cdot \prod_{i=1}^{K-1} \sum_{u_i=0}^{1} \delta_{r_i-r_{i+1}-u_i} \otimes H_{(-1)^{u_i}(r_i+2u_i-1)}, K \in \{2,3,\cdots,N-3,N-2\}, N \geq 3 \end{cases}$$

(6.57)

or as

$$r_K^{(K)} = \kappa(r_0, r_1) \cdot \left[ \delta_{K-1} + H_{K-2} \prod_{i=1}^{K-1} \sum_{u_i=0}^{1} \delta_{r_i-r_{i+1}-u_i} \otimes H_{(-1)^{u_i}(r_i+2u_i-1)} \right], K \in \{1,2,3,\cdots,N-3,N-2\}, N \geq 3.$$

(6.58)

The rest of the indices can be obtained using (6.31-34) as previously explained.

For a 5 x 5 matrix, for instance, the relevant primed indices are found from (6.22) as $r_{N-\rho}^{(N-\rho)}, s_{N-\rho}^{(N-\rho)}, (3-r_{N-2})^{(N-2)}$, and $(3-s_{N-2})^{(N-2)}$ with $N = 5$, and $\rho$ from 2 to $(N-1) = 4$:

$$r_1^{(1)} = r_1' = r_1 + 1 - H_{r_0-r_1-1} = \kappa(r_0, r_1) \tag{6.59}$$

$$s_1' = s_1^{(1)} = r_1^{(1)} \cdot \prod_{k=0}^{1} \delta_{r_k-s_k} = s_1 + 1 - H_{s_0-s_1-1} \tag{6.60}$$

$$r_2^{(2)} = r_2'' = \kappa(r_0, r_1) \cdot \prod_{i=1}^{1} \sum_{u_i=0}^{1} \delta_{r_i-r_{i+1}-u_i} \otimes H_{(-1)^{u_i}(r_i+2u_i-1)} = \sum_{u_1=0}^{1} \kappa(r_0, r_2+u_1) H_{(-1)^{u_1}(r_1-r_2+u_1-1)} \xrightarrow{u_1=u-1} \sum_{u=1}^{2} \kappa(r_0, r_2+u-1) H_{(-1)^u(r_2-r_1-u+2)}$$

(6.61)

$$s_2^{(2)} = s_2'' = r_2^{(2)} \cdot \prod_{k=0}^{2} \delta_{r_k-s_k} = r_2'' \cdot \prod_{k=0}^{2} \delta_{r_k-s_k} = \sum_{u=1}^{2} \kappa(s_0, s_2+u-1) H_{(-1)^u(s_2-s_1-u+2)} \tag{6.62}$$

$$r_3^{(3)} = r_3''' = \kappa(r_0, r_1) \cdot \prod_{i=1}^{2} \sum_{u_i=0}^{1} \delta_{r_i-r_{i+1}-u_i} \otimes H_{(-1)^{u_i}(r_i+2u_i-1)} = \left\{ \left[ \left( \kappa(r_0, r_1) \cdot \sum_{u_1=0}^{1} \delta_{r_1-r_2-u_1} \right) \otimes H_{(-1)^{u_1}(r_1+2u_1-1)} \right] \cdot \sum_{u_2=0}^{1} \delta_{r_2-r_3-u_2} \right\} \otimes H_{(-1)^{u_2}(r_2+2u_2-1)}$$

$$= \left\{ \sum_{u_1=0}^{1} \kappa(r_0, r_2+u_1) H_{(-1)^{u_1}(r_1-r_2+u_1-1)} \cdot \sum_{u_2=0}^{1} \delta_{r_2-r_3-u_2} \right\} \otimes H_{(-1)^{u_2}(r_2+2u_2-1)} = \sum_{u_2=0}^{1} \sum_{u_1=0}^{1} \kappa(r_0, r_3+u_1+u_2) H_{(-1)^{u_1}(r_1-r_3+u_1-u_2-1)} H_{(-1)^{u_2}(r_2-r_3+u_2-1)}$$

$$\xrightarrow{u_1=u-1, u_2=v} \sum_{v=0}^{1} \sum_{u=1}^{2} \kappa(r_0, r_3+u+v-1) H_{(-1)^u(r_3-r_1-u+v+2)} H_{(-1)^v(r_2-r_3+v-1)}$$

(6.63)

$$s_3^{(3)} = s_3''' = r_3^{(3)} \cdot \prod_{k=0}^{3} \delta_{r_k - s_k} = r_3''' \cdot \prod_{k=0}^{3} \delta_{r_k - s_k} = \sum_{v=0}^{1} \sum_{u=1}^{2} \kappa(s_0, s_3 + u + v - 1) H_{(-1)^u (s_3 - s_1 - u + v + 2)} H_{(-1)^v (s_2 - s_3 + v - 1)}$$
(6.64)

$$(3 - r_3)^{(3)} = (3 - r_3)''' = r_3^{(3)} \otimes \delta_{3 - r_3} = \sum_{v=0}^{1} \sum_{u=1}^{2} \kappa(r_0, 2 - r_3 + u + v) H_{(-1)^u (-r_3 - r_1 - u + v + 5)} H_{(-1)^v (r_2 + r_3 + v - 4)}$$
(6.65)

and

$$(3 - s_3)''' = (3 - s_3)^{(3)} = (3 - r_3)^{(3)} \cdot \prod_{k=0}^{3} \delta_{r_k - s_k} = \sum_{v=0}^{1} \sum_{u=1}^{2} \kappa(r_0, 2 - s_3 + u + v) H_{(-1)^u (-s_3 - s_1 - u + v + 5)} H_{(-1)^v (s_2 + s_3 + v - 4)}.$$
(6.66)

Eqs. (6.59-66) reproduce (5.22, 23, 26, 27, 30, 32-34), after applying (6.59) in each case, and are in agreement with (6.39-46) using the previous method.

## 7. Generalized functions in terms of standard or elementary functions

This report makes extensive use of the Kronecker delta-function, and the *discrete* Heaviside step-function. As stated in previous sections, they are not conventional, standard functions, but they may be considered to be the discrete analogues of the Dirac delta-function $\delta(\cdot)$, and the Heaviside step-function $H(\cdot)$, which are distributions, or generalized functions [39]. Historically, it appears that the Kronecker delta-function was the first and only one to have been introduced with a mathematical foundation. Although informally used throughout the nineteenth century, the Dirac delta-function and the Heaviside step-function were not rigorously understood, till the development of the theory of distributions by L. Schwartz [40, 41] in the middle of the twentieth century.

In one-dimensional Euclidean real space R, $\delta$ is a linear functional from the space of distributions $S'(R)$, on a test function $f$ from the Schwartz Space $S(R)$, that results in the assignment of a specific value to that function [42-45],

$$\left\langle f, \delta_{x - \xi} \right\rangle = \int_R f(x) \delta(x - \xi) dx = f(\xi); \quad x, \xi \in R$$
(7.1)

and where, by a definition adapted from one due to P.A.M. Dirac [46],

$$\delta(x - \xi) = \begin{cases} \infty, & x = \xi \\ 0, & x \neq \xi \end{cases}$$
(7.2)

The Dirac delta-function may also be viewed as a generalized function, which is an extension of the concept of a conventional function [47]. Many functions are capable of identical behavior in the limit, an instance of which is the scaled Gaussian function [39]

$$\delta(x - \xi) = \lim_{\lambda \to \infty} \left( \frac{\lambda}{4\pi} \right)^{1/2} \exp\left( -\frac{\lambda}{4} (x - \xi)^2 \right); \quad \xi, \lambda \in R.$$
(7.3)

The Kronecker delta-function, which is due to L. Kronecker, maps to the set {0, 1}, the Cartesian product of a set of integers with itself, taking the value of zero everywhere, except on the diagonal subset of that product, for which it is identical with unity [48]. It is demonstrated on a test function $g$ as follows

$$g[r]\delta[r-n] = g[n], \quad n \in Z \tag{7.4}$$

with the Kronecker delta-function defined, in various ways, as

$$\delta[r-n] = \delta_{r-n} = \delta_{rn} = \begin{cases} 1, & r = n \\ 0, & r \neq n \end{cases} \tag{7.5}$$

Like the Dirac delta-function, the Kronecker delta-function may also be viewed as a generalized function, an instance of which is the discrete analogue of (7.3), the discrete Gaussian function [49],

$$\delta[r-n] = \lim_{\kappa \to \infty} \exp\left(-\kappa^2[r-n]^2\right); \quad n, \kappa \in Z. \tag{7.6}$$

In one-dimensional Euclidean real space R, H is a linear functional, from the space of distributions $S'(R)$, on a test function $f$ from the Schwartz Space $S(R)$, that maps that function to its integral over a sub-space of the real-number line [50],

$$\langle f, H_{x-\xi} \rangle = \int_R f(x) H(x-\xi) dx = \int_\xi^\infty f(x) dx; \quad x, \xi \in R \tag{7.7}$$

and where, by a definition due to O. Heaviside [51 - 53],

$$H(x-\xi) = \begin{cases} 0, & x < \xi \\ 1/2, & x = \xi \\ 1, & x > \xi \end{cases} \tag{7.8}$$

but it can also be construed as the distributional derivative of the Dirac delta-function (7.2),

$$\langle f', H_{x-\xi} \rangle = \int_R f'(x) H(x-\xi) dx = -\int_R f(x) H'(x-\xi) dx = -f(\xi) \tag{7.9}$$

obtained through integration-by-parts, and according to (7.1), leads to the conclusion that

$$H'(x-\xi) = \delta(x-\xi). \tag{7.10}$$

Many functions are capable of identical behavior in the limit, an instance of which is the scaled complementary error-function [54], which is the integral of the scaled Gaussian (7.3),

$$H(x-\xi) = \frac{1}{2}\lim_{\lambda \to \infty}\operatorname{erfc}(-\lambda(x-\xi)); \quad \xi, \lambda \in \mathrm{R}. \tag{7.11}$$

The *discrete* Heaviside step-function is the discrete analogue of (7.8),

$$H[r-n] = \begin{cases} 0, & r < n \\ 1/2, & r = n \\ 1, & r > n \end{cases} \tag{7.12}$$

In this case, however, since the function's domain is the set of integers, the inequality $r > n$ effectively implies that $r \geq n + 1$, whereas $r < n$, implies that $r \leq n - 1$. Consequently, the function exhibits *ramp*-like behavior over the sub-domain $\{n\text{-}1, n, n\text{+}1\}$. This is not the case for the continuous Heaviside step-function defined over the set of real numbers R, for which its variable $x$ can be infinitesimally smaller or larger than the constant $\xi$ as required by (7.8), which can render the function step-like for all practical purposes. The discrete Heaviside step-function is instead defined as [55]

$$H[r-n] = H_{r-n} = \begin{cases} 0, & r < n \\ 1, & r \geq n \end{cases}, \quad n \in \mathrm{Z} \tag{7.13}$$

which ensures a step-like behavior over the sub-domain $\{n\text{-}1, n, n\text{+}1\}$. Thus (7.13) is the actual definition of the *discrete* Heaviside step-function used in this report.

The connection between the Kronecker delta-function (7.5) and the discrete Heaviside step-function (7.13), can be found from the discrete analogue of (7.10), using the forward difference scheme [56]

$$\delta[r-n] = H'[r-n] = \lim_{\Delta \to 1}\frac{H[r-n+\Delta]-H[r-n]}{\Delta} = H[r-n+1] - H[r-n] \tag{7.14}$$

with a similar result, using the backward difference scheme. However, the Heaviside step-function is also obtained from the Dirac delta-function by integration from (7.10), followed by (7.1), as

$$\int_{-\infty}^{x}\delta(z-\xi)dz = \int_{-\infty}^{\infty}H(x-z)\delta(z-\xi)dz = H(x-\xi) \tag{7.15}$$

whose discrete analogue yields, using (7.5),

$$H[r-n] = \sum_{k=-\infty}^{r}\delta[k-n] \xrightarrow{k=n+r-s} \sum_{s=n}^{\infty}\delta[r-s]. \tag{7.16}$$

The middle equality uses (7.5), which is only non-zero when $(k\text{-}n)$ is allowed to assume a value of $\geq 0$, or equivalently $r \geq n$, in accordance with the requirements of (7.13).

It is not difficult to implement the Kronecker delta-function and the discrete Heaviside step-function in a computer program [57, 58]. Although dedicated function modules can be created for that purpose, it is unnecessary for the smaller matrices, since the implementation may be efficiently carried out with mere logical operations. This is best illustrated with an example. In §**4**, it was found that the general expression for the minor-matrix of a 4 x 4 matrix **B** is given by (4.5), reproduced here as

$$\mathbf{C_B}(m,n) = \mathbf{B}\left(1+\delta_{m-1}, 2+\delta_{m-1}+\delta_{m-2}, 3+\delta_{m-1}+\delta_{m-2}+\delta_{m-3}; 1+\delta_{n-1}, 2+\delta_{n-1}+\delta_{n-2}, 3+\delta_{n-1}+\delta_{n-2}+\delta_{n-3}\right)$$
(7.17)

In Matlab for instance, after defining the 4 x 4 matrix B, and initializing the elements of a 3 x 3 matrix C to zeros, (7.17) is encoded as

$$\begin{aligned} \text{C} = \text{B}([1+(m==1), 2+(m==1)+(m==2), 3+(m==1)+(m==2)+(m==3)], \ldots \\ [1+(n==1), 2+(n==1)+(n==2), 3+(n==1)+(n==2)+(n==3)]) \end{aligned}$$
(7.18)

Alternatively, (7.17) was found to have the more compact expression given by (4.7),

$$\mathbf{C_B}(m,n) = \mathbf{B}(2-\mathrm{H}_{m-2}, 3-\mathrm{H}_{m-3}, 4-\mathrm{H}_{m-4}; 2-\mathrm{H}_{n-2}, 3-\mathrm{H}_{n-3}, 4-\mathrm{H}_{n-4}),$$
(7.19)

which unlike (7.17), is exclusively in terms of Heaviside step-functions. In Matlab then,

$$\text{C=B}([2-(m>=2), 3-(m>=3), 4-(m>=4)], [2-(n>=2), 3-(n>=3), 4-(n>=4)])$$
(7.20)

Although clearly not required for computer programming implementations, a goal of this report is the expression of the inverse of a matrix using standard functions, which has yet to be demonstrated, and which is now the focus of the remainder of this section. The following section only applies to inverse-matrix expressions which are explicitly in terms of the elements of their respective original matrices themselves.

One method of converting the generalized functions to real, standard functions, is to first limit the domain of these functions to the set of values required by the summations or the indices, for a given formula. The validity of the conversion outside the constrained domain is not relevant to this report. Achieving this conversion is possible by utilizing the gamma function for instance, generally defined as Euler's Integral for $\mathrm{Re}(z) > 0$,

$$\Gamma(z) = \int_0^\infty s^{z-1} e^{-s} ds,$$
(7.21)

but for strictly positive integral arguments, which makes it synonymous with the factorial function [59]

$$\Gamma(\eta+1) = 1 \cdot 2 \cdot 3 \cdots (\eta-2) \cdot (\eta-1) \cdot \eta = \eta!, \quad \eta \in \mathbb{Z}^+.$$
(7.22)

The gamma function has the peculiar, useful property that $\Gamma(1) = \Gamma(2) = 1$, whereas $\Gamma(3) = 2$. However it is infinite for all negative integers, as well as for zero, all of which are irrelevant to this report. The conversion is illustrated beginning with the simplest case, that of the inverse of the 2 x 2 matrix (2.3),

$$\mathbf{D}^{-1} \sim \delta_{ji} = \frac{(-1)^{i+j} d_{1+\delta_{i-1},1+\delta_{j-1}}}{\sum_{j=1}^{2} \delta_{i-1} d_{ij} (-1)^{i+j} d_{1+\delta_{i-1},1+\delta_{j-1}}}, \{i,j\},\{j,i\} \in \{\{1,2\} \times \{1,2\}\} \quad (7.23)$$

with the *italicized* matrix element $\delta_{ji}$ used on the LHS, not to be confused with the non-italicized Kronecker delta-function $\delta_{i-1}$ used on the RHS. Using the simple relation

$$\delta_{z-1} = \Gamma(z) - z + 1, \ z \in \{1,2,3\} \quad (7.24)$$

which need be valid only for the values of the indices $i$ and $j$ required by (7.23),

$$\mathbf{D}^{-1} \sim \delta_{ji} = -\frac{(-1)^{i+j} d_{\Gamma(i)-i+2,\Gamma(j)-j+2}}{\sum_{j=1}^{2} (-1)^j d_{1,j} d_{2,\Gamma(j)-j+2}}, \{i,j\} \in \{\{1,2\} \times \{1,2\}\}. \quad (7.25)$$

An alternative to (7.23) may be obtained by applying to (7.23) the relation

$$\sum_{u=1}^{P} (1+\delta_{z-u}) = P + 1 - H_{z-P-1}, \ P \in Z^+ \quad (7.26)$$

with $P = 1$, which yields

$$\mathbf{D}^{-1} \sim \delta_{ji} = -\frac{(-1)^{i+j} d_{2-H_{i-2}, 2-H_{j-2}}}{\sum_{j=1}^{2} (-1)^j d_{1,j} d_{2,2-H_{j-2}}}, \{i,j\} \in \{\{1,2\} \times \{1,2\}\}. \quad (7.27)$$

After a comparison with (7.25), it is concluded that

$$H_{z-2} = z - \Gamma(z), \ z \in \{1,2,3\}. \quad (7.28)$$

The formula for the inverse of a 2 x 2 matrix is actually trivial, and is given by (2.3). Eq. (7.25) however, which is clearly less efficient than (2.3), or maybe even (7.23), merely illustrates the conversion to standard functions using the simplest of examples.

A more complicated case is that of the inverse of the 3 x 3 matrix **C**. Its general minor-matrix was found to be (3.8). Using (7.28), augmented with

$$H_{z-3} = \Gamma(z) - 1, \ z \in \{1,2,3\} \quad (7.29)$$

yields, in terms of the original matrix elements of **C**, and in terms of standard functions,

$$\mathbf{D_C}(k,l) = \begin{bmatrix} c_{\Gamma(k)-k+2,\Gamma(l)-l+2} & c_{\Gamma(k)-k+2, 4-\Gamma(l)} \\ c_{4-\Gamma(k),\Gamma(l)-l+2} & c_{4-\Gamma(k), 4-\Gamma(l)} \end{bmatrix} \sim c_{2i+(i-2)k-(-1)^i \Gamma(k), 2j+(j-2)l-(-1)^j \Gamma(l)}; \quad (7.30)$$

$$\{k,l\} \in \{\{1,2,3\} \times \{1,2,3\}\}, \{i,j\} \in \{\{1,2\} \times \{1,2\}\}$$

which is the general minor-matrix (3.8) for **C**, and valid for any couple $\{k, l\}$. The use of the gamma function should present no difficulty, since it is known *a priori* due to (7.22) that $\Gamma(1) = \Gamma(2) = 1$, whereas $\Gamma(3) = 2$. For the inverse matrix, beginning with (3.26)

$$\mathbf{C}^{-1} \sim \gamma_{lk} = \frac{(-1)^{k+l}\left(c_{1+\delta_{k-1},1+\delta_{l-1}}c_{3-H_{k-3},3-H_{l-3}} - c_{1+\delta_{k-1},3-H_{l-3}}c_{3-H_{k-3},1+\delta_{l-1}}\right)}{\sum_{l=1}^{3}(-1)^{1+l}c_{1,l}\left(c_{2,1+\delta_{l-1}}c_{3,3-H_{l-3}} - c_{2,3-H_{l-3}}c_{3,1+\delta_{l-1}}\right)}, \quad \{l,k\} \in \{\{1,2,3\} \times \{1,2,3\}\}, \tag{7.31}$$

or alternatively with (7.26),

$$\mathbf{C}^{-1} \sim \gamma_{lk} = \frac{(-1)^{k+l}\left(c_{2-H_{k-2},2-H_{l-2}}c_{3-H_{k-3},3-H_{l-3}} - c_{2-H_{k-2},3-H_{l-3}}c_{3-H_{k-3},2-H_{l-2}}\right)}{\sum_{l=1}^{3}(-1)^{1+l}c_{1,l}\left(c_{2,2-H_{l-2}}c_{3,3-H_{l-3}} - c_{2,3-H_{l-3}}c_{3,2-H_{l-2}}\right)} \tag{7.32}$$

and using (7.28) along with (7.29), yields

$$\mathbf{C}^{-1} \sim \gamma_{lk} = \frac{(-1)^{k+l}\left(c_{\Gamma(k)-k+2,\Gamma(l)-l+2}c_{4-\Gamma(k),4-\Gamma(l)} - c_{\Gamma(k)-k+2,4-\Gamma(l)}c_{4-\Gamma(k),\Gamma(l)-l+2}\right)}{\sum_{l=1}^{3}(-1)^{1+l}c_{1,l}\left(c_{2,\Gamma(l)-l+2}c_{3,4-\Gamma(l)} - c_{2,4-\Gamma(l)}c_{3,\Gamma(l)-l+2}\right)} \tag{7.33}$$

or more succinctly as

$$\mathbf{C}^{-1} \sim \gamma_{lk} = \frac{\sum_{j=0}^{1}(-1)^{j+k+l}c_{\Gamma(k)-k+2,(-1)^{j}\Gamma(l)+j(l+2)-l+2} \cdot c_{4-\Gamma(k),4-(-1)^{j}\Gamma(l)-j(l+2)}}{\sum_{l=1}^{3}\sum_{j=0}^{1}(-1)^{j+l}c_{1,l} \cdot c_{2,\,4-(-1)^{j}\Gamma(l)-j(l+2)} \cdot c_{3,(-1)^{j}\Gamma(l)+j(l+2)-l+2}}, \quad \{l,k\} \in \{\{1,2,3\} \times \{1,2,3\}\}. \tag{7.34}$$

It is an analytical formula for any element $\gamma_{lk}$ in the inverse of a 3 × 3 matrix **C**, explicitly in terms of the elements of **C**, and entirely in terms of standard functions, whereas

$$|\mathbf{C}| = \sum_{l=1}^{3}\sum_{j=0}^{1}(-1)^{j+l}c_{1,l} \cdot c_{2,\,4-(-1)^{j}\Gamma(l)-j(l+2)} \cdot c_{3,(-1)^{j}\Gamma(l)+j(l+2)-l+2} \tag{7.35}$$

is the determinant, and is the denominator of (7.34). Beyond a 3 × 3 matrix, it is difficult to make an equivalent, similarly compact statement about the inverse of a matrix, as the indicial functions become relatively large. This is demonstrated by the next example.

It is difficult to extend (7.24), (7.28) and (7.29) to larger square matrices without significantly sacrificing compactness. A new scheme is thus adopted, beginning with

$$\delta_{z-1} = \left((-1)^{\Gamma(z+2)} - (-1)^{\Gamma(z+1)}\right)/2, \quad z \in Z^{+} \tag{7.36}$$

$$\mathrm{H}_{z-p} = \left(1+(-1)^{\Gamma(z-p+3)}\right)\!/2\,;\ \ p\in\{2,3\},\ \ z\in Z^{+} \tag{7.37}$$

which are both valid for all positive integers by contrast to (7.24, 28, 29), although slightly more complicated, and not as compact. Other, equivalent expressions are also possible. For a 4 x 4 matrix **B**, its general minor-matrix was found to be (4.8). Using (7.37) along with

$$\mathrm{H}_{z-4} = \left(1-(-1)^{\Gamma(6-z)}\right)\!/2,\ \ z\in\{1,2,3,4\} \tag{7.38}$$

there results, wholly in terms of standard functions, the general, 3 x 3 minor-matrix of **B**,

$$\mathbf{C_B}(m,n) = \begin{bmatrix} b_{(3-(-1)^{\Gamma(m+1)})/2,(3-(-1)^{\Gamma(n+1)})/2} & b_{(3-(-1)^{\Gamma(m+1)})/2,(5-(-1)^{\Gamma(n)})/2} & b_{(3-(-1)^{\Gamma(m+1)})/2,(7+(-1)^{\Gamma(6-n)})/2} \\ b_{(5-(-1)^{\Gamma(m)})/2,(3-(-1)^{\Gamma(n+1)})/2} & b_{(5-(-1)^{\Gamma(m)})/2,(5-(-1)^{\Gamma(n)})/2} & b_{(5-(-1)^{\Gamma(m)})/2,(7+(-1)^{\Gamma(6-n)})/2} \\ b_{(7+(-1)^{\Gamma(6-m)})/2,(3-(-1)^{\Gamma(n+1)})/2} & b_{(7+(-1)^{\Gamma(6-m)})/2,(5-(-1)^{\Gamma(n)})/2} & b_{(7+(-1)^{\Gamma(6-m)})/2,(7+(-1)^{\Gamma(6-n)})/2} \end{bmatrix},$$

$$\sim b_{\left(2k+1-(3-2\Gamma(k))\cos\left(\pi\Gamma\left[m(3-2\Gamma(k))+7\Gamma(k)-k-5\right]\right)\right)/2,\,\left(2l+1-(3-2\Gamma(l))\cos\left(\pi\Gamma\left[n(3-2\Gamma(l))+7\Gamma(l)-l-5\right]\right)\right)/2},$$

$$\{k,l\}\in\{\{1,2,3\}\times\{1,2,3\}\}\ \&\ \{m,n\}\in\{\{1,2,3,4\}\times\{1,2,3,4\}\}, \tag{7.39}$$

which itself is indexed to $(k,l)$, but is a function of the indices $\{m,n\}$ of **B**. In §4, the inverse of **B** was found to have the form of (4.37), which with the use of (7.37) and (7.38), becomes

$$\mathbf{B}^{-1} \sim \beta_{nm} = \frac{(-1)^{m+n}\displaystyle\sum_{l=1}^{3}\sum_{j=1}^{2}(-1)^{j+l}\,b_{(3-(-1)^{\Gamma(m+1)})/2,\,\kappa(l,n)}\cdot b_{(5-(-1)^{\Gamma(m)})/2,\,\lambda(j,l,n)}\cdot b_{(7+(-1)^{\Gamma(6-m)})/2,\,\mu(j,l,n)}}{\displaystyle\sum_{n=1}^{4}\sum_{l=1}^{3}\sum_{j=1}^{2}(-1)^{1+j+l+n}\,b_{1,n}\cdot b_{2,\kappa(l,n)}\cdot b_{3,\lambda(j,l,n)}\cdot b_{4,\mu(j,l,n)}},$$

$$\{n,m\}\in\{\{1,2,3,4\}\times\{1,2,3,4\}\}, \tag{7.40}$$

with the expression for the determinant of the 4 x 4 matrix **B** found in the denominator. In terms of the gamma function, (4.39-41) are expressed as

$$M - \mathrm{H}_{m-M} = M - \frac{1}{2}\left(1-(-1)^{\Gamma\left[(m-M+8-5\Gamma(M-1))(3-2\Gamma(M-1))\right]+\Gamma(M-1)}\right),\ \ M\in\{2,3,4\}, \tag{7.41}$$

$$\kappa(l,n) = l+1-\mathrm{H}_{n-l-1} = l+1-\frac{1}{2}\left(1-(-1)^{\Gamma\left[(n-l+7-5\Gamma(l))(3-2\Gamma(l))\right]+\Gamma(l)}\right), \tag{7.42}$$

$$\lambda(j,l,n) = \sum_{u=1}^{2}\left(j+u-\mathrm{H}_{n-j-u}\right)\mathrm{H}_{(-1)^{u}(j-l-u+2)} = \frac{1}{4}\sum_{u=0}^{1}\left(1+2(j+u)+(-1)^{j+\Gamma\left[6(j-1)-(2j-3)(n-u+1)\right]}\right)\left(1+(-1)^{u+\Gamma(l-j+2)}\right) \tag{7.43}$$

$$\mu(j,l,n) = \sum_{u=1}^{2}\left(-j+u+3-\mathrm{H}_{n+j-u-3}\right)\mathrm{H}_{(-1)^{u}(-j-l-u+5)} = \frac{1}{4}\sum_{u=0}^{1}\left(7+2(u-j)-(-1)^{j+\Gamma\left[6(2-j)+(2j-3)(n-u+1)\right]}\right)\left(1+(-1)^{u+\Gamma(l+j-1)}\right) \tag{7.44}$$

with *j*, *l*, *m,* and *n* specified in (7.40). Eq. (7.41) is the general expression for the row indices of the elements under the summands in the numerator. The last 2 equations yield 4! distinct expressions for each value of *u*. These equations may also be re-expressed in terms of trigonometric functions using Euler's Identity, i.e. that (-1) = exp(iπ). However, the resultant expressions will be even more cumbersome, as the arguments of the trigonometric functions will themselves be dependent on gamma functions.

The discrete generalized functions (7.5) and (7.13) can also be expressed in terms of standard, orthogonal functions, such as Bessel functions. The inverse of a 2 x 2 matrix (7.23) required the Kronecker delta-function given by (7.24), or alternatively with (7.36). It can be re-expressed as

$$\delta_{z-1} = \frac{1}{2}(1-z)(z-2)J_0(2z_1) - (-1)^z J_0((z-1)z_1), \quad z \in \{1,2,3\} \tag{7.45}$$

which is less compact, and is also precision-dependent, since the 1st zero ($z_1$) of the zero-th order Bessel function $J_0(\cdot)$, is irrational,

$$z_1 = 2.4048255576957727\cdots \tag{7.46}$$

The definition (7.45) is valid over the values of *z* required by (7.24). In general,

$$\delta_{z-n} = \frac{(n-2)^{n+1}}{2}(n-z)(z-2)J_0(2z_1) + (-1)^{n+z} J_0((z-n)z_1); \quad n, z \in \{1,2,3\}. \tag{7.47}$$

Along with (7.24), the inverse of a 3 x 3 matrix (7.31, 32) also required the use of the discrete Heaviside step-functions (7.37). In terms of the zero-th order Bessel function, it can be generally expressed as

$$H_{z-n} = \frac{1}{2}\left[z - J_0((2-2)z_1) + (-1)^{n+z} J_0((z-2)z_1)\right]; \quad n \in \{2,3\}, \quad z \in \{1,2,3\}. \tag{7.48}$$

The case for *n* = 1 is redundant.

The preceding discussion is focused on the use of *standard* functions, *in lieu* of discrete generalized functions. Although the gamma and the zero-th order Bessel functions may be considered to be standard functions, they are not *elementary* functions according to the Risch algorithm [60, 61]. All elementary functions are classified as standard functions, but the reverse is frequently untrue. The erf(c) and gamma functions are not elementary functions, although they are standard functions. An elementary function is a function of a single variable, found by some algebraic composition of a finite number of exponential, hyperbolic, logarithmic, polynomial, rational, and/or trigonometric functions. The composition may potentially entail the respective inverses of some of these functions. Constants (without restrictions), power functions (inclusive of fractional powers or *n*-th roots), and the absolute-value function (for real arguments) are also considered to be elementary functions. An elementary function may also be multi-valued. However, the product of an elementary function with a non-elementary function such as the gamma function, would be non-elementary. It is possible to express the discrete generalized functions used in this and the previous section, in terms of elementary functions, but it would not lead to significantly more compact expressions. This will be illustrated with some examples

Eq. (7.24) or (7.36), can be expressed solely in terms of elementary functions, and also generalized as follows

$$\delta_{z-n} = \frac{(n-2)^{n+1}}{2}(n-z)(z-2)\cos(2\pi/2) + (-1)^{n+z}\cos((z-n)\pi/2) \, ; \, n, z \in \{1,2,3\}$$

(7.49)

The expression is actually less compact than (7.24, 36) for $n = 1$. Although it uses a trigonometric function which is less esoteric than a Bessel function, like (7.47) which involves the irrational constant $z_1$ (7.46), (7.49) involves another irrational constant $\pi/2$. Eq. (7.49) is functionally identical to (7.47), but with the zero-th order Bessel function and its first zero $z_1$, being respectively replaced by the cosine and its first zero $\pi/2$. In this formulation, some even orthogonal polynomials such as those of Hermite and Legendre, which are considered elementary functions, are also feasible. In terms of an even Hermite polynomial $He_2(z) = z^2 - 1$, for instance

$$\delta_{z-n} = -\frac{(n-2)^{n+1}}{2}(n-z)(z-2)He_2(2z_1) - (-1)^{n+z}He_2((z-n)z_1) \, ; \, n, z \in \{1,2,3\}, z_1 = 1$$

(7.50)

The expression is negated to attain the correct functional convexity. Moreover, the lowest order, even polynomial is used to maintain simplicity and compactness for the more explicit version of (7.50).

Similarly, (7.48) can be alternatively expressed in terms of elementary functions as

$$H_{z-n} = \frac{1}{2}\left[z - \cos((2-2)\pi/2) + (-1)^{n+z}\cos(\pi(z-2)/2)\right], \, n \in \{2,3\}, \, z \in \{1,2,3\}. \quad (7.51)$$

Like the Kronecker delta-function in (7.50), it can also be expressed in terms of even orthogonal polynomials. For instance, an expression equivalent to (7.51) is the following:

$$H_{z-n} = \frac{1}{2}\left[z + He_2((2-2)z_1) - (-1)^{n+z}He_2((z-2)z_1)\right]; n \in \{2,3\}, z \in \{1,2,3\}, z_1 = 1.$$

(7.52)

The above equations fail for $z > 3$, unlike (7.36) and (7.37), and would require modifications. However, for the inverse of a 3 x 3 matrix (7.31, 32), these equations need not be valid for $z > 3$. These examples also demonstrate that the use of elementary functions can result in a significant sacrifice of compactness.

## 8. Applications

The curl of a vector field **F** is frequently encountered in science and engineering, and is defined in an arbitrary orthogonal curvilinear coordinate system as [62]

$$\nabla \times \mathbf{F} = \frac{1}{h_1 h_2 h_3} \begin{vmatrix} h_1\mathbf{u}_1 & h_2\mathbf{u}_2 & h_3\mathbf{u}_3 \\ \frac{\partial}{\partial u_1} & \frac{\partial}{\partial u_2} & \frac{\partial}{\partial u_3} \\ h_1 F_1 & h_2 F_2 & h_3 F_3 \end{vmatrix} = \frac{\mathbf{u}_1}{h_2 h_3}\left(\frac{\partial h_3 F_3}{\partial u_2} - \frac{\partial h_2 F_2}{\partial u_3}\right) + \frac{\mathbf{u}_2}{h_1 h_3}\left(\frac{\partial h_1 F_1}{\partial u_3} - \frac{\partial h_3 F_3}{\partial u_1}\right) + \frac{\mathbf{u}_3}{h_1 h_2}\left(\frac{\partial h_2 F_2}{\partial u_1} - \frac{\partial h_1 F_1}{\partial u_2}\right).$$

(8.1)

Invoking (7.35) for the determinant of a 3 x 3 matrix, making the bijective assignment of

$$\{c_{1,l}, c_{2,\,4-(-1)^j\Gamma(l)-j(l+2)}, c_{3,\,2+(-1)^j\Gamma(l)+j(l+2)-l}\} \longleftrightarrow \{h_l\mathbf{u}_l, \partial/\partial u_{4-(-1)^j\Gamma(l)-j(l+2)}, (hF)_{2+(-1)^j\Gamma(l)+j(l+2)-l}\}$$

(8.2)

the RHS of (8.1) may now be alternatively, and compactly expressed as

$$\nabla \times \mathbf{F} = \sum_{l=1}^{3}\sum_{j=0}^{1}(-1)^{j+l}\frac{h_l\mathbf{u}_l}{h_1 h_2 h_3}\frac{\partial (hF)_{2+(-1)^j\Gamma(l)+j(l+2)-l}}{\partial u_{4-(-1)^j\Gamma(l)-j(l+2)}}.$$

(8.3)

It is comparatively smaller in size than the sum of the first 2 terms, on the RHS of (8.1). The numerator index is not independent of the denominator index. It can be shown that they are related by a discrete convolution ($\otimes$),

$$4-(-1)^j\Gamma(l)-j(l+2) = \delta_{1-j}\otimes\left(2+(-1)^j\Gamma(l)+j(l+2)-l\right).$$

(8.4)

Another application may be that of the volume of a parallel-piped [63], composed of a parallelogram for a base defined by the vectors **a** and **b**, and with a height defined by another vector **c**,

$$V = |(\mathbf{a}\times\mathbf{b})\bullet\mathbf{c}| = \begin{vmatrix} a_1 & a_2 & a_3 \\ b_1 & b_2 & b_3 \\ c_1 & c_2 & c_3 \end{vmatrix} = a_1(b_2 c_3 - b_3 c_2) - a_2(b_1 c_3 - b_3 c_1) + a_3(b_1 c_2 - b_2 c_1). \quad (8.5)$$

Using (7.35) once again, but this time with the assignment of

$$\{c_{1,l}, c_{2,\,4-(-1)^j\Gamma(l)-j(l+2)}, c_{3,\,2+(-1)^j\Gamma(l)+j(l+2)-l}\} \longleftrightarrow \{a_l, b_{4-(-1)^j\Gamma(l)-j(l+2)}, c_{2+(-1)^j\Gamma(l)+j(l+2)-l}\} \quad (8.6)$$

yields the alternative compact expression

$$V = \sum_{l=1}^{3}\sum_{j=0}^{1}(-1)^{j+l}a_l\cdot b_{4-(-1)^j\Gamma(l)-j(l+2)}\cdot c_{2+(-1)^j\Gamma(l)+j(l+2)-l}.$$

(8.7)

## 9. Summary and Conclusions

A telescoping method has been developed for the inverse of a square, $N \times N$ matrix $\mathbf{A}$, based on a general expression of its $(N-1) \times (N-1)$ primary minor-matrix. Another, $(N-2) \times (N-2)$ minor-matrix is then extracted from the $(N-1) \times (N-1)$ minor-matrix, and the process is continued until a $2 \times 2$ minor-matrix is attained. This procedure is illustrated mathematically beginning with the primary minor-matrix $\mathbf{M}_{r_1,s_1}^{(N-1)}(r_0,s_0)$, whose determinants are the entries of the numerator of the inverse of $\mathbf{A}$ (6.2), and for $N \geq 3$,

$$\mathbf{M}_{r_1,s_1}^{(N-1)} = \mathbf{A}\left(\bigcup_{\sigma_1=2}^{N}(\sigma_1 - \mathrm{H}_{r_0-\sigma_1}); \bigcup_{\tau_1=1}^{N}(\tau_1 - \mathrm{H}_{s_0-\tau_1})\right), \{r_1,s_1\} \in \{\{1,\cdots,N-1\}\times\{1,\cdots,N-1\}\},$$

$$\mathbf{M}_{r_2,s_2}^{(N-2)} = \mathbf{M}_{r_1,s_1}^{(N-1)}\left(\bigcup_{\sigma_2=2}^{N-1}(\sigma_2 - \mathrm{H}_{r_1-\sigma_2}); \bigcup_{\tau_2=2}^{N-1}(\tau_2 - \mathrm{H}_{s_1-\tau_2})\right), \{r_2,s_2\} \in \{\{1,\cdots,N-2\}\times\{1,\cdots,N-2\}\},$$

$$\vdots$$

$$\mathbf{M}_{r_{N-3},s_{N-3}}^{(3)} = \mathbf{M}_{r_{N-4},s_{N-4}}^{(4)}\left(\bigcup_{\sigma_{N-3}=2}^{4}(\sigma_{N-3} - \mathrm{H}_{r_{N-4}-\sigma_{N-3}}); \bigcup_{\tau_{N-3}=2}^{4}(\tau_{N-3} - \mathrm{H}_{s_{N-4}-\tau_{N-3}})\right), \{r_{N-3},s_{N-3}\} \in \{\{1,2,3\}\times\{1,2,3\}\},$$

$$\mathbf{M}_{r_{N-2},s_{N-2}}^{(2)} = \mathbf{M}_{r_{N-3},s_{N-3}}^{(3)}\left(\bigcup_{\sigma_{N-2}=2}^{3}(\sigma_{N-2} - \mathrm{H}_{r_{N-2}-\sigma_{N-2}}); \bigcup_{\tau_{N-2}=2}^{3}(\tau_{N-2} - \mathrm{H}_{s_{N-2}-\tau_{N-2}})\right), \{r_{N-2},s_{N-2}\} \in \{\{1,2\}\times\{1,2\}\}.$$

(9.1)

The dimensions of the matrix $\mathbf{A}$ is stated as $N \times N = N^2$, and in the above progression, the *superscript* of a minor-matrix $\mathbf{M}$ indicates the square-root of its dimension. The *index subscript* indicates the number of telescoping steps used to attain the corresponding minor-matrix $\mathbf{M}$. For instance, if $\mathbf{A}$ were a $3 \times 3$ matrix, the couple $(r_{N-2},s_{N-2}) = (r_1,s_1)$ in which the indices are both subscripted with a 1 implies that its corresponding $2 \times 2$ minor-matrix $\mathbf{M}_{r_{N-2},s_{N-2}}^{(2)} = \mathbf{M}_{r_1,s_1}^{(2)}$ is obtained after 1 telescoping step-from $\mathbf{A}$. In general, the index subscript is obtained by negating the superscript and adding the square-root of the dimensions of $\mathbf{A}$. Once all the minor-matrices are found by the telescoping process, they are processed in the formula for the inverse of the $N \times N$ matrix $\mathbf{A}$

$$\mathbf{A}^{-1} \sim \alpha_{s_0,r_0} = \frac{(-1)^{r_0+s_0}\prod_{\rho=2}^{N-1}\sum_{s_{N-\rho}=1}^{\rho}\delta_{r_{N-\rho}-1}\cdot(-1)^{r_{N-\rho}+s_{N-\rho}}m_{r_{N-\rho},s_{N-\rho}}^{(\rho)}\cdot\left(\mathrm{H}_{\rho-3}+\delta_{\rho-2}m_{3-r_{N-2},3-s_{N-2}}^{(2)}\right)}{\prod_{\rho=2}^{N}\sum_{s_{N-\rho}=1}^{\rho}\delta_{r_{N-\rho}-1}\cdot(-1)^{r_{N-\rho}+s_{N-\rho}}\left(\delta_{\rho-N}a_{r_0,s_0}+\mathrm{H}_{N-1-\rho}m_{r_{N-\rho},s_{N-\rho}}^{(\rho)}\right)\left(\mathrm{H}_{\rho-3}+\delta_{\rho-2}m_{3-r_{N-2},3-s_{N-2}}^{(2)}\right)},$$

$$\{s_0,r_0\} \in \{\{1,\cdots,N\}\times\{1,\cdots,N\}\}, \ N \geq 3.$$

(9.2)

where the functional dependence of the various matrix elements is as follows:

$$m_{r_{N-\rho},s_{N-\rho}}^{(\rho)} \sim m_{r_{N-\rho},s_{N-\rho}}^{(\rho)}(r_{N-\rho-1},s_{N-\rho-1}),$$

(9.3)

$$m^{(2)}_{3-r_{N-2},3-s_{N-2}} \sim m^{(2)}_{3-r_{N-2},3-s_{N-2}}(r_{N-3},s_{N-3}). \tag{9.4}$$

The procedure, which is termed the surrogate matrix approach, relies on the discrete generalized functions explained in §**7**, specifically the Kronecker delta-function and the *discrete* Heaviside step-function,

$$\delta_{z\text{-}n} = \delta_{zn} = \delta[z-n] = \begin{cases} 1, & z = n \\ 0, & z \neq n \end{cases} \tag{9.5}$$

$$H_{z-n} = H[z-n] = \begin{cases} 0, & z < n \\ 1, & z \geq n \end{cases}. \tag{9.6}$$

This procedure may be specialized for the respective inverses of the 3 x 3 **C**, 4 x 4 **B** and 5 x 5 **A** matrices

$$\mathbf{C}^{-1} \sim \gamma_{lk} = \frac{\sum_{j=1}^{2}(-1)^{1+j+k+l}\left(d_{1,j}d_{2,3-j}\right)(k,l)}{\sum_{l=1}^{3}\sum_{j=1}^{2}(-1)^{j+l}c_{1,l}\cdot\left(d_{1,j}d_{2,3-j}\right)(1,l)}, \tag{9.7}$$

$$\{l,k\} \in \{\{1,2,3\}\times\{1,2,3\}\},$$

$$\mathbf{B}^{-1} \sim \beta_{nm} = \frac{\sum_{l=1}^{3}\sum_{j=1}^{2}(-1)^{j+l+m+n}c_{1,l}(m,n)\cdot\left(d_{1,j}d_{2,3-j}\right)(1,l)}{\sum_{n=1}^{4}\sum_{l=1}^{3}\sum_{j=1}^{2}(-1)^{1+j+l+n}b_{1,n}\cdot c_{1,l}(1,n)\cdot\left(d_{1,j}d_{2,3-j}\right)(1,l)}, \tag{9.8}$$

$$\{n,m\} \in \{\{1,2,3,4\}\times\{1,2,3,4\}\},$$

$$\mathbf{A}^{-1} \sim \alpha_{qp} = \frac{\sum_{n=1}^{4}\sum_{l=1}^{3}\sum_{j=1}^{2}(-1)^{1+j+l+n+p+q}b_{1,n}(p,q)\cdot c_{1,l}(1,n)\cdot\left(d_{1,j}d_{2,3-j}\right)(1,l)}{\sum_{q=1}^{5}\sum_{n=1}^{4}\sum_{l=1}^{3}\sum_{j=1}^{2}(-1)^{j+l+n+q}a_{1,q}\cdot b_{1,n}(1,q)\cdot c_{1,l}(1,n)\cdot\left(d_{1,j}d_{2,3-j}\right)(1,l)} \tag{9.9}$$

$$\{n,m\} \in \{\{1,2,3,4,5\}\times\{1,2,3,4,5\}\}.$$

In §**6**, it was found possible to express the inverse of a $N$ x $N$ matrix **A** as

$$\mathbf{A}^{-1} \sim \alpha_{s_0,r_0} = \frac{(-1)^{r_0+s_0}\prod_{\rho=2}^{N-1}\sum_{s_{N-\rho}=1}^{\rho}\delta_{r_{N-\rho}-1}\cdot(-1)^{r_{N-\rho}+s_{N-\rho}}a_{r_{N-\rho}^{(N-\rho)},s_{N-\rho}^{(N-\rho)}}\cdot\left(H_{\rho-3}+\delta_{\rho-2}a_{(3-r_{N-\rho})^{(N-\rho)},(3-s_{N-\rho})^{(N-\rho)}}\right)}{\prod_{\rho=2}^{N}\sum_{s_{N-\rho}=1}^{\rho}\delta_{r_{N-\rho}-1}\cdot(-1)^{r_{N-\rho}+s_{N-\rho}}\left(\delta_{\rho-N}a_{r_0,s_0}+H_{N-1-\rho}a_{r_{N-\rho}^{(N-\rho)},s_{N-\rho}^{(N-\rho)}}\right)\left(H_{\rho-3}+\delta_{\rho-2}a_{(3-r_{N-\rho})^{(N-\rho)},(3-s_{N-\rho})^{(N-\rho)}}\right)}, \tag{9.10}$$

$$\{s_0,r_0\} \in \{\{1,\cdots,N\}\times\{1,\cdots,N\}\}, \ N \geq 3$$

and explicitly in terms of the elements of the original matrix **A**. The functional dependence of each of the elements in (9.10) is as follows:

$$a_{r_{N-\rho}^{(N-\rho)}, s_{N-\rho}^{(N-\rho)}} \sim a_{r_{N-\rho}^{(N-\rho)}, s_{N-\rho}^{(N-\rho)}}(r_{N-\rho-1}, s_{N-\rho-1}), \tag{9.11}$$

$$a_{(3-r_{N-\rho})^{(N-\rho)}, (3-s_{N-\rho})^{(N-\rho)}} \sim a_{(3-r_{N-\rho})^{(N-\rho)}, (3-s_{N-\rho})^{(N-\rho)}}(r_{N-\rho-1}, s_{N-\rho-1}) \tag{9.12}$$

and with the following general indicial functions, all valid over $K \in \{1, 2, 3, \cdots, N-3, N-2\}$, $N \geq 3$, beginning with the ($K$)-primed row index $r_K$ derived in §**6**,

$$r_K^{(K)} = \prod_{k=1}^{K-H_{K-2}} \sum_{u_k=\delta_{k-1}}^{\delta_{k-1}+H_{K-2}} \left[ H_{k-2} + \delta_{k-1}\left(r_K + \sum_{l=1}^{K-H_{K-2}} u_l - H\left(r_0 - r_K - \sum_{l=1}^{K-H_{K-2}} u_l\right)\right) \right] H\left((-1)^{u_k+\delta_{k-1}}\left(r_k - r_K + 2u_k - \sum_{l=k}^{K-H_{K-2}} u_l - \delta_{k-1} - H_{K-2}\right)\right) \tag{9.13}$$

or its alternative

$$r_K^{(K)} = \kappa(r_0, r_1) \bullet \left[\delta_{K-1} + H_{K-2}\prod_{i=1}^{K-1}\sum_{u_i=0}^{1}\delta_{r_i-r_{i+1}-u_i} \otimes H_{(-1)^{u_i}(r_i+2u_i-1)}\right], \tag{9.14}$$

and along with the other indices,

$$s_K^{(K)} = r_K^{(K)} \cdot \prod_{k=0}^{K} \delta_{r_k-s_k} \tag{9.15}$$

$$(3-r_K)^{(K)} = r_K^{(K)} \otimes \delta_{3-r_K} \tag{9.16}$$

$$(3-s_K)^{(K)} = (3-r_K)^{(K)} \cdot \prod_{k=0}^{K} \delta_{r_k-s_k}. \tag{9.17}$$

Eq. (9.14) demonstrates that the general row index $r_K^{(K)}$ may also be obtained from successive operations on the row index $\kappa$ of the primary minor-matrix $\mathbf{M}_{r_1, s_1}^{(N-1)}(r_0, s_0)$ found in (9.1). Eq. (9.16) requires a discrete convolution. Only one of the indicial functions need be derived, whereas the rest may be obtained from simple operations on the resultant expression (9.15-17). Although the complexity of (9.10) is similar to that of (9.2), insisting on the explicit dependence of the inverse matrix on the original elements of the general $N \times N$ matrix **A** has resulted in complicated expressions for the indicial functions. This approach may be adapted to the inverses of the 2 x 2 **D**, the 3 x 3 **C**, the 4 x 4 **B** and the 5 x 5 **A** matrices, now in terms of the original elements of those matrices,

$$\mathbf{D}^{-1} \sim \delta_{ji} = \frac{(-1)^{i+j} d_{1+\delta_{i-1}, 1+\delta_{j-1}}}{\sum_{j=1}^{2} \delta_{i-1} d_{ij} (-1)^{i+j} d_{1+\delta_{i-1}, 1+\delta_{j-1}}}, \{i,j\}, \{j,i\} \in \{\{1,2\} \times \{1,2\}\}, \tag{9.18}$$

$$\mathbf{C}^{-1} \sim \gamma_{lk} = \frac{\sum_{j=1}^{2} (-1)^{1+j+k+l} c_{2-\mathrm{H}_{k-2}, j+1-\mathrm{H}_{l-j-1}} \cdot c_{3-\mathrm{H}_{k-3}, 4-j-\mathrm{H}_{l-4+j}}}{\sum_{l=1}^{3} \sum_{j=1}^{2} (-1)^{2+j+l} c_{1,l} \cdot c_{2, j+1-\mathrm{H}_{l-j-1}} \cdot c_{3, 4-j-\mathrm{H}_{l-4+j}}}, \{l,k\} \in \{\{1,2,3\} \times \{1,2,3\}\}, \tag{9.19}$$

$$\begin{cases}
\mathbf{B}^{-1} \sim \beta_{nm} = \dfrac{(-1)^{m+n} \sum_{l=1}^{3} \sum_{j=1}^{2} (-1)^{j+l} b_{2-\mathrm{H}_{m-2}, \kappa(l,n)} \cdot b_{3-\mathrm{H}_{m-3}, \lambda(j,l,n)} \cdot b_{4-\mathrm{H}_{m-4}, \mu(j,l,n)}}{\sum_{n=1}^{4} \sum_{l=1}^{3} \sum_{j=1}^{2} (-1)^{1+j+l+n} b_{1,n} \cdot b_{2, \kappa(l,n)} \cdot b_{3, \lambda(j,l,n)} \cdot b_{4, \mu(j,l,n)}}, \{n,m\} \in \{\{1,2,3,4\} \times \{1,2,3,4\}\}, \\[4pt]
\kappa(l,n) = l+1-\mathrm{H}_{n-l-1}, \quad \lambda(j,l,n) = \sum_{u=1}^{2} (j+u-\mathrm{H}_{n-j-u}) \mathrm{H}_{(-1)^u (j-u-l+2)}, \\[4pt]
\mu(j,l,n) = \sum_{u=1}^{2} (-j+u+3-\mathrm{H}_{n+j-u-3}) \mathrm{H}_{(-1)^u (-j-l-u+5)},
\end{cases} \tag{9.20}$$

$$\begin{cases}
\mathbf{A}^{-1} \sim \alpha_{qp} = \dfrac{(-1)^{p+q} \sum_{n=1}^{4} \sum_{l=1}^{3} \sum_{j=1}^{2} (-1)^{1+j+l+n} a_{2-\mathrm{H}_{p-2}, \kappa(n,q)} \cdot a_{3-\mathrm{H}_{p-3}, \lambda(l,n,q)} \cdot a_{4-\mathrm{H}_{p-4}, \mu(j,l,n,q)} \cdot a_{5-\mathrm{H}_{p-5}, \nu(j,l,n,q)}}{\sum_{q=1}^{5} \sum_{n=1}^{4} \sum_{l=1}^{3} \sum_{j=1}^{2} (-1)^{j+l+n+q} a_{1q} \cdot a_{2, \kappa(n,q)} \cdot a_{3, \lambda(l,n,q)} \cdot a_{4, \mu(j,l,n,q)} \cdot a_{5, \nu(j,l,n,q)}}, \\[4pt]
\{q,p\} \in \{\{1,2,3,4,5\} \times \{1,2,3,4,5\}\}, \\[4pt]
\kappa(m,p) = m+1-\mathrm{H}_{p-m-1}, \quad \lambda(k,m,p) = \sum_{u=1}^{2} (k+u-\mathrm{H}_{p-k-u}) \mathrm{H}_{(-1)^u (k-m-u+2)}, \\[4pt]
\mu(i,k,m,p) = \sum_{u=1}^{2} \sum_{v=0}^{1} (i+u+v-\mathrm{H}_{p-i-u-v}) \mathrm{H}_{(-1)^u (i-m-u+v+2)} \mathrm{H}_{(-1)^v (k-i+v-1)}, \\[4pt]
\nu(j,l,n,q) = \sum_{u=1}^{2} \sum_{v=0}^{1} (-j+u+v+3-\mathrm{H}_{q+j-u-v-3}) \mathrm{H}_{(-1)^u (-j-n-u+v+5)} \mathrm{H}_{(-1)^v (j+l+v-4)}.
\end{cases} \tag{9.21}$$

along with their attendant indicial functions, where noted within the brace of each equation. The denominators of (9.18-21) yield the respective determinants of the matrices, expressed explicitly in terms of the elements of these matrices. The number of possible expressions for the indicial function $\mu$, for instance, increases with the matrix dimension $N \geq 3$ as $N!2^{N-3}$.

The indicial functions were derived using the Kronecker delta-function (9.5) and the discrete Heaviside step-function (9.6), which are discrete generalized functions. They may be regarded as discrete analogs of functionals. Unlike a function, whose domain is a set of numbers, the domain of a functional is a set of functions. Regardless, the expression of discrete generalized functions in terms of standard functions was also explored in this report. For instance, in terms of the gamma function $\Gamma(\bullet)$ (7.22), and in §**7**, it was found that

$$\delta_{z-1} = \Gamma(z) - z + 1, \quad z \in \{1,2,3\}, \tag{9.22}$$

$$\mathrm{H}_{z-p} = (-1)^p \left( p - 2 - (p-3)z - \Gamma(z) \right), \quad p \in \{2,3\}, \quad z \in \{1,2,3\}. \tag{9.23}$$

In §**7**, relations (9.22, 23) were used to find the general minor-matrix of a 3 x 3 matrix **C**,

$$\mathbf{D_C}(k,l) = \begin{bmatrix} c_{2-\mathrm{H}_{k-2},2-\mathrm{H}_{l-2}} & c_{2-\mathrm{H}_{k-2},3-\mathrm{H}_{l-3}} \\ c_{3-\mathrm{H}_{k-3},2-\mathrm{H}_{l-2}} & c_{3-\mathrm{H}_{k-3},3-\mathrm{H}_{l-3}} \end{bmatrix} \sim c_{i+1-\mathrm{H}_{k-i-1},\,j+1-\mathrm{H}_{l-j-1}}, \tag{9.24}$$

$$\{i,j\} \in \{\{1,2\} \times \{1,2\}\}, \& \{k,l\} \in \{\{1,2,3\} \times \{1,2,3\}\},$$

entirely in terms of the original elements of **C**. In terms of standard functions,

$$\mathbf{D_C}(k,l) = \begin{bmatrix} c_{\Gamma(k)-k+2,\Gamma(l)-l+2} & c_{\Gamma(k)-k+2,\,4-\Gamma(l)} \\ c_{4-\Gamma(k),\,\Gamma(l)-l+2} & c_{4-\Gamma(k),\,4-\Gamma(l)} \end{bmatrix} \sim c_{2i+(i-2)k-(-1)^i\Gamma(k),\,2j+(j-2)l-(-1)^j\Gamma(l)} \tag{9.25}$$

which involves the gamma function $\Gamma(\bullet)$. It yields any minor-matrix of **C**, given $\{k, l\}$. Applying (9.23) to (9.19) yields a formula for any element in the inverse of **C**,

$$\mathbf{C}^{-1} \sim \gamma_{lk} = \frac{\sum_{j=0}^{1}(-1)^{j+k+l} c_{\Gamma(k)-k+2,(-1)^j\Gamma(l)+j(l+2)-l+2} \cdot c_{4-\Gamma(k),4-(-1)^j\Gamma(l)-j(l+2)}}{\sum_{l=1}^{3}\sum_{j=0}^{1}(-1)^{j+l} c_{1,l} \cdot c_{2,4-(-1)^j\Gamma(l)-j(l+2)} \cdot c_{3,(-1)^j\Gamma(l)+j(l+2)-l+2}}, \quad \{l,k\} \in \{\{1,2,3\} \times \{1,2,3\}\}. \tag{9.26}$$

The denominator yields the determinant of the 3 x 3 matrix **C**,

$$|\mathbf{C}| = \sum_{l=1}^{3}\sum_{j=0}^{1}(-1)^{j+l} c_{1,l} \cdot c_{2,4-(-1)^j\Gamma(l)-j(l+2)} \cdot c_{3,(-1)^j\Gamma(l)+j(l+2)-l+2}. \tag{9.27}$$

The 2nd and 3rd column indices are not independent, and it can be shown that the 2nd index is obtainable from the 3rd index by a discrete convolution ($\otimes$), or vice versa,

$$4 - (-1)^j \Gamma(l) - j(l+2) = \delta_{1-j} \otimes \left( 2 + (-1)^j \Gamma(l) + j(l+2) - l \right) \tag{9.28}$$

which entails a reflection about $j = 0$, followed by a translation by one unit.

Without significantly sacrificing compactness, it was not possible to extend (9.22) and (9.23) to larger square matrices. In their stead, a new formulation was adopted in §**7**,

$$\delta_{z-p} = \left((-1)^{\Gamma(z-p+3)} - (-1)^{\Gamma(z-p+2)}\right)/2 \,; \quad p \in \{1,2\}\,, \tag{9.29}$$

$$\mathrm{H}_{z-p} = \left(1 + (-1)^{\Gamma(z-p+3)}\right)/2 \,; \quad p \in \{2,3\}\,, \tag{9.30}$$

which are both valid for all positive integers, in contrast to (9.22) and (9.23), although slightly more complicated, and less compact. Similar, equivalent expressions are also possible. In §**4**, the general minor-matrix of a 4 x 4 matrix **B** was found to be (4.8),

$$\mathbf{C}_{\mathbf{B}}(m,n) = \begin{bmatrix} b_{2-\mathrm{H}_{m-2},2-\mathrm{H}_{n-2}} & b_{2-\mathrm{H}_{m-2},3-\mathrm{H}_{n-3}} & b_{2-\mathrm{H}_{m-2},4-\mathrm{H}_{n-4}} \\ b_{3-\mathrm{H}_{m-3},2-\mathrm{H}_{n-2}} & b_{3-\mathrm{H}_{m-3},3-\mathrm{H}_{n-3}} & b_{3-\mathrm{H}_{m-3},4-\mathrm{H}_{n-4}} \\ b_{4-\mathrm{H}_{m-4},2-\mathrm{H}_{n-2}} & b_{4-\mathrm{H}_{m-4},3-\mathrm{H}_{n-3}} & b_{4-\mathrm{H}_{m-4},4-\mathrm{H}_{n-4}} \end{bmatrix} \sim b_{k+1-\mathrm{H}_{m-k-1},l+1-\mathrm{H}_{n-l-1}}\,,$$

$$\tag{9.31}$$

$$\{k,l\} \in \{\{1,2,3\}\times\{1,2,3\}\} \ \& \ \{m,n\} \in \{\{1,2,3,4\}\times\{1,2,3,4\}\}.$$

Using (9.30) along with

$$\mathrm{H}_{z-4} = \left(1 - (-1)^{\Gamma(6-z)}\right)/2\,, \quad z \in \{1,2,3,4\} \tag{9.32}$$

there results, wholly in terms of the original elements, and standard functions, the general, 3 x 3 minor-matrix of **B**,

$$\mathbf{C}_{\mathbf{B}}(m,n) = \begin{bmatrix} b_{(3-(-1)^{\Gamma(m+1)})/2,(3-(-1)^{\Gamma(n+1)})/2} & b_{(3-(-1)^{\Gamma(m+1)})/2,(5-(-1)^{\Gamma(n)})/2} & b_{(3-(-1)^{\Gamma(m+1)})/2,(7+(-1)^{\Gamma(6-n)})/2} \\ b_{(5-(-1)^{\Gamma(m)})/2,(3-(-1)^{\Gamma(n+1)})/2} & b_{(5-(-1)^{\Gamma(m)})/2,(5-(-1)^{\Gamma(n)})/2} & b_{(5-(-1)^{\Gamma(m)})/2,(7+(-1)^{\Gamma(6-n)})/2} \\ b_{(7+(-1)^{\Gamma(6-m)})/2,(3-(-1)^{\Gamma(n+1)})/2} & b_{(7+(-1)^{\Gamma(6-m)})/2,(5-(-1)^{\Gamma(n)})/2} & b_{(7+(-1)^{\Gamma(6-m)})/2,(7+(-1)^{\Gamma(6-n)})/2} \end{bmatrix},$$

$$\sim b_{\left(2k+1-(3-2\Gamma(k))\cos\left(\pi\Gamma\left[m(3-2\Gamma(k))+7\Gamma(k)-k-5\right]\right)\right)/2,\,\left(2l+1-(3-2\Gamma(l))\cos\left(\pi\Gamma\left[n(3-2\Gamma(l))+7\Gamma(l)-l-5\right]\right)\right)/2}\,,$$

$$\tag{9.33}$$

and is valid for any couple $\{m, n\}$. In §**4**, the inverse of **B** was found to have the form of (4.37), which with the use of (9.30) and (9.32), becomes

$$\mathbf{B}^{-1} \sim \beta_{nm} = \frac{(-1)^{m+n}\sum_{l=1}^{3}\sum_{j=1}^{2}(-1)^{j+l}b_{(3-(-1)^{\Gamma(m+1)})/2,\,\kappa(l,n)}\cdot b_{(5-(-1)^{\Gamma(m)})/2,\,\lambda(j,l,n)}\cdot b_{(7+(-1)^{\Gamma(6-m)})/2,\,\mu(j,l,n)}}{\sum_{n=1}^{4}\sum_{l=1}^{3}\sum_{j=1}^{2}(-1)^{1+j+l+n}b_{1,n}\cdot b_{2,\,\kappa(l,n)}\cdot b_{3,\,\lambda(j,l,n)}\cdot b_{4,\,\mu(j,l,n)}}\,,$$

$$\{n,m\} \in \{\{1,2,3,4\}\times\{1,2,3,4\}\}.$$

$$\tag{9.34}$$

Beyond a 3 x 3 matrix **C**, it is concluded that the indicial functions can become large and unwieldy. For the inverse of a 4 x 4 matrix **B** (9.34) for instance, it was found that

$$M - \mathrm{H}_{m-M} = \frac{1}{2}\left(2M - 1 + (-1)^{\Gamma\left[(m-M+8-5\Gamma(M-1))(3-2\Gamma(M-1))\right]+\Gamma(M-1)}\right), \ M \in \{2,3,4\}, \quad (9.35)$$

$$\kappa(l,n) = l + 1 - \mathrm{H}_{n-l-1} = l + 1 - \frac{1}{2}\left(1 - (-1)^{\Gamma\left[(n-l-5\Gamma(l)+7)(3-2\Gamma(l))\right]+\Gamma(l)}\right), \quad (9.36)$$

$$\lambda(j,l,n) = \sum_{u=1}^{2}\left(j+u-\mathrm{H}_{n-j-u}\right)\mathrm{H}_{(-1)^u(j-l-u+2)} = \frac{1}{4}\sum_{u=0}^{1}\left(1+2(j+u)+(-1)^{j+\Gamma[6(j-1)-(2j-3)(n-u+1)]}\right)\left(1+(-1)^{u+\Gamma(l-j+2)}\right)$$
$$(9.37)$$

$$\mu(j,l,n) = \sum_{u=1}^{2}\left(-j+u+3-\mathrm{H}_{n+j-u-3}\right)\mathrm{H}_{(-1)^u(-j-l-u+5)} = \frac{1}{4}\sum_{u=0}^{1}\left(7+2(u-j)-(-1)^{j+\Gamma[6(2-j)+(2j-3)(n-u+1)]}\right)\left(1+(-1)^{u+\Gamma(l+j-1)}\right)$$
$$(9.38)$$

Standard alternatives to the gamma function were also explored for $z \in \{1,2,3\}$. For a 3 x 3 matrix for instance, for which the indices $z$ are limited to just $\{1, 2, 3\}$, and in terms of the zero-th order Bessel functions, it was found in §**7** that the Kronecker delta-function (9.22) could be expressed as

$$\delta_{z-1} = \frac{1}{2}(1-z)(z-2)J_0(2z_1) - (-1)^z J_0\left((z-1)z_1\right), \ z \in \{1,2,3\} \quad (9.39)$$

with $z_1$ being the 1st zero of the zero-th order Bessel function. In general, it was found

$$\delta_{z-n} = \frac{(n-2)^{n+1}}{2}(n-z)(z-2)J_0(2z_1) + (-1)^{n+z}J_0\left((z-n)z_1\right); \ n,z \in \{1,2,3\} \quad (9.40)$$

and in terms of elementary functions in accordance with the Risch Algorithm [61],

$$\delta_{z-n} = \frac{(n-2)^{n+1}}{2}(n-z)(z-2)\cos(2\pi/2) + (-1)^{n+z}\cos\left((z-n)\pi/2\right); \ n,z \in \{1,2,3\}. \quad (9.41)$$

The equation is functionally identical to (9.40), but with the zero-th order Bessel function and its first zero $z_1$, being respectively replaced by the cosine and its first zero, $\pi/2$. In this formulation, most even orthogonal polynomials $\Lambda$, such as those of Hermite (*He*), of Laguerre (*L*), and of Jacobi, are feasible. Those of the latter include the Gegenbauer ($C^{(1)}$), the Chebyshev (*T*), and the Legendre (*P*) polynomials. In general,

$$\delta_{z-n} = \frac{1}{\Lambda_2(-z_0 z_1)}\left[\frac{(n-2)^{n+1}}{2}(n-z)(z-2)\Lambda_2\left((2-z_0)z_1\right) + (-1)^{n+z}\Lambda_2\left((z-z_0-n)z_1\right)\right];$$
$$\Lambda_2 \in \{He_2, L_2; C_2^{(1)}, P_2, T_2\}; \ n,z \in \{1,2,3\}; \ z_0, z_1 \in \mathbb{R}$$
$$(9.42)$$

where $\Lambda_2$ is chosen to be of the lowest (2nd-) even order, and $z_0$ is an optional offset that ensures that the extremum of $\Lambda_2$ coincides with $z = 0$, which lies outside the domain $z$ in (9.42). Furthermore, $z_1$ is the 1st positive zero of $\Lambda_2$. The expression is normalized by $\Lambda_2((z-z_0)z_1)$ at $z = 0$ to ensure that it is convex over $\{1, 2, 3\}$, like $J_0$ in (9.40) and the cosine in (9.41). The 2nd-order Laguerre polynomial in particular, would be $\Lambda_2(z) = L_2(z) = [(z-2)^2 - 2]/2$, which when used in (9.42) would require that $z_0 = -2$, and $z_1 = 2^{1/2}$.

Similarly, the discrete Heaviside step-function could be alternatively expressed as

$$H_{z-n} = \frac{1}{2}\left[z - J_0\left((2-2)z_1\right) + (-1)^{n+z} J_0\left((z-2)z_1\right)\right]; \; n \in \{2,3\}, \; z \in \{1,2,3\} \quad (9.43)$$

which is also in terms of standard functions. In terms of elementary functions, however,

$$H_{z-n} = \frac{1}{2}\left[z - \cos\left((2-2)\pi/2\right) + (-1)^{n+z} \cos\left(\pi(z-2)/2\right)\right], \; n \in \{2,3\}, \; z \in \{1,2,3\}. \quad (9.44)$$

The arguments of the 2nd terms in (9.43, 44) were not simplified in order to clarify their origins. Moreover and in general, for an even, orthogonal polynomial as for (9.42),

$$H_{z-n} = \frac{z}{2} - \frac{\Lambda_2\left((2-z_0-2)z_1\right) - (-1)^{n+z}\Lambda_2\left((z-z_0-2)z_1\right)}{2\Lambda_2(-z_0 z_1)}, \; n \in \{2,3\} \quad (9.45)$$

along with the same parameter list as that of (9.42), except that for $n$. Expressions (9.43) and (9.44) are identical in form, due to the fact that both the zero-th order Bessel and the cosine functions are even function. Based on the above observations, it is concluded that the gamma function produces the most compact *general* expressions, especially as the domain $z$ gets larger for the larger matrices. The gamma function seems to have no equivalent elementary function that would qualitatively yield a similar degree of compactness. Eqs. (9.43), (9.44), or (9.45) can be used instead of (9.23) and (9.30), in (9.24) and (9.31), which are entirely in terms of the discrete Heaviside step-function. Based on these observations, it is concluded that the indices of the matrix elements in a determinant, and those in a general minor-matrix, may be treated as generalized discrete functions, since the Heaviside step-function itself is one such function. In **Appendix A**, it is demonstrated using Matlab, that the determinant of a 3 x 3 matrix may be computed using a variety of functions, all of which yield identical results.

The primary goal of this report has been the exploration of an analytical formula for the inverse of an invertible $N$ x $N$ square matrix, explicitly in terms of the elements of that matrix. However, such a formula is not intended to be an alternative for a computationally efficient algorithm such as the LU-decomposition, which is often used as one of many possible expedient algorithms for finding the inverse of a matrix [1]. This report illustrated for the first time, the complications involved in the attainment of a general analytical formula, which may however be useful for smaller square matrices. It has also been shown for the first time that the indicial functions of an inverse matrix expression can treated as generalized discrete functions. The indicial functions become increasingly complicated for larger $N$, but alternative, even simpler formulations, may be possible.

There are two appendices that follow the **References** section, which demonstrate Matlab implementations of some of the concepts reported herein, such as that of the determinant of a 3 x 3 matrix, and the inverse of a 5 x 5 matrix.

# Appendix A

The following Matlab code may be used to find the determinant (3.25),

$$|\mathbf{C}| = -\sum_{l=1}^{3}\sum_{j=2}^{3}(-1)^{j+l}c_{1,l}\cdot c_{2,\,j-H_{l-j}}\cdot c_{3,\,5-j-H_{l-5+j}}$$

of a randomly generated, 3 x 3 (complex) matrix **C**, using one of four methods, demonstrating that the indices of the matrix elements used in the determinant, are treatable as generalized discrete functions. In each method, the indices of the matrix element used in the determinant of **C**, is expressed using different functions. The METHOD-input is not case-sensitive. However, a spelling error would cause the switch statement to default to the otherwise-condition, which generates the determinant using cosine functions.

```
METHOD=input('select method (e.g. Bessel, or cosine, or gamma, or Heaviside): ','s');
z1=2.4048255576957727; %First zero of the zero-th order Bessel function
C=randn(3)+sqrt(-1)*randn(3) %Generates 3 x 3 matrix C of pseudo-random elements
Suml=0; %Initialize the l-summation
for l=1:1:3,
   Sumj=0; %Initialize the j-summation
   for j=0:1:1,
      switch lower(METHOD)
         case 'heaviside'
            %Using the Heaviside step-function (Eq. (3.25) with j=j+2):
            s='using Heaviside step-functions:';
            Sumj=Sumj-[(-1)^(j+2)]*C(2, 2+j-[l>=2+j])*C(3, 3-j-[l>=3-j]);
         case 'gamma'
            s='using gamma functions:';
            %Using the gamma function: Eq. (7.35)
            Sumj=Sumj+((-1)^j)*C(2, 4-(-1)^j*gamma(l)-j*(l+2))*...
                              C(3, (-1)^j*gamma(l)+j*(l+2)-l+2);
         case 'bessel'
            s='using Bessel functions:';
            %with the Heaviside step-functions in terms of Bessel functions (Eq. (7.48))
            Sumj=Sumj+[(-1)^(j+1)]*C(2,0.5*[5+2*j-l-(-1)^(j+l)*bessel(0,(l-2)*z1)])*...
                              C(3,0.5*[7-2*j-l+(-1)^(-j+l)*bessel(0,(l-2)*z1)]);
         otherwise
            % with the Heaviside step-functions in terms of cosine functions (Eq. (7.51))
            s='using cosine functions:';
            Sumj=Sumj+[(-1)^(j+1)]*C(2,0.5*[5+2*j-l-(-1)^(j+l)*cos((l-2)*pi/2)])*...
                              C(3,0.5*[7-2*j-l+(-1)^(-j+l)*cos((l-2)*pi/2)]);
      end
   end
   Suml=Suml+[(-1)^l]*C(1,l)*Sumj;
end
disp(s) %display method selected
determinant=Suml %determinant of C, computed using selected method
Matlabdet=det(C) %determinant of C, computed using Matlab det-function
```

# Appendix B

After re-arranging (5.45), into the form of

$$\mathbf{A}^{-1} \sim \alpha_{qp} = \frac{(-1)^{p+q}\sum_{n=1}^{4}(-1)^n a_{2-\mathrm{H}_{p-2},\kappa(n,q)}\sum_{l=1}^{3}(-1)^l a_{3-\mathrm{H}_{p-3},\lambda(l,n,q)}\cdot\sum_{j=1}^{2}(-1)^{1+j}a_{4-\mathrm{H}_{p-4},\mu(j,l,n,q)}\cdot a_{5-\mathrm{H}_{p-5},\nu(j,l,n,q)}}{\sum_{q=1}^{5}(-1)^q a_{1q}\sum_{n=1}^{4}(-1)^n a_{2,\kappa(n,q)}\sum_{l=1}^{3}(-1)^l a_{3,\lambda(l,n,q)}\sum_{j=1}^{2}(-1)^j a_{4,\mu(j,l,n,q)}\cdot a_{5,\nu(j,l,n,q)}},$$

$$\{q,p\}\in\{\{1,2,3,4,5\}\times\{1,2,3,4,5\}\}$$

it can be implemented using the following Matlab code, which generates the entire inverse $\mathbf{A}^{-1}$, but is only valid for a 5 x 5 (complex) matrix. The code, which may be Matlab-version-dependent, can be quickly copied, pasted, and executed in Matlab's command window. Using this approach, it can be surmised that the inverse of a *N* x *N* matrix will in general require *N* for-loops.

```
A=randn(5)+sqrt(-1)*randn(5) %generate a 5 x 5 matrix of pseudo-random elements
Aadj=zeros(size(A)); %initialize to zeros, the elements of the adjugate matrix of A
Sumq=0; %initialize the q-summation
tic
for p=1:1:5
for q=1:1:5,
  Sumn=0; %initialize the n-summation
   for n=1:1:4
     Suml=0; %initialize the l-summation
     for l=1:1:3
        Sumj=0; %initialize the j-summation
        for j=1:1:2
           mu_r =4-(p>=4);%...Eq.(5.41)
           mu_c =[(j+1-(q>=j+1))*(n>=j+1)+(j+2-(q>=j+2))*(j>=n)]*(l>=j+1)+ ...
                 [(j+2-(q>=j+2))*(n>=j+2)+(j+3-(q>=j+3))*(n<=j+1)]*(j>=l);%...Eq.(5.32)
           nu_r =5-(p>=5);%...Eq.(5.43)
           nu_c =[(4-j-(q>=(4-j)))*(n>=(4-j))+(5-j-(q>=(5-j)))*(j<=(3-n))]*(l>=(4-j))+ ...
                 [(5-j-(q>=(5-j)))*(n>=(5-j))+(6-j-(q>=(6-j)))*(n<=(4-j))]*(l<=(3-j));%...Eq.(5.34)
           Sumj =Sumj+[(-1)^(1+j)]*A(mu_r, mu_c)*A(nu_r,nu_c);
        end
        lambda_r =3-(p>=3);%...Eq.(5.39)
        lambda_c =[l+1-(q>=l+1)]*(n>=l+1)+[l+2-(q>=l+2)]*(l>=n);%...Eq.(5.27)
        Suml =Suml+[(-1)^(l)]*A(lambda_r, lambda_c)*Sumj;
     end
     kappa_r =2-(p>=2);%...Eq.(5.37)
     kappa_c =n+1-(q>=n+1);%...Eq.(5.23)
     Sumn=Sumn+[(-1)^(n)]*A(kappa_r, kappa_c)*Suml;
  end
  if p == 1%compute determinant:
  Sumq=Sumq+[(-1)^(p+q)]*A(p,q)*Sumn;%...denominator of Eq.(5.45)=Eq.(5.46)
  end
   Aadj(q,p)=[(-1)^(p+q)]*Sumn;%....numerator of Eq.(5.45)
end
end
toc
DetA=Sumq %...Eq.(5.46)
MatlabDeterminant=det(A) %compute determinant of A using Matlab's det-function
InvA=Aadj/DetA %...Eq.(5.45)
MatlabInv=inv(A) %compute the inverse of A using Matlab's inv-function
A*InvA % yields the 5 x 5 identity matrix
```